\documentclass[amscd, leqno]{amsart}
\usepackage{amscd}
\usepackage{amsmath}
\usepackage{amssymb}
\usepackage{url}
\usepackage{color}
\definecolor{blue}{rgb}{0.0, 0.0, 1.0}
\definecolor{green}{rgb}{0.0, 0.5, 0.0}

\theoremstyle{plain}
\newtheorem{theorem}{Theorem}[section]
\newtheorem{proposition}[theorem]{Proposition}
\newtheorem{lemma}[theorem]{Lemma}

\newtheorem{corollary}[theorem]{Corollary}
\theoremstyle{definition}
\newtheorem{definition}[theorem]{Definition}

\theoremstyle{remark}
\newtheorem{remark}[theorem]{Remark}
\newtheorem{problem}[theorem]{Problem}

\newenvironment{pf}{\begin{proof}}{\end{proof}}
\makeatletter

\@addtoreset{equation}{section}
\makeatother

\begin{document}

\title{Analytic torsion for log-Enriques surfaces and Borcherds product}

\author{Xianzhe Dai}
\address{
Department of Mathematics, 
University of Californai, Santa Barbara
CA93106,
USA}
\email{dai@math.ucsb.jedu}
\author{Ken-Ichi Yoshikawa}
\address{
Department of Mathematics,
Faculty of Science,
Kyoto University,
Kyoto 606-8502,
Japan}
\email{yosikawa@math.kyoto-u.ac.jp}

\begin{abstract}
We introduce a holomorphic torsion invariant of log-Enriques surfaces of index two with cyclic quotient singularities of type $\frac{1}{4}(1,1)$.
The moduli space of such log-Enriques surfaces with $k$ singular points is a modular variety of orthogonal type 
associated with a unimodular lattice of signature $(2,10-k)$.
We prove that the invariant, viewed as a function on the modular variety, is given by the Petersson norm of an explicit Borcherds product.
We note that this torsion invariant is essentially the BCOV invariant in the complex dimension $2$. As a consequence, the BCOV invariant
in this case is not a birational invariant, unlike the Calabi-Yau case.
\end{abstract}

\maketitle
\tableofcontents

\section
{Introduction}
\label{sec:Introduction}
\par
The analytic torsion, which is a certain
combination of the determinants of Hodge Laplacians on differential forms, is an
invariant of Riemannian manifolds defined by Ray and Singer
\cite{RS} as an analytic analog of the Reidemeister torsion, the first topological invariant which is not a homotopy invariant. It was proved independently by
 Cheeger \cite{C} and M\"uller \cite{Mu} that the analytic torsion and the Reidemeister torsion
agree on closed manifolds (Ray-Singer conjecture). Ray and Singer \cite{RaySinger73} also introduced a version of the analytic torsion for complex manifolds, usually referred as the holomorphic torsion. The holomorphic torsion has found significant applications in Arakelov theory, canonical metrics, and mirror symmetry. Unlike its real analogue, it depends on the geometry and complex structure of the underlying complex manifold \cite{BGS88} (the anomaly formula), which gives rise to interesting functions on moduli spaces. In this paper we focus on this aspect of holomorphic torsion, i.e., its connection with modular forms.

In fact, Ray and Singer already noticed the remarkable connection. Using Kronecker's first limit formula, Ray and Singer \cite{RaySinger73} computed the analytic torsion for elliptic curves and found it to be given in terms of the Jacobi $\Delta$-function, a modular form of weight $12$ on $\mathbf H/{\rm SL}(2, \mathbf Z)$.
Their result has since then been extended to higher genus Riemann surfaces by Zograf \cite{Zograf97}, McIntyre-Takhtajan \cite{McIntyreTakhtajan06},
Kokotov-Korotkin \cite{KK} and McIntyre-Park \cite{McIntyrePark14}; Zograf and McIntyre-Takhtajan studied the analytic torsion of Riemann surfaces 
with respect to the hyperbolic metric, while Kokotov-Korotkin and McIntyre-Park studied it with respect to the (degenerate) flat metric attached to an abelian differential of the Riemann surface.

In dimension two, motivated by string duality, the second author \cite{Yoshikawa04} studied the case of 
$2$-elementary $K3$ surfaces,
i.e., the pairs consisting of a $K3$ surface $X$ and a holomorphic involution $\iota:\, X\longrightarrow X$  (acting nontrivially on holomorphic two forms) 
and introduced a (equivariant) holomorphic torsion invariant for those surfaces.
By the global Torelli theorem for $K3$ surfaces, the moduli space of 2-elementary K3 surfaces of fixed topological type is a modular variety 
of orthogonal type, so the holomorphic torsion invariant is viewed as a function on such modular varieties. 
On orthogonal modular varieties, Borcherds \cite{Borcherds98} constructed a class of automorphic forms
with remarkable properties as singular theta lifts of elliptic modular forms. These automorphic forms are called Borcherds products. 
It is shown that the holomorphic torsion invariant of $2$-elementary $K3$ surfaces is given by the Petersson norm 
of a certain series of Borcherds products \cite{Yoshikawa13}, \cite{MaYoshikawa15}.

If $\iota$ is fixed point free, then the quotient $Y=X/\iota$ is an Enriques surface
, whose holomorphic torsion invariant is given by one of the most remarkable Borcherds products, the Borcherds $\Phi$-function of rank $10$. 
In this paper we extend this result to a class of singular rational surfaces called {\em log-Enriques surfaces} introduced by D.-Q. Zhang \cite{Zhang91}. 
As in the case of Enriques surfaces, a log-Enriques surface $Y$ is expressed as a quotient $Y=X/\iota$, 
where $X$ is a $K3$ surface with rational double points, called the canonical covering of $Y$, and $\iota$ is an anti-symplectic involution on $X$ 
free from fixed points outside the singular points. To be precise, our log-Enriques surfaces are those of index two in the sense of Zhang \cite{Zhang91}.
To obtain a nice moduli space, we restrict ourself to the case where $X$ has only nodes as its singular points. A log-Enriques surface with this property
is called {\em good} in this paper. Then a good log-Enriques surface can admit at most $10$ singular points, any of which is a cyclic quotient singularity 
of type $\frac{1}{4}(1,1)$.
It turns out that the moduli space of good log-Enriques surfaces of $k$ singular points is again a Zariski open subset of a modular variety of orthogonal 
type of dimension $10-k$ attached to a unimodular lattice of signature $(2,10-k)$. Let us write ${\mathcal M}_{k}$ for this modular variety.
When $k=2$, we have ${\mathcal M}_{2}^{\rm odd}$ and ${\mathcal M}_{2}^{\rm even}$, according to the parity of the unimodular lattice 
of signature $(2,2)$. For simplicity, we write ${\mathcal M}_{2}$ for ${\mathcal M}_{2}^{\rm odd}$ and ${\mathcal M}_{2}^{\rm even}$ when there is no possibility of confusion.
For a good log-Enriques surface $Y$ with $k$ singular points, we write ${\varpi}(Y)\in {\mathcal M}_k$ for the isomorphism class of $Y$.
Interestingly enough, ${\mathcal M}_{k}$ can be identified with a Zariski open subset of the K\"ahler moduli of a Del Pezzo surface $V$ of degree 
$\deg V = k$, the modular variety given by $\mathcal{KM}(V)=\Omega_{H(V, {\mathbf Z})}/O^{+}(H(V, {\mathbf Z}))$, where $H(V, {\mathbf Z})$ is 
the total cohomology lattice of $V$, $O^{+}(H(V, {\mathbf Z}))$ is its automorphism group, and $\Omega_{H(V, {\mathbf Z})}$ is the domain of 
type IV attached to $H(V, {\mathbf Z})$. (See Theorem~\ref{thm:moduli:log:Enr}.)

Analogously to the Enriques lattice, the Del Pezzo lattice $H(V, {\mathbf Z})$ admits a reflective modular form $\Phi_{V}$ on 
$\Omega_{H(V, {\mathbf Z})}$ for $O^{+}(H(V, {\mathbf Z}))$ of weight $\deg V+4$, which is nowhere vanishing on the Zariski open subset 
corresponding to ${\mathcal M}_{k}$ and characterizes the Heegner divisor of norm $(-1)$-vectors \cite{Yoshikawa09}. 
In addition, $\Phi_{V}$ is the denominator function of a generalized Kac-Moody algebra 
with explicit Fourier series expansion by Gritsenko and Nikulin \cite{Gritsenko18}, \cite{GritsenkoNikulin18}. 
(See \S 8 for more about $\Phi_{V}$.)

On the other hand, even though they are rational surfaces, every log-Enriques surface $Y$ admits a Ricci flat K\"ahler orbifold metric 
\cite{Kobayashi85}. 
Let $\tau(Y)$ denote the analytic torsion of $Y$ in the sense of X. Ma \cite{Ma05} (suitably normalized by volume, see Section~\ref{inv}, 
especially Theorem~\ref{thm:comparison:torsion:logEnriques} and Theorem~\ref{thm:torsion:logEnriques} for the precise definition).
Then our main result says that $\tau(Y)$ is given by some power of the Petersson norm of the Borcherds product $\Phi_{V}$.

\begin{theorem} 
\label{MainTheorem}
There exists a constant $C_{k}>0$ depending only on $k$ such that
for every good log-Enriques surface $Y$ with $k$ singular points,
$$
\tau(Y) = C_{k} \left\|\Phi_{V}({\varpi}(Y))\right\|^{-1/4},
$$
where $V$ is a Del Pezzo surface of degree $k$.
\end{theorem}

It is important to note that our torsion invariant is essentially the complex $2$-dimensional analogue of the BCOV invariant 
(See \cite{BCOV94}, \cite{FLY08}, \cite{EFM18}, \cite{FuZhang20}).
In higher dimensions, Bershadsky, Cecotti, Ooguri and Vafa \cite{BCOV94} introduced a certain combination of holomorphic torsions, 
called the BCOV torsion, and predicted the mirror symmetry at genus one as an equivalence of the BCOV torsion and certain curve 
counting invariants at genus one. The corresponding holomorphic torsion invariant of Calabi-Yau threefolds, called the BCOV invariant, 
was introduced by Fang, Lu and the second author \cite{FLY08}, who verified some prediction in \cite{BCOV94}. Very recently, 
the BCOV invariant is extended to Calabi-Yau manifolds of arbitrary dimension by Eriksson, Freixas i Montplet and Mourougane \cite{EFM18}, who have established the mirror symmetry at genus one for the Dwork family in arbitrary dimension \cite{EFM19}. 
The notion of BCOV invariant is further extended to certain class of pairs by Y. Zhang \cite{YZhang19}, who, together with L. Fu,
has established the birational invariance of the BCOV invariants \cite{YZhang20}, \cite{FuZhang20}.
According to mirror symmetry, the BCOV invariants correspond to the topological string amplitudes whose modular properties are important features.
In the final section, we will interpret Theorem~\ref{MainTheorem} in terms of the BCOV torsion, so that the BCOV invariant
of good log-Enriques surfaces is expressed as the Borcherds product $\Phi_{V}$, an infinite product of expected type in mirror symmetry.
As log-Enriques surfaces are rational, the BCOV invariant is {\em not} a birational invariant in this case.

We remark that the equivalence of the analytic torsion of Ricci flat Enriques surfaces and the Borcherds $\Phi$-function \cite{Yoshikawa04}
may be viewed as the limiting case $k=0$. Since $\tau(Y)$ is the analytic torsion of a resolution of $Y$ with respect to 
a degenerate Ricci flat metric, our theorem may be viewed as a two-dimensional analogue of the theorems of Kokotov-Korotkin \cite{KK} and 
McIntyre-Park \cite{McIntyrePark14} as mentioned above.
Because of the isomorphism between the complex structure moduli of good log-Enriques surfaces and
the K\"ahler moduli of Del Pezzo surfaces, in view of mirror symmetry at genus one as mentioned above,
it may be  
worth asking if the Fourier coefficients of the elliptic modular form appearing in the infinite product expansion of $\Phi_{V}$ are interpreted 
as some counting invariants of Del Pezzo surfaces. 
We also remark that by Theorem~\ref{MainTheorem} and the recent result of S. Ma \cite{Ma19},
the analytic torsion of good log-Enriques surfaces is obtained from the Borcherds $\Phi$-function of rank $10$ 
by manipulating quasi-pullbacks successively. See Section~8.3 for the details.

Our method of proof, which should have independent interest and which carries out the program proposed in \cite[Question 5.18]{Yoshikawa09} 
for $2$-elementary $K3$ surfaces, is to de-singularize the double covering of $Y$ via the Eguchi-Hanson instanton 
to obtain a $2$-elementary $K3$ surface $(\widetilde{X},\theta)$ and study the limiting behavior of the (equivariant) analytic torsion of 
$(\widetilde{X},\theta)$, as well as other constituents of the invariant $\tau(\widetilde{X}, \theta)$ of 
$(\widetilde{X}, \theta)$, as $\widetilde{X}$ degenerates into the orbifold double covering $X$ of $Y$.
As a result, the ratio $\tau(Y)/\tau(\widetilde{X}, \theta)^{1/2}$ may be viewed as the {\em (equivariant) analytic torsion of the Eguchi-Hanson instanton} 
(cf. Theorem~\ref{lopat2}). In \cite{Bismut97}, Bismut computed the behavior of Quillen metrics when the exceptional divisor is blown down to a smooth point.
In this paper, we study the same type of problem, where the blowing-up of ${\mathbf C}^{2}$ will be replaced by the Eguchi-Hanson instanton.
We remark that Theorem~\ref{MainTheorem} would be proved in the same way as in \cite{Yoshikawa04} by making use of
the fundamental theorems for Quillen metrics such as the curvature formula, anomaly formula, and the embedding formula
\cite{Bismut95}, \cite{BGS88}, \cite{BismutLebeau91}, \cite{MaX00}, whose extension to orbifolds were obtained by X. Ma 
\cite{Ma05}, \cite{MaX19}, if we could understand degenerations of log Enriques surfaces. 
On the other hand, it would be difficult to understand the geometric meaning 
of the ratio $\tau(Y)/\tau(\widetilde{X}, \theta)^{1/2}$ by this method.
In the final section, we will observe that $\tau(Y)/\tau(\widetilde{X}, \theta)^{1/2}$ is the key factor
in the exact comparison formula for the BCOV invariants for certain Calabi-Yau orbifolds.

This paper is organized as follows.

In Section~2, we recall log-Enriques surfaces and study their moduli space.
In Section~3, we recall the notion of analytic torsion and also the holomorphic torsion invariant $\tau(\widetilde{X}, \theta)$ 
for $2$-elementary $K3$ surfaces \cite{Yoshikawa04}. 
In Theorem~\ref{thm:holomorphic:torsion:K3}, we will give an explicit formula for the analytic torsion of a $K3$ surface 
with respect to an arbitrary K\"ahler metric.
In Section~4, we recall the Eguchi-Hanson instanton and construct a family of K\"ahler metrics $\{ \gamma_{\epsilon,\delta} \}$ on $\widetilde{X}$ 
converging to an orbifold metric with uniformly bounded Ricci curvature.
In Section~5, we study the behavior of some constituents of the invariant $\tau(\widetilde{X}, \theta)$ with respect to the metric 
$\gamma_{\epsilon,\delta}$ as $\epsilon\to0$.
In Section~6, we derive some estimates for the heat kernels of $(\widetilde{X}, \gamma_{\epsilon,\delta})$. 
In Section~7, we determine the behavior of (equivariant) analytic torsion of $(\widetilde{X}, \theta)$ with respect to the metric 
$\gamma_{\epsilon,\delta}$ as $\epsilon\to0$ and $\delta\to0$.
In Section~8, we introduce a holomorphic torsion invariant for good log-Enriques surfaces and prove the main theorem.
In Section~9, we study the relation between the invariant $\tau(Y)$ and the BCOV invariant.

\bigskip
\noindent
{\bf Acknowledgements }

The first author is partially supported by  NSF-DMS1611915 and the Simons Foundation. The second author is partially supported by JSPS KAKENHI Grant Numbers 16H03935, 16H06335.
This research was initiated when the authors met in the conference ``Analytic torsion and its applications'' held 
in the department of mathematics of Paris-Sud University in June 2012. The authors thank the organizer of this conference.

\section
{log-Enriques surfaces}
\par

\subsection
{\bf log-Enriques surfaces}
\par
Following D.-Q. Zhang \cite{Zhang91}, \cite{Zhang98}, we recall the notion of log-Enriques surfaces (of index $2$) and its basic properties.

\begin{definition}
An irreducible normal projective complex surface $Y$ is called a {\em log-Enriques surface} if the following conditions are satisfied:
\begin{itemize}
\item[(1)]
$Y$ is {\em singular} and has at most quotient singularities except 
rational double points. 
In particular, $Y$ has the structure of a compact complex orbifold.
\item[(2)]
The irregularity of $Y$ vanishes, i.e., $H^{1}(Y,{\mathcal O}_{Y})=0$.
\item[(3)]
Let $K_{Y}$ be the canonical line bundle of $Y$ in the sense of orbifolds. Then
$$
K_{Y}\not\cong{\mathcal O}_{Y},
\qquad
K_{Y}^{\otimes2}={\mathcal O}_{Y}.
$$
\end{itemize}
\end{definition}

\begin{remark}
For ${\mathfrak p} \in {\rm Sing}\,Y$, there exist a neighborhood $U_{\mathfrak p}$ of ${\mathfrak p}$ in $Y$, a finite group 
$G_{\mathfrak p} \subset {\rm GL}({\mathbf C}^{2})$ and a $G_{\mathfrak p}$-invariant neighborhood $V$ of $0$ in ${\mathbf C}^{2}$ such that
$(U_{\mathfrak p}, {\mathfrak p}) \cong (V/G_{\mathfrak p}, 0)$. Then $K_{Y}|_{U_{\mathfrak p}}$ is defined as $(V \times {\mathbf C})/G_{\mathfrak p}$,
where the $G_{\mathfrak p}$-action is given by $g\cdot(z,\zeta) = (g\cdot z, \det(g)\zeta)$.
\end{remark}

\begin{remark}
Logarithmic Enriques surfaces in this paper are those of {\em index two} in Zhang's papers \cite{Zhang91}, \cite{Zhang98}.
We only deal with log-Enriques surfaces of index two in this paper.
\end{remark}

If a {\em smooth} complex surface satisfies conditions (2), (3), then it is an Enriques surface.
For this reason, we impose that {\em log-Enriques surfaces are singular}. Then a log-Enriques surface is rational \cite[Lemma 3.4]{Zhang91}.
By Zhang \cite[Lemma 3.1]{Zhang91}, every singularity of a log-Enriques surface $Y$ is the quotient of a rational double point by 
${\mathbf Z}/2{\mathbf Z}$ and hence non-Gorenstein.
Indeed, if ${\mathfrak p}\in{\rm Sing}\,Y$, then there exists by (1) an isomorphism of germs of analytic spaces 
$(Y,{\mathfrak p})\cong({\mathbf C}^{2}/G,0)$, where $G\subset{\rm GL}({\mathbf C}^{2})$ is a finite group. 
By (3), the image of the homomorphism $\det\colon G\ni g\to\det g\in{\mathbf C}^{*}$ is $\pm1$.
If $G_{0}:=\ker\det\subset G$, then $G_{0}\subset{\rm SL}({\mathbf C}^{2})$ is a normal subgroup of $G$ of index $2$,
so that $(X,0)=({\mathbf C}^{2}/G_{0},0)$ is a rational double point.
If $p\colon(X,0)\to(Y,0)$ denotes the projection induced by the inclusion of groups $G_{0}\subset G$, then
$p$ induces an isomorphism of germs $(X/(G/G_{0}),0)\to(Y,0)$, where $G/G_{0}\cong\{\pm1\}\cong{\mathbf Z}/2{\mathbf Z}$.
By \cite[Lemma 3.1]{Zhang91}, $(X, 0)$ is a rational double point of type $A_{2n-1}$ for some $n$. 
Since the homomorphism $\det^{2} \colon G \to {\mathbf C}^{*}$ is trivial, $K_{Y}^{\otimes 2}$ is a holomorphic line bundle on $Y$ 
in the ordinary sense.

\subsection
{\bf The canonical double covering}
\par
Let $Y$ be a log-Enriques surface and let $\varXi\in H^{0}(Y,K_{Y}^{\otimes2})\setminus\{0\}$ be a nowhere vanishing bicanonical form on $Y$
in the sense of orbifolds. The canonical double covering of $Y$ is defined as
$$
X:=\{(y,\xi)\in K_{Y};\,\xi\otimes\xi=\varXi\}\subset K_{Y},
$$
which is equipped with the projection $p\colon X\to Y$ induced from the projection $K_{Y}\to Y$.
Then $p\colon X\to Y$ is a double covering, which ramifies only over ${\rm Sing}\,Y$. (Since $K_{Y, {\mathfrak p}} = {\mathbf C}/\pm1$
for ${\mathfrak p} \in {\rm Sing}(Y)$, $p^{-1}({\mathfrak p})$ consists of a single point.)
The {\em canonical involution} $\iota\colon X\to X$ is defined as the non-trivial covering transformation:
$$
\iota(y,\xi)=(y,-\xi).
$$
Since the ramification locus of $p \colon X \to Y$ is ${\rm Sing}\,X$, we have $X^{\iota} = {\rm Sing}\,X$ and that
$\iota$ has no fixed points on $X\setminus {\rm Sing}\,X$.
\par
Let $\pi\colon \widetilde{X}\to X$ be the minimal resolution and let ${\theta}\colon\widetilde{X}\to\widetilde{X}$ be the involution
induced by the canonical involution $\iota$. The involution $\theta$ is also called the canonical involution on $\widetilde{X}$.
We have the following commutative diagram:
\begin{equation}
\label{diagram:1}
\begin{CD}
\widetilde{X} & @>\pi>> X & @>p>> & Y=X/\iota
\\
@V\theta VV & @V\iota VV &\, & @VV{\rm id} V
\\
\widetilde{X} & @>>\pi> X & @>>p> & Y=X/\iota
\end{CD}
\end{equation}
Here the projection $p\colon X\to Y$ ramifies only at ${\rm Sing}\,Y$. In what follows, we denote by $X^{\iota}$ and $\widetilde{X}^{\theta}$
the sets of fixed points of $\iota$ and $\theta$, respectively. Since $\iota$ has no fixed points on $X\setminus {\rm Sing}\,X$, $\theta$ has
no fixed points on $\widetilde{X}\setminus \pi^{-1}({\rm Sing}\,X)$. 
Hence $\widetilde{X}\setminus \pi^{-1}({\rm Sing}\,X) \subset \widetilde{X} \setminus \widetilde{X}^{\theta}$. In other words,
$\widetilde{X}^{\theta} \subset \pi^{-1}({\rm Sing}\,X)$.

\begin{lemma}
In the commutative diagram \eqref{diagram:1}, the following hold:
\begin{itemize}
\item[(1)]
$X$ is a $K3$ surface with rational double points and
$$
X^{\iota}={\rm Sing}\,X=p^{-1}({\rm Sing}\,Y),
\qquad
\iota^{*}|_{H^{0}(X,K_{X})}=-1.
$$
\item[(2)]
$(\widetilde{X},\theta)$ is a $2$-elementary $K3$ surface. Namely, $\theta$ acts non-trivially on holomorphic $2$-forms on $\widetilde{X}$.
Moreover, there exists an integer $k\in\{1,\ldots,10\}$ such that
$$
\widetilde{X}^{\theta}=E_{1}\amalg\ldots\amalg E_{k},
\qquad
E_{i}\cong{\bf P}^{1}.
$$
The pair $(\widetilde{X},\theta)$ is called the $2$-elementary $K3$ surface associated to $Y$.
\end{itemize}
\end{lemma}

\begin{pf}
See \cite[Lemma 3.1, Th.\,3.6]{Zhang91} for (1) and \cite[Lemma 2.1]{Zhang98} for (2). 
\end{pf}

\begin{lemma}
\label{lemma:map:log:Enriques}
Let $Y$, $Y'$ be log Enriques surfaces with canonical double coverings $p' \colon X' \to Y'$ and $p \colon X \to Y$, respectively.
Let $\varphi \colon Y' \to Y$ be a birational holomorphic map. 
Then the following hold:
\begin{itemize}
\item[(1)]
$\varphi^{*}$ induces an isomorphism from $H^{0}(Y, K_{Y}^{\otimes2})$ to $H^{0}(Y', K_{Y'}^{\otimes2})$.
\item[(2)]
$\varphi( {\rm Sing}\,Y') \subset {\rm Sing}\,Y$.
\item[(3)]
$\varphi$ lifts to a holomorphic map $f \colon X' \to X$ of canonical double coverings.
\end{itemize}
\end{lemma}

\begin{pf}
(1) 
Let $\Xi \in H^{0}(Y, K_{Y}^{\otimes2})\setminus\{0\}$ and $\Xi' \in H^{0}(Y', K_{Y'}^{\otimes2})\setminus\{0\}$.
Then $\varphi^{*}\Xi$ is a bicanonical from on $Y'\setminus({\rm Sing}\,Y' \cup \varphi^{-1}({\rm Sing}\,Y))$, 
and $\Xi'$ is nowhere vanishing.
We get $\varphi^{*}\Xi/\Xi' \in {\mathcal O}(Y'\setminus({\rm Sing}\,Y' \cup \varphi^{-1}({\rm Sing}\,Y)))
= {\mathcal O}(Y'\setminus\varphi^{-1}({\rm Sing}\,Y)) = {\mathcal O}(Y\setminus {\rm Sing}\,Y) = {\mathcal O}(Y)={\mathbf C}$,
where the first and the third equalities follow from the normality of $Y'$ and $Y$ and the second equality follows from
the Zariski Main Theorem. Hence $\varphi^{*}\Xi = c\,\Xi'$ with some $c \in {\mathbf C}\setminus\{0\}$, and
$\varphi^{*}$ is an isomorphism.
\par
(2)
Let $o \in {\rm Sing}\,Y'$. Assume $\varphi(o) \in Y \setminus {\rm Sing}\,Y$. 
There exist a neighborhood $U$ of $\varphi(o)$ and a nowhere vanishing canonical form $\eta \in H^{0}(U, K_{Y})$.
We can express $\Xi |_{U} = F\cdot\eta^{\otimes2}$, $F \in {\mathcal O}^{*}(U)$. 
Since $\varphi^{*}\Xi$ and $\varphi^{*}F$ are nowhere vanishing on $\varphi^{-1}(U)$, so is $\varphi^{*}\eta^{\otimes2}$.
Hence $\varphi^{*}\eta$ is nowhere vanishing.
Since any singular point of $Y'$ is non-Gorenstein, we get a contradiction. 
Thus $\varphi(o) \in {\rm Sing}\,Y$.
\par
(3) 
Since $\varphi^{*}\Xi$ is nowhere vanishing on $Y'\setminus\varphi^{-1}({\rm Sing}\,Y)$, $\varphi$ has no critical points
on $Y'\setminus\varphi^{-1}({\rm Sing}\,Y)$. Since the restriction of $\varphi$ to $Y'\setminus\varphi^{-1}({\rm Sing}\,Y)$ is a closed map,
$\varphi \colon Y'\setminus\varphi^{-1}({\rm Sing}\,Y) \to Y \setminus {\rm Sing}\,Y$ is an \'etale covering of degree one, i.e., an isomorphism.
$\varphi$ induces a holomorphic map $f \colon X' \setminus (p')^{-1}\varphi^{-1}({\rm Sing}\,Y) \to X\setminus p^{-1}({\rm Sing}\,Y)$
such that $p\circ f = \varphi \circ p'$. Since $p^{-1}(y)$ consists of a unique point for any $y\in {\rm Sing}\,Y$, $f$ extends to a map 
from $X'$ to $X$ by setting $f(x') := p^{-1}(\varphi(p'(x')))$ for $x'\in (p')^{-1}\varphi^{-1}({\rm Sing}\,Y)$. 
By construction, $p\circ f = \varphi \circ p'$. By this equality and the bijectivity of the map $p\colon{\rm Sing}\,X \to {\rm Sing}\,Y$, 
$f$ is continuous. Since $f$ is holomorphic on a Zariski open subset, $f \colon X'\to X$ is holomorphic by the normality of $X'$.
\end{pf}

\subsection
{\bf The good model of a log-Enriques surface}
\par
The group ${\bf Z}/4{\bf Z}$ acts on ${\bf C}^{2}$ as the multiplication by $i=\sqrt{-1}$, i.e., $i(z_{1},z_{2}):=(iz_{1},iz_{2})$.
We define the cyclic quotient singularity of type $\frac{1}{4}(1,1)$ by
$$
({\bf C}^{2}/\langle i\rangle,0).
$$
Its minimal resolution is the total space of the line bundle ${\mathcal O}_{{\bf P}^{1}}(-4)$:
$$
\varpi\colon({\mathcal O}_{{\bf P}^{1}}(-4),E) \to ({\bf C}^{2}/\langle i\rangle,0),
$$
where the exceptional divisor $E=\varpi^{-1}(0)$ is a $(-4)$-curve, i.e., $E^{2}=-4$.

\begin{definition}
A log-Enriques surface $Y$ is {\em good} if $Y$ has only cyclic quotient singularities of type $\frac{1}{4}(1,1)$.
\end{definition}

\par
Let $Y$ be a log-Enriques surface, $p\colon X\to Y$ be its canonical double covering, and  $\pi\colon\widetilde{X}\to X$ be 
the minimal resolution. Then $X$ and $\widetilde{X}$ are equipped with the canonical involutions $\iota$ and $\theta$, respectively.
Let $E=\pi^{-1}({\rm Sing}\,X)$ be the exceptional divisor of $\pi\colon\widetilde{X}\to X$.
Then $E \supset \widetilde{X}^{\theta} =\amalg_{i=1}^{k}E_{i}$ with $1\leq k\leq10$. 
Since $E_{i}$ is a $(-2)$-curve of $\widetilde{X}$, it is a $(-4)$-curve of $\widetilde{X}/\theta$
and its contraction produces a cyclic quotient singularity of type $\frac{1}{4}(1,1)$.

\begin{definition}
The good model of a log-Enriques surface $Y$, denoted by $Y^{\natural}$, is defined as the contraction of the disjoint union of $(-4)$-curves
$\widetilde{X}^{\theta}$ in $\widetilde{X}/\theta$, where $(\widetilde{X},\theta)$ is the $2$-elementary $K3$ surface
associated to $Y$.
\end{definition}

Another construction of $Y^{\natural}$ from $Y$ is as follows \cite[Th.\,3.6]{Zhang91}, \cite[Lemmas\,1.4 and 2.1]{Zhang98}.
Let $\widetilde{Y}$ be the minimal resolution of $Y$ with exceptional divisor $D \subset \widetilde{Y}$. Let $Y^{\#}$ be the blowing-up of $\widetilde{Y}$
at ${\rm Sing}\,D$. Then the proper transform of $D$ consists of disjoint $(-4)$-curves, say $\widetilde{D}_{1},\ldots, \widetilde{D}_{k}$. 
Then $Y^{\#} \cong \widetilde{X}/\theta$ and $Y^{\natural}$ is obtained from $Y^{\#}$ by contracting the $\widetilde{D}_{i}$'s.
(Notice that $\widetilde{Y}$ and $Y^{\#}$ are not log-Enriques surfaces.) As is verified easily, the composition of the rational map
$Y^{\natural} \dashrightarrow Y^{\#}$ and the blowing-down $Y^{\#} \to Y$ extends to a holomorphic map from $Y^{\natural}$ to $Y$. 
\par
By construction, $Y^{\natural}$ has at most cyclic quotient singularities of type $\frac{1}{4}(1,1)$. 
If $Y$ is a good log-Enriques surface, then $Y=Y^{\natural}$.

\begin{proposition}
\label{prop:good:log:Enr}
Let $Y$ be a log-Enriques surface. 
If there is a birational holomorphic map from a good log-Enriques surface $Y'$ to $Y$, then $Y'\cong Y^{\natural}$.
\end{proposition}

\begin{pf}
Let $X^{\natural}$ (resp. $X'$) be the canonical double covering of $Y^{\natural}$ (resp. $Y'$) and
let $\widetilde{X}^{\natural}$ (resp. $\widetilde{X'}$) be the minimal resolution of $X^{\natural}$ (resp. $X'$).
The birational morphism $Y'\to Y$ induces a birational morphism $\psi\colon(X',\iota') \to (X,\iota)$ 
by Lemma~\ref{lemma:map:log:Enriques} (3),
and this $\psi$ induces an isomorphism $f\colon(\widetilde{X}',\theta')\to(\widetilde{X},\theta)=(\widetilde{X}^{\natural},\theta)$,
by the minimality of $K3$ surfaces. 
Hence $(\widetilde{X}'/\theta',(\widetilde{X}')^{\theta'})\cong(\widetilde{X}^{\natural}/\theta,(\widetilde{X}^{\natural})^{\theta})$.
Since the projection $\widetilde{X}'/\theta'\to Y'$ (resp. $\widetilde{X}^{\natural}/\theta\to Y^{\natural}$) is obtained by contracting 
every component of $(\widetilde{X}')^{\theta'}$ (resp. $(\widetilde{X}^{\natural})^{\theta}$) to a cyclic quotient singularity of type $\frac{1}{4}(1,1)$, 
$f$ induces an isomorphism from $Y'$ to $Y$.
\end{pf}

By Proposition~\ref{prop:good:log:Enr}, every log-Enriques surface has a unique good model. 
By Zhang \cite[Th.\,3.6]{Zhang91}, \cite[Th.4, Cor.\,5, Lemma 2.3]{Zhang98}, one can associate to a log-Enriques surface another log-Enriques surface
with a {\em unique} singular point in the canonical way. So log-Enriques surfaces of this type form another class to be studied.
Because of the uniqueness (up to a scaling) of the Ricci-flat ALE hyperk\"ahler metric on the minimal resolution of $A_{1}$-singularity, 
in this paper, we focus on good log-Enriques surfaces.

In the rest of this section, we study the moduli space of good log-Enriques surfaces.
Throughout this paper, we mean by lattice a free ${\bf Z}$-module of finite rank equipped with a non-degenerate integral symmetric bilinear form.
We often identify a lattice with its Gram matrix.

\subsection
{\bf $2$-elementary $K3$ surfaces and log-Enriques surfaces}
\par
A pair $(Z,\iota)$ is called a $2$-elementary $K3$ surface if $Z$ is a $K3$ surface and if $\iota\colon Z\to Z$ is a holomorphic 
anti-symplectic involution. For a $2$-elementary $K3$ surface $(Z,\iota)$, we define
$$
H^{2}(Z,{\bf Z})^{\pm}=\{l\in H^{2}(Z,{\bf Z});\,\iota^{*}(l)=\pm l\},
$$
which is equipped with the integral bilinear form induced from the intersection pairing on $H^{2}(Z,{\bf Z})$.
Then $H^{2}(Z,{\bf Z})$ is isometric to the $K3$-lattice (cf. \cite{BPV84})
$$
{\Bbb L}_{K3}:={\mathbb U}\oplus{\mathbb U}\oplus{\mathbb U}\oplus{\mathbb E}_{8}(-1)\oplus{\mathbb E}_{8}(-1),
$$
where ${\Bbb U}=\binom{0\,1}{1\,0}$ and ${\Bbb E}_{8}(-1)$ is the negative-definite even unimodular lattice of rank $8$ 
whose Gram matrix is given by the Cartan matrix of type $E_{8}$.
If $r$ denotes the rank of $H^{2}(Z,{\bf Z})^{+}$, then $H^{2}(Z,{\bf Z})^{+}$ (resp. $H^{2}(Z,{\bf Z})^{-}$) has signature $(1,r-1)$ 
(resp. $(2,20-r)$). For a $2$-elementary $K3$ surface $(Z,\iota)$, the topological type of $Z^{\iota}$ is determined by the isometry class 
of the lattice $H^{2}(Z,{\bf Z})^{-}$.
\par
Let $Y$ be a {\em good} log-Enriques surface and let $(\widetilde{X},\theta)$ be the corresponding $2$-elementary $K3$ surface.
Hence $(\widetilde{X}/\theta,\widetilde{X}^{\theta})\to(Y,{\rm Sing}(Y))$ is the minimal resolution of the cyclic quotient singularities 
of type $\frac{1}{4}(1,1)$ of $Y$. We set
$$
k:=\#{\rm Sing}(Y)
$$
and define $\Lambda_{k}$ as  the unimodular lattice of signature $(2,10-k)$. Under the identification with a lattice with its Gram matrix,
we have
$$
\Lambda_{k}
=
\begin{pmatrix}
I_{2}&0
\\
0&-I_{10-k}
\end{pmatrix}
\quad
(k\not=2),
\qquad
\Lambda_{2} = 
\begin{pmatrix}
I_{2}&0
\\
0&-I_{2}
\end{pmatrix}
\hbox{ or }
{\mathbb U}\oplus{\mathbb U}
\quad
(k=2).
$$
According to the parity of $\Lambda_{2}$, we set $\Lambda_{2}^{\rm odd} := I_{2}\oplus - I_{2}$ and 
$\Lambda_{2}^{\rm even} := {\mathbb U}\oplus{\mathbb U}$.
Since $\widetilde{X}^{\theta}$ consists of smooth rational curves, we deduce from Nikulin \cite[Th.\,4.2.2]{Nikulin83} that there is an isometry 
of lattices $\alpha\colon H^{2}(\widetilde{X},{\bf Z}) \cong {\Bbb L}_{K3}$ with
\begin{equation}
\label{eqn:2:marking:1}
\alpha\colon H^{2}(\widetilde{X},{\bf Z})^{-} \cong \Lambda_{k}(2).
\end{equation}
Here $\Lambda_{k}(2)$ stands for the rescaling of $\Lambda_{k}$, whose bilinear form is the double of that of $\Lambda_{k}$.
An isometry of lattices $\alpha\colon H^{2}(\widetilde{X},{\bf Z})\cong{\Bbb L}_{K3}$ satisfying \eqref{eqn:2:marking:1} is called a marking 
of $(\widetilde{X},\theta)$.
We set
$$
M_{k}:=\Lambda_{k}(2)^{\perp},
$$
where the orthogonal complement is considered in the $K3$-lattice ${\Bbb L}_{K3}$. A $2$-elementary $K3$ surface is of type $M_{k}$
if its invariant lattice is isometric to $M_{k}$.
\par
We define
$$
\Omega_{k}:=\{[\eta]\in{\bf P}(\Lambda_{k}\otimes{\bf C});\,\langle\eta,\eta\rangle=0,\,\langle\eta,\overline{\eta}\rangle>0\}.
$$
Then $\Omega_{k}$ consists of two connected components $\Omega_{k}^{+}$ and $\Omega_{k}^{-}$, each of which is isomorphic to
bounded symmetric domain of type IV of dimension $10-k$.
Let $O(\Lambda_{k})$ be the automorphism group of $\Lambda_{k}$ and
let $O^{+}(\Lambda_{k})\subset O(\Lambda_{k})$ be the subgroup of index $2$ consisting of elements preserving $\Omega_{k}^{\pm}$.
We define the orthogonal modular variety associated with $\Lambda_{k}$ by
$$
{\mathcal M}_{k}:=\Omega_{k}/O(\Lambda_{k})=\Omega_{k}^{+}/O^{+}(\Lambda_{k}).
$$
When $k=2$, we define ${\mathcal M}_{2}^{\rm odd} := \Omega_{2}/O(\Lambda_{2}^{\rm odd})$ and 
${\mathcal M}_{2}^{\rm even} := \Omega_{2}/O(\Lambda_{2}^{\rm even})$. When there is no possibility of confusion, we write
${\mathcal M}_{2}$ for ${\mathcal M}_{2}^{\rm odd}$ and ${\mathcal M}_{2}^{\rm even}$.
\par
Since $\theta$ acts non-trivially on $H^{0}(\widetilde{X},\Omega^{2}_{\widetilde{X}})$, we deduce the inclusion from the Hodge decomposition
$H^{0}(\widetilde{X},\Omega^{2}_{\widetilde{X}})\subset H^{2}(\widetilde{X},{\bf C})^{-}$.
Since $H^{0}(\widetilde{X},\Omega^{2}_{\widetilde{X}})$ is a complex line,
it follows from the Riemann-Hodge bilinear relations that
$$
\varpi(\widetilde{X},\theta,\alpha)
:=
[\alpha(H^{0}(\widetilde{X},\Omega^{2}_{\widetilde{X}})]\in\Omega_{k}.
$$
The point $\varpi(\widetilde{X},\theta,\alpha)\in\Omega_{k}$ is called the period of $(\widetilde{X},\theta,\alpha)$.
We define the period of $(\widetilde{X},\theta)$ as the $O(\Lambda_{k})$-orbit of $\varpi(\widetilde{X},\theta,\alpha)$, i.e.,
$$
\overline{\varpi}(\widetilde{X},\theta):=O(\Lambda_{k})\cdot[\alpha(H^{0}(\widetilde{X},\Omega^{2}_{\widetilde{X}})]\in{\mathcal M}_{k}.
$$
By \cite[Th.\,1.8]{Yoshikawa04}, the coarse moduli space of $2$-elementary $K3$ surfaces of type $M_{k}$ is isomorphic via the period map
to the analytic space ${\mathcal M}_{k}^{o}:={\mathcal M}_{k}\setminus{\mathcal D}_{k}$,
where ${\mathcal D}_{k}$ is the discriminant divisor
$$
{\mathcal D}_{k}=(\bigcup_{d\in\Lambda_{k},\,d^{2}=-1}d^{\perp})/O(\Lambda_{k}),
\qquad
d^{\perp}:=\{[\eta]\in\Omega_{k};\,\langle\eta,d\rangle=0\}.
$$

\subsection
{\bf The period mapping for log-Enriques surfaces}
\par

\begin{definition}
The period of a good log-Enriques surface $Y$ with $k$ singular points is defined as the period of the corresponding $2$-elementary 
$K3$ surface $(\widetilde{X},\theta)$:
$$
\overline{\varpi}(Y):=\overline{\varpi}(\widetilde{X},\theta)\in{\mathcal M}_{k}.
$$
When $k=2$, we define the parity of $Y$ as that of the lattice $\Lambda_{2}$ defined by \eqref{eqn:2:marking:1}.
\end{definition}

\begin{theorem}
\label{thm:moduli:log:Enr}
The period mapping induces a bijection between the isomorphism classes of good log-Enriques surfaces with $k$ singular points 
(and fixed parity when $k=2$) and ${\mathcal M}_{k}^{o}$. 
\end{theorem}

\begin{pf}
Let ${\mathcal N}_{k}$ be the isomorphism classes of good log-Enriques surfaces with $k$ singular points (and fixed parity when $k=2$).
By \cite[Th.\,1.8]{Yoshikawa04}, 
we can identify ${\mathcal M}_{k}^{o}$ with the isomorphism classes of $2$-elementary $K3$ surfaces of type $M_{k}$ via  the period mapping.
We define a map $f\colon{\mathcal N}_{k}\to{\mathcal M}_{k}^{o}$ by setting $f(Y)=(\widetilde{X},\theta)$,
where $(\widetilde{X},\theta)$ is the $2$-elementary $K3$ surface associated to $Y$.
Similarly, we define a map $g \colon {\mathcal M}_{k}^{o} \to {\mathcal N}_{k}$ by sending $(Z,\sigma) \in {\mathcal M}_{k}^{o}$ to the surface obtained
from $Z/\sigma$ by blowing down $Z^{\sigma}$. 
Since $Z^{\sigma}$ consists of $k$ disjoint $(-2)$-curves, its image in $Z/\sigma$ consists of $k$ disjoint $(-4)$-curves, so that $g(Z, \sigma)$ 
is a good log-Enriques surface with $k$ singular points. 
Since $g=f^{-1}$ by \cite[Lemmas 1.4 and 2.1]{Zhang98}, $f$ is a bijection. 
\end{pf}

Since the (locally defined) family of $2$-elementary $K3$ surfaces of type $M_{k}$ associated to a holomorphic family of good log-Enriques surfaces 
with $k$-singular points is again holomorphic, 
the period mapping for any holomorphic family of good log-Enriques surfaces with $k$-singular points is holomorphic.
In what follows, we regard ${\mathcal M}_{k}^{o}$ as a coarse moduli space of good log-Enriques surfaces with $k$ singular points
(and fixed parity when $k=2$).

\section
{Analytic torsion for $K3$ surfaces and $2$-elementary $K3$ surfaces}
\par

\subsection
{Analytic torsion}
\par
Let $Z$ be a compact complex orbifold of dimension $n$ and let $\gamma$ be a K\"ahler form on $Z$ in the sense of orbifolds.
Let $\iota\colon Z\to Z$ be a holomorphic involution and assume that $\iota$ preserves $\gamma$.
Let $A^{0,q}_{Z}$ be the space of smooth $(0,q)$-forms on $Z$ in the sense of orbifolds.
Let $\square_{q}=(\bar{\partial}+\bar{\partial}^{*})^{2}$ be the Hodge-Kodaira Laplacian acting on $A^{0,q}_{Z}$.
Let
$$
\zeta_{q}(s):=\sum_{\lambda\in\sigma(\square_{q})\setminus\{0\}}\lambda^{-s}\,\dim E(\lambda;\square_{q})
$$
be the spectral zeta function of $\square_{q}$, where $E(\lambda;\square_{q})$ is the eigenspace of $\square_{q}$
corresponding to the eigenvalue $\lambda$.
Similarly, let
$$
\zeta_{q}(s)(\iota):=\sum_{\lambda\in\sigma(\square_{q})\setminus\{0\}}\lambda^{-s}\,{\rm Tr}\left[\iota^{*}|_{E(\lambda;\square_{q})}\right]
$$
be the equivariant spectral zeta function of $\square_{q}$.
Since $(Z,\gamma)$ is a K\"ahler orbifold, $\zeta_{q}(s)$ and $\zeta_{q}(s)(\iota)$ converge absolutely when $\Re s>\dim Y$,
extend to meromorphic functions on ${\bf C}$, and are holomorphic at $s=0$.
After Ray-Singer \cite{RaySinger73} and Bismut \cite{Bismut95}, we make the following:

\begin{definition}
The analytic torsion of the K\"ahler orbilod $(Z,\gamma)$ is defined as
$$
\tau(Z,\gamma):=\exp[-\sum_{q=0}^{n}(-1)^{q}q\,\zeta'_{q}(0)].
$$
The equivariant analytic torsion of $(Z,\iota,\gamma)$ is defined as
$$
\tau_{{\bf Z}_{2}}(Z,\gamma)(\iota):=\exp[-\sum_{q=0}^{n}(-1)^{q}q\,\zeta'_{q}(0)(\iota)].
$$
\end{definition}

\subsection
{\bf Analytic torsion for $K3$ surfaces}
\par

\begin{theorem}
\label{thm:holomorphic:torsion:K3}
Let $Z$ be a $K3$ surface and let $\eta\in H^{0}(Z,K_{Z})\setminus\{0\}$ and let $\gamma$ be a K\"ahler form on $Z$.
Then the following formula holds:
$$
\tau(Z,\gamma)
=
\exp
\left[
-\frac{1}{24}\int_{Z}
\log\left\{\frac{\eta\wedge\overline{\eta}}{\gamma^{2}/2!}\cdot\frac{{\rm Vol}(Z,\gamma)}{\|\eta\|_{L^{2}}^{2}}\right\}\,c_{2}(Z,\gamma)
\right],
$$
where $c_{i}(Z,\gamma)$ denotes the $i$-th Chern form of $(TZ,\gamma)$ and $\|\eta\|_{L^{2}}^{2}:=\int_{Z}\eta\wedge\overline{\eta}$.
\end{theorem}

\begin{pf}
Let $\omega$ be a Ricci-flat K\"ahler form on $Z$ such that
\begin{equation}
\label{eqn:4:Ricci:flat}
\frac{\omega^{2}}{2!}=\eta\wedge\overline{\eta}.
\end{equation}
Since the $L^{2}$-metric on $H^{2}(Z,{\mathcal O}_{Z})=H^{0}(Z,K_{Z})^{\lor}$ is independent of the choice of a K\"ahler metric on $Z$,
we get by the anomaly formula for Quillen metrics \cite{BGS88}
\begin{equation}
\label{eqn:4:anomaly:1}
\log\left(\frac{\tau(Z,\gamma)\,{\rm Vol}(Z,\gamma)}{\tau(Z,\omega)\,{\rm Vol}(Z,\omega)}\right)
=
\frac{1}{24}\int_{Z}\widetilde{c_{1}c_{2}}(TZ;\gamma,\omega),
\end{equation}
where $\widetilde{c_{1}c_{2}}(TZ;\gamma,\omega)$ is the Bott-Chern secondary class \cite{BGS88} such that
$$
-dd^{c}\widetilde{c_{1}c_{2}}(TZ;\gamma,\omega)=c_{1}(Z,\gamma)c_{2}(Z,\gamma)-c_{1}(Z,\omega)c_{2}(Z,\omega).
$$
Since $c_{1}(Z,\omega)=0$ by the Ricci-flatness of $\omega$, and $\widetilde{c_{1}}(L;h,h')=\log(h/h')$ for a holomorphic line bundle $L$ and Hermitian metrics
$h$ and $h'$ on $L$, and since
$$
\widetilde{c_{1}c_{2}}(TZ;\gamma,\omega)=\widetilde{c_{1}}(TZ;\gamma,\omega)c_{2}(Z,\gamma)+c_{1}(Z,\omega)\widetilde{c_{2}}(TZ;\gamma,\omega)
$$
by \cite{GilletSoule90}, we get by \eqref{eqn:4:Ricci:flat}
\begin{equation}
\label{eqn:4:BottChern}
\widetilde{c_{1}c_{2}}(TZ;\gamma,\omega)
=
\widetilde{c_{1}}(TZ;\gamma,\omega)c_{2}(Z,\gamma)
=
\log\left(\frac{\gamma^{2}}{\omega^{2}}\right)c_{2}(Z,\gamma)
=
\log\left(\frac{\gamma^{2}/2!}{\eta\wedge\overline{\eta}}\right)c_{2}(Z,\gamma).
\end{equation}
Since ${\rm Vol}(Z,\gamma)/{\rm Vol}(Z,\omega)={\rm Vol}(Z,\gamma)/\|\eta\|_{L^{2}}^{2}$,
we get by substituting \eqref{eqn:4:BottChern} into \eqref{eqn:4:anomaly:1}
\begin{equation}
\label{eqn:4:anomaly:2}
\begin{aligned}
\log\left(\frac{\tau(Z,\gamma)}{\tau(Z,\omega)}\right)
&=
-\log\left(\frac{{\rm Vol}(Z,\gamma)}{\|\eta\|_{L^{2}}^{2}}\right)
-
\frac{1}{24}\int_{Z}
\log\left(\frac{\eta\wedge\overline{\eta}}{\gamma^{2}/2!}\right)\,c_{2}(Z,\gamma)
\\
&=
-\frac{1}{24}\int_{Z}
\log\left\{\frac{\eta\wedge\overline{\eta}}{\gamma^{2}/2!}\cdot\frac{{\rm Vol}(Z,\gamma)}{\|\eta\|_{L^{2}}^{2}}\right\}\,c_{2}(Z,\gamma),
\end{aligned}
\end{equation}
where we used the Gauss-Bonnet-Chern formula for $Z$ to get the second equality.
\par
Since $\omega$ is Ricci-flat, the Laplacians $\square_{0}$ and $\square_{2}$ are isospectral via the map $A^{0,0}_{Y}\ni f\mapsto f\overline{\eta}\in A^{0,2}_{Y}$.
Hence, for the Ricci-flat metric $\omega$, we get the equality of meromorphic functions
\begin{equation}
\label{eqn:4:zeta:function:1}
\zeta_{0}(s)=\zeta_{2}(s)
\end{equation}
Since the Dolbeault complex is exact on the orthogonal complement of harmonic forms, we get the equality of meromorphic functions
\begin{equation}
\label{eqn:4:zeta:function:2}
\zeta_{0}(s)-\zeta_{1}(s)+\zeta_{2}(s)=0.
\end{equation}
By \eqref{eqn:4:zeta:function:1} and \eqref{eqn:4:zeta:function:2}, we get
\begin{equation}
\label{eqn:4:vanishing:torsion}
\tau(Z,\omega)=1.
\end{equation}
The result follows from \eqref{eqn:4:anomaly:2} and \eqref{eqn:4:vanishing:torsion}.
\end{pf}

\subsection
{\bf Equivariant analytic torsion for $2$-elementary $K3$ surfaces}
\par
Let $Z$ be a $K3$ surface and let $\iota\colon Z\to Z$ be an anti-symplectic holomorphic involution.
Let $Z^{\iota}=\amalg_{\alpha}C_{\alpha}$ be the decomposition into the connected components.
By Nikulin \cite[Th.\,4.2.2]{Nikulin83}, every $C_{\alpha}$ is a compact Riemann surface unless $Z^{\iota}=\emptyset$.
\par
Let $\gamma$ be an $\iota$-invariant K\"ahler form on $Z$ and let $\eta\in H^{2}(Z,K_{Z})\setminus\{0\}$.
Let
$$
M:=H^{2}(Z,{\bf Z})^{+}
$$
be the invariant sublattice of $H^{2}(Z,{\bf Z})$ with respect to the $\iota$-action. We define
$$
\tau_{M}(Z,\iota)
:=
{\rm Vol}(Z,\gamma)^{\frac{14-r(M)}{4}}\tau_{{\bf Z}_{2}}(Z,\gamma)(\iota)\,
A_{M}(Z,\iota,\gamma)\,
{\rm Vol}(Z^{\iota},\gamma|_{Z^{\iota}})\tau(Z^{\iota},\gamma|_{Z^{\iota}}),
$$
where we define
$$
\tau(Z^{\iota},\gamma|_{Z^{\iota}}):=\prod_{\alpha}\tau(C_{\alpha},\gamma|_{C_{\alpha}}),
\qquad
{\rm Vol}(Z^{\iota},\gamma|_{Z^{\iota}}):=\prod_{\alpha}{\rm Vol}(C_{\alpha},\gamma|_{C_{\alpha}})
$$
and
$$
A_{M}(Z,\iota,\gamma)
:=
\exp
\left[
\frac{1}{8}\int_{Z^{\iota}}
\left.\log\left\{\frac{\eta\wedge\overline{\eta}}{\gamma^{2}/2!}\cdot\frac{{\rm Vol}(Z,\gamma)}{\|\eta\|_{L^{2}}^{2}}\right\}\right|_{Z^{\iota}}\,
c_{1}(Z^{\iota},\gamma|_{Z^{\iota}})
\right].
$$
As before, $c_{1}(Z^{\iota},\gamma|_{Z^{\iota}})$ is the first Chern form of $(TZ^{\iota},\gamma|_{Z^{\iota}})$.

\begin{theorem}
The number $\tau_{M}(Z,\iota)$ is independent of the choice of an $\iota$-invariant K\"ahler form on $Z$.
\end{theorem}

\begin{pf}
See \cite[Th.\,5.7]{Yoshikawa04}.
\end{pf}

For an explicit formula for $\tau_{M}$ as a function on the moduli space of $2$-elementary $K3$ surfaces, 
see \cite{Yoshikawa04}, \cite{Yoshikawa13}, \cite{MaYoshikawa15}.
By Theorem \ref{thm:holomorphic:torsion:K3}, we can rewrite $\tau_{M}(Z,\iota)$ as follows
\begin{equation}
\label{eqn:4:torsion:invariant}
\begin{aligned}
\tau_{M}(Z,\iota)
&=
{\rm Vol}(Z,\gamma)^{\frac{14-r(M)}{4}}\tau(Z,\gamma)\tau_{{\bf Z}_{2}}(Z,\gamma)(\iota)\,
{\rm Vol}(Z^{\iota},\gamma|_{Z^{\iota}})\tau(Z^{\iota},\gamma|_{Z^{\iota}})
\\
&\quad\times
A_{M}(Z,\iota,\gamma)\,
\exp
\left[
\frac{1}{24}\int_{Z}
\log\left\{\frac{\eta\wedge\overline{\eta}}{\gamma^{2}/2!}\cdot\frac{{\rm Vol}(Z,\gamma)}{\|\eta\|_{L^{2}}^{2}}\right\}\,c_{2}(Z,\gamma)
\right].
\end{aligned}
\end{equation}

\section
{A degenerating family of K\"ahler metrics}
\par
Let $Y$ be a good log-Enriques surface. For an orbifold K\"ahler form  $\gamma$ on $Y$,
we write ${\rm Vol}(Y,\gamma)=\int_{Y}\gamma^{2}/2!$ for the volume of $(Y,\gamma)$. We set
$$
k:=\#{\rm Sing}(Y)\in\{1,\ldots,10\}.
$$
Let $(\widetilde{X},\theta)$ be the $2$-elementary $K3$ surface associated to $Y$ such that
$$
\widetilde{X}^{\theta}=\amalg_{{\frak p}\in{\rm Sing}(Y)}E_{\frak p},
\qquad
E_{\frak p}\cong{\bf P}^{1}.
$$
Let
$$
\pi\colon(\widetilde{X},\widetilde{X}^{\theta})\to(X,{\rm Sing}\,X)
$$
be the blowing-down of the disjoint union of $(-2)$-curves. Then
$$
{\frak p}=\pi(E_{\frak p}).
$$
In this section, we construct a two parameter family of K\"ahler metrics $\{ \gamma_{\epsilon,\delta} \}$ on $\widetilde{X}$
converging to an orbifold K\"ahler metric on $X$, which is obtained by gluing the Eguchi-Hanson instanton at each ${\mathfrak p}$ 
and a K\"ahler metric on $X$. In the subsequent sections, we study the limiting behavior of various geometric quantities of
$(\widetilde{X}, \gamma_{\epsilon,\delta})$ to construct an invariant of the log-Enriques surface $Y$.

\subsection
{\bf Eguchi-Hanson instanton}
\par
For $\epsilon\geq0$, let $F_{\epsilon}(z)$ be the function on ${\bf C}^{2}\setminus\{0\}$ defined by
$$
F_{\epsilon}(z)
:=
\sqrt{\|z\|^{4}+\epsilon^{2}}+\epsilon\,\log\left(\frac{\|z\|^{2}}{\sqrt{\|z\|^{4}+\epsilon^{2}}+\epsilon}\right).
$$
On every compact subset of ${\bf C}^{2}\setminus\{0\}$, we have $\lim_{\epsilon\to0}F_{\epsilon}(z)=\|z\|^{2}$.
For all $\epsilon\geq0$ and $\delta>0$,
$$
F_{\epsilon}(\delta z)=\delta^{2}F_{\epsilon\delta^{-2}}(z).
$$
Let $T^{*}{\bf P}^{1}$ be the holomorphic cotangent bundle of the projective line and let $E\subset T^{*}{\bf P}^{1}$ be its zero section.
Let
$$
\varPi\colon(T^{*}{\bf P}^{1},E)\to({\bf C}^{2}/\{\pm1\},0)
$$
be the blowing-down of the zero section. Since
$$
i\partial\bar{\partial}F_{\epsilon}(z)
=
i
\left(
\frac{\epsilon\,\partial\|z\|^{2}\wedge\bar{\partial}\|z\|^{2}}{\sqrt{(\|z\|^{4}+\epsilon^{2}}+\epsilon)\sqrt{\|z\|^{4}+\epsilon^{2}}}
+
\frac{\|z\|^{2}\partial\bar{\partial}\|z\|^{2}}{\sqrt{\|z\|^{4}+\epsilon^{2}}+\epsilon}
+
\epsilon\,\partial\bar{\partial}\log\|z\|^{2}
\right)
$$
is a positive $(1,1)$-form on $({\bf C}^{2}\setminus\{0\})/\pm1$ satisfying
$$
\frac{(i\partial\bar{\partial}F_{\epsilon})^{2}}{2!}=(\sqrt{-1})^{2}dz_{1}\wedge d\bar{z}_{1}\wedge dz_{2}\wedge d\bar{z}_{2},
$$
its pull-back to $T^{*}{\bf P}^{1}$
$$
\gamma^{\rm EH}_{\epsilon}:=\varPi^{*}(i\partial\bar{\partial}F_{\epsilon})
$$
extends to a Ricci-flat K\"ahler form on $T^{*}{\bf P}^{1}$ for $\epsilon>0$, called the {\em Eguchi-Hanson instanton}.
We write $\gamma^{\rm EH}$ for $\gamma^{\rm EH}_{1}$. The coordinate change $z\mapsto \sqrt{\epsilon}z$ on ${\bf C}^{2}$ induces
an isometry of K\"ahler manifolds 
\begin{equation}
\label{eqn:scaling:property:EH}
(T^{*}{\bf P}^{1},\gamma^{\rm EH}_{\epsilon})\cong(T^{*}{\bf P}^{1},\epsilon\gamma^{\rm EH}).
\end{equation}
When $\epsilon=0$,
$$
i\partial\bar{\partial}F_{0}=i\partial\bar{\partial}\|z\|^{2}
$$
is the Euclidean K\"ahler form on ${\bf C}^{2}/\{\pm1\}$, and
$\gamma^{\rm EH}_{0}=\varPi^{*}\{i\partial\bar{\partial}F_{0}(z)\}=\pi^{*}(i\partial\bar{\partial}\|z\|^{2})$
is a degenerate K\"ahler form on $T^{*}{\bf P}^{1}$.
\par
Let $\omega_{\rm FS}$ be the Fubini-Study form on ${\bf P}^{1}$ such that
$$
[\omega_{\rm FS}]=c_{1}({\mathcal O}_{{\bf P}^{1}}(1)).
$$
By the definition of $F_{\epsilon}$, we get
$$
\gamma^{\rm EH}_{\epsilon}|_{E}
=
\epsilon\varPi^{*}(i\partial\bar{\partial}\log\|z\|^{2})|_{E}
=
2\pi\epsilon\,\omega_{\rm FS}.
$$

\subsection
{\bf Glueing of the Eguchi-Hanson instanton}
\par

\subsubsection
{A modification of the Eguchi-Hanson instanton}
\par
Let $B(r)\subset{\bf C}^{2}$ be the ball of radius $r>0$ centered at $0\in{\bf C}^{2}$ and set
$$
V(r):=B(r)/\{\pm1\}.
$$
Let $\varPi\colon(\widetilde{V}(r),E)\to(V(r),0)$ be the blowing-up at the origin. 
Then $V(\infty)={\bf C}^{2}/\pm1$ and $\widetilde{V}(\infty)=T^{*}{\bf P}^{1}$.
For $z\in({\bf C}^{2}\setminus\{0\})/\pm1$ and $\epsilon\geq0$, we define
$$
E(z,\epsilon):=F_{\epsilon}(z)-\|z\|^{2}.
$$
Since the error term $E(z,\epsilon)$ is a $C^{\omega}$ function on $(V(4)\setminus V(1))\times[0,1]$ with $E(z,0)=0$,
there is a constant $C_{k}$ for all $k\geq0$ with
\begin{equation}
\label{eqn:5:estimate:error}
\sup_{z\in V(4)\setminus V(1)}|\partial_{z}^{k}E(z,\epsilon)|\leq C_{k}\,\epsilon.
\end{equation}
\par
Let $\rho(t)$ be a $C^{\infty}$ function on ${\bf R}$ such that $0\leq\rho(t)\leq1$ on ${\bf R}$,  $\rho(t)=1$ for $t\leq1$ and $\rho(t)=0$ for $t\geq2$.
We set
$$
\phi_{\epsilon}(z)
:=
\rho(\|z\|)\,F_{\epsilon}(z)+\{1-\rho(\|z\|)\}\,\|z\|^{2}
=
\|z\|^{2}+\rho(\|z\|)\,E(z,\epsilon)
$$
and we define a $(1,1)$ form on $V(\infty)\setminus\{0\}$ by
$$
\kappa_{\epsilon}:=i\partial\bar{\partial}\phi_{\epsilon}.
$$
Since $\phi_{\epsilon}(z)=F_{\epsilon}(z)$ on $V(1)$, $\kappa_{\epsilon}$ extends to a real $(1,1)$-form on $T^{*}{\bf P}^{1}$,
which is positive on $\widetilde{V}(1)$. 
Since $\phi_{\epsilon}(z)=\|z\|^{2}+\rho(\|z\|)\,E(z,\epsilon)$ on $V(2)\setminus V(1)$,
there exists by \eqref{eqn:5:estimate:error} a constant $\epsilon(\rho)\in(0,1)$ depending only on the choice of the cut-off function $\rho$
such that $\kappa_{\epsilon}$ is a positive $(1,1)$-form on $V(2)\setminus V(1)$ for $0<\epsilon\leq\epsilon(\rho)$.
As a result, $\{\kappa_{\epsilon}\}_{0<\epsilon\leq\epsilon(\rho)}$ is a family of K\"ahler forms on $T^{*}{\bf P}^{1}$ such that
$\kappa_{\epsilon}=i\partial\bar{\partial}\|z\|^{2}$ on $T^{*}{\bf P}^{1}\setminus\widetilde{V}(2)$.
\newline

We have the following slightly refined estimate for the error term $E(z, \epsilon)$. Set 
$$
E(z):=E(z,1)=E_{1}(z)+E_{2}(z),
$$
where 
$$
E_{1}(z)=\sqrt{\|z\|^{4}+1}-\|z\|^{2}=\frac{1}{\sqrt{\|z\|^{4}+1}+\|z\|^{2}},
\quad
E_{2}(z)=\log\frac{\|z\|^{2}}{\sqrt{\|z\|^{4}+1}+1}.
$$
Then, for any nonnegative integer $k$, there exists a constant $C_{k}>0$ such that 
\begin{itemize}
\item[(i)]
$|\partial_{z}^{k}E_{1}(z)|\leq C_{k}(1+\|z\|)^{-(2+k)}$ for all $z\in V(\infty)\setminus\{0\}$;
\item[(ii)]
$|\partial_{z}^{k}E_{2}(z)|\leq C_{k}(1+\|z\|)^{-(2+k)}$ for all $z\in V(\infty)\setminus V(2)$;
\item[(iii)]
$|\partial_{z}^{k}E_{2}(z)|\leq C_{k}\|z\|^{-k}$ for all $z\in V(2)\setminus\{0\}$ $k\geq 1$; $C_0 \log \|z\|^2$ for $k=0$.
\end{itemize}
From these inequalities, we get
\begin{equation}
\label{eqn:5:estimate1:error}
|\partial_{z}^{k}E(z)|
\leq 
\begin{cases}
\begin{array}{ll}
C_{k}\|z\|^{-k} \ (k\geq 1; C_0 \log \|z\|^2, k=0)
&
(\forall\,z\in V(2)\setminus\{0\}),
\\
C_{k}(1+\|z\|)^{-(2+k)}
&
(\forall\,z\in V(\infty)\setminus V(2)).
\end{array}
\end{cases}
\end{equation}
Since $E(z,\epsilon)=\epsilon\,E(\frac{z}{\sqrt{\epsilon}},1)=\epsilon\,E(\frac{z}{\sqrt{\epsilon}})$ and hence
$\partial_{z}^{k}E(z,\epsilon)=\epsilon^{1-\frac{k}{2}}(\partial_{z}^{k}E)(\frac{z}{\sqrt{\epsilon}})$, we get by \eqref{eqn:5:estimate1:error}
\begin{equation}
\label{eqn:5:estimate2:error}
|\partial_{z}^{k}E(z,\epsilon)|
\leq 
\begin{cases}
\begin{array}{ll}
C_{k}\epsilon\|z\|^{-k} \ (k\geq 1; C_0 \epsilon (\log \|z\|^2 + \log \epsilon), k=0)
&
(\forall\,z\in V(2)\setminus\{0\}),
\\
C_{k}\epsilon^{2}(\sqrt{\epsilon}+\|z\|)^{-(2+k)}
&
(\forall\,z\in V(\infty)\setminus V(2)).
\end{array}
\end{cases}
\end{equation}
Here, to get the estimate on $V(2)\setminus\{0\}$, we used the fact $\epsilon^{2}(\sqrt{\epsilon} + \|z\|)^{-(2+k)} < \epsilon\|z\|^{-k}$
on $V(2) \setminus V(2\sqrt{\epsilon})$.
Replacing $\epsilon(\rho)$ by a smaller constant if necessary, we may assume by \eqref{eqn:5:estimate2:error} 
the following inequality of Hermitian matrices for all $0<\epsilon\leq\epsilon(\rho)$ and $z\in V(\infty)\setminus V(2)$:
\begin{equation}
\label{eqn:5:estimate4:error}
\frac{1}{2}(\delta_{ij})
\leq
(\delta_{ij}+\frac{\partial^{2}E(z,\epsilon)}{\partial z_{i}\partial\bar{z}_{j}})
\leq
2(\delta_{ij}).
\end{equation}

Moreover, for $\|z\| \leq 2$,
\begin{equation}
\label{eqn:5:estimate3:error}
|\partial_{\epsilon}E(z,\epsilon)|\leq C\,\epsilon \|z\|^{-2}.
\end{equation}

\begin{lemma}
\label{lemma:comparison:1}
There exist constants $C_{1},C_{2}>0$ such that 
the following inequality of $(1,1)$-forms on $T^{*}{\bf P}^{1}$ hold for all $0<\epsilon\leq\epsilon(\rho)$:
$$
C_{1}\gamma_{\epsilon}^{\rm EH}\leq\kappa_{\epsilon}\leq C_{2}\gamma_{\epsilon}^{\rm EH}.
$$
\end{lemma}

\begin{pf}
{\em (Step 1) } 
On $\widetilde{V}(1)$, we have $\kappa_{\epsilon}=\gamma_{\epsilon}^{\rm EH}$.
On $\widetilde{V}(2)\setminus\widetilde{V}(1)$, it follows from \eqref{eqn:5:estimate:error} that there exist constants $C_{1},C_{2}>0$ 
independent of $\epsilon\in(0,\epsilon(\rho)]$ with
$C_{1}\gamma_{\epsilon}^{\rm EH}\leq\kappa_{\epsilon}\leq C_{2}\gamma_{\epsilon}^{\rm EH}$. 
Combining these two estimates, we get $C_{1}\gamma_{\epsilon}^{\rm EH}\leq\kappa_{\epsilon}\leq C_{2}\gamma_{\epsilon}^{\rm EH}$ 
on $\widetilde{V}(2)$.
\par{\em (Step 2) } 
We compare $\kappa_{\epsilon}$ and $\gamma_{\epsilon}^{\rm EH}$ on $T^{*}{\bf P}^{1}\setminus\widetilde{V}(2)$.
On $T^{*}{\bf P}^{1}\setminus\widetilde{V}(2)$, we have $\kappa_{\epsilon}=\gamma^{\rm EH}_{0}$.
By \eqref{eqn:5:estimate4:error}, we have
$\frac{1}{2}\gamma^{\rm EH}_{\epsilon}\leq\gamma^{\rm EH}_{0}\leq2\gamma^{\rm EH}_{\epsilon}$ 
on $T^{*}{\bf P}^{1}\setminus\widetilde{V}(2)$.
Since $\kappa_{\epsilon}=\gamma^{\rm EH}_{0}$ on $T^{*}{\bf P}^{1}\setminus\widetilde{V}(2)$,
We get the desired estimate on $T^{*}{\bf P}^{1}\setminus\widetilde{V}(2)$.
This completes the proof.
\end{pf}

\subsubsection
{A family of K\"ahler metrics on $\widetilde{X}$}
\par
Since $E_{\frak p}$ is a $(-2)$-curve on $\widetilde{X}$, there exist a neighborhood $U_{\frak p}$ of $E_{\frak p}$ in $\widetilde{X}$ 
and an isomorphism of pairs
$$
\psi_{\frak p}\colon(U_{\frak p},E_{\frak p})\cong(\widetilde{V}(1),E).
$$
We may and will assume that $\psi_{\frak p}$ extends to
an isomorphism between an open subset of $\widetilde{X}$ containing $U_{\frak p}$ and $\widetilde{V}(4)$.
We write $V(r)_{{\frak p}}$ for $V(r)$ viewed as a neighborhood of ${\frak p}\in{\rm Sing}(X)$.
In what follows, we identify $\widetilde{V}(r)_{{\frak p}}$ with $\psi_{\frak p}^{-1}(\widetilde{V}(r)_{{\frak p}})$.
\par
Let $\gamma$ be a $\theta$-invariant K\"ahler form on $X$ in the sense of orbifolds, which has a potential function on every $V(4)_{\frak p}$.
By modifying the potential of $\gamma$ on each $V(4)_{\frak p}$ (cf. \cite[Proof of Lemma 6.2]{Yoshikawa04}),
there exists a K\"ahler form $\gamma_{0}$ on $X$ in the sense of orbifolds
such that 
\begin{equation}
\label{eqn:initial:metric}
\gamma_{0}|_{X\setminus\bigcup_{{\frak p}\in{\rm Sing}(X)}V(2)_{\frak p}}=\gamma,
\qquad
\gamma_{0}|_{V(2)_{\frak p}}=i\partial\bar{\partial}\|z\|^{2}
\quad
(\forall\,{\frak p}\in{\rm Sing}(X)).
\end{equation}
In particular, $\|z\|^{2}\in C^{\omega}(V(2)_{\frak p})$ is a potential function of $\gamma_{0}$ on every $V(2)_{\frak p}$.
Since $\phi_{\epsilon}(z)=\|z\|^{2}$ near $\partial V(2)_{\frak p}$, we can glue the K\"ahler form $\kappa_{\epsilon}$ on 
$\bigcup_{{\frak p}\in{\rm Sing}(X)}\widetilde{V}(2)_{\frak p}$ and the K\"ahler form $\gamma_{0}$ on 
$X\setminus\bigcup_{{\frak p}\in{\rm Sing}(X)}\widetilde{V}(2)_{\frak p}$ by setting
\begin{equation}
\label{eqn:family:Kahler:metric:1}
\gamma_{\epsilon}
:=
\begin{cases}
\begin{array}{ll}
\kappa_{\epsilon}&\hbox{on }\bigcup_{{\frak p}\in{\rm Sing}(X)}\widetilde{V}(2)_{\frak p},
\\
\gamma_{0}&\hbox{on }X\setminus\bigcup_{{\frak p}\in{\rm Sing}(X)}\widetilde{V}(2)_{\frak p}.
\end{array}
\end{cases}
\end{equation}
By construction, $\{\gamma_{\epsilon}\}_{0<\epsilon\leq\epsilon(\rho)}$ is a family of $\theta$-invariant K\"ahler forms on $\widetilde{X}$.

\begin{lemma}
\label{lemma:Ricci:1}
The family of K\"ahler forms $\{\gamma_{\epsilon}\}_{0<\epsilon\leq\epsilon(\rho)}$ on $\widetilde{X}$ satisfies the following:
\begin{itemize}
\item[(1)]
For all ${\frak p}\in{\rm Sing}(X)$, $\gamma_{0}|_{V(2)_{\frak p}}=i\partial\bar{\partial}\|z\|^{2}$.
\item[(2)]
For all ${\frak p}\in{\rm Sing}(X)$, $\gamma_{\epsilon}|_{\widetilde{V}(1)_{\frak p}}=\psi_{\frak p}^{*}\gamma^{\rm EH}_{\epsilon}$.
\item[(3)]
On $\widetilde{X}$, $\gamma_{\epsilon}$ converges to $\pi^{*}\gamma_{0}$ in the $C^{\infty}$-topology.
\item[(4)]
There exist constants $C, C'>0$ independent of $\epsilon$ (but depending on $\rho$) such that
$|{\rm Ric}(\gamma_{\epsilon})|_{\gamma_{\epsilon}}\leq C\cdot\epsilon$ on $\bigcup_{{\frak p}\in{\rm Sing}\,X}\widetilde{V}(2)_{\frak p}$ and
$|{\rm Ric}(\gamma_{\epsilon})|_{\gamma_{\epsilon}}\leq C'$ on $\widetilde{X}$.
\end{itemize}
\end{lemma}

\begin{pf}
By construction, (1), (2), (3) are obvious. Let us see (4).
Since $\gamma_{\epsilon}^{\rm EH}$ is Ricci-flat and since $\kappa_{\epsilon}=\gamma_{\epsilon}^{\rm EH}$ on $\widetilde{V}(1)_{\frak p}$, 
we get ${\rm Ric}(\kappa_{\epsilon})={\rm Ric}(\gamma_{\epsilon}^{\rm EH})=0$ on $\widetilde{V}(1)_{\frak p}$. 
On $\widetilde{V}(2)_{\frak p}\setminus\widetilde{V}(1)_{\frak p}$, we get
$|{\rm Ric}(\gamma_{\epsilon})|_{\gamma_{\epsilon}}=|{\rm Ric}(\kappa_{\epsilon})|_{\kappa_{\epsilon}}\leq C\cdot\epsilon$ 
by \eqref{eqn:5:estimate:error}. This proves the first estimate. Since $\gamma_{\epsilon}=\gamma_{0}$ on
$X\setminus\bigcup_{{\frak p}\in{\rm Sing}(X)}\widetilde{V}(2)_{\frak p}$, we get the second estimate.
\end{pf}

\subsubsection
{A two parameter family of K\"ahler metrics on $T^{*}{\bf P}^{1}$}
\par
For later use we introduce another small parameter $\delta>0$. Instead of gluing in the Eguchi-Hanson instanton in the region 
$\widetilde{V}(2)-\widetilde{V}(1)$ we now do it in the region $\widetilde{V}(2\delta)-\widetilde{V}(\delta)$. 
This is effected by replacing the cut-off function $\rho(t)$ by $\rho_{\delta}(t)=\rho(\frac{t}{\delta})$ in defining the K\"ahler potential $\phi_{\epsilon}$ 
for the K\"ahler metric $\gamma_{\epsilon}$ such that $\rho_{\delta}(t)=1$ for $t\leq \delta$ and $\rho_{\delta}(t)=0$ for $t\geq 2\delta$. 
This gives us the family of real $(1,1)$-forms on $T^{*}{\bf P}^{1}$
$$
\kappa_{\epsilon, \delta}:= i \partial \overline{\partial} \phi_{\epsilon, \delta},
$$ 
where 
$$
\phi_{\epsilon, \delta}(z):=\|z\|^2 + \rho_{\delta}(\|z\|)E(z, \epsilon).
$$
To verify the positivity of $\kappa_{\epsilon, \delta}$, we see the relation between $\phi_{\epsilon}$ and $\phi_{\epsilon, \delta}$.
Since $F_{\epsilon}(\delta z)=\delta^{2}F_{\epsilon/\delta^{2}}(z)$, we get $E(\delta z,\epsilon)=\delta^{2}E(z,\epsilon/\delta^{2})$.
Since $\phi_{\epsilon,1}(z)=\phi_{\epsilon}(z)$ and
$\phi_{\epsilon, \delta}(z)=\|\delta\cdot\frac{z}{\delta}\|^{2}+\rho(\frac{\|z\|}{\delta})E(\delta\cdot\frac{z}{\delta},\epsilon)$, 
this implies that
$$
\phi_{\epsilon, \delta}(z)=\delta^{2}\phi_{\epsilon/\delta^{2}}(z/\delta).
$$
Hence if $0<\epsilon/\delta^{2}\leq\epsilon(\rho)$, then $\kappa_{\epsilon, \delta}=i\partial\bar{\partial}\phi_{\epsilon, \delta}$
is a positive $(1,1)$-form on $T^{*}{\bf P}^{1}$. In what follows, we define $\phi_{\epsilon, \delta}$ for 
$\epsilon,\delta\in(0,1]$ with $0<\epsilon/\delta^{2}\leq\epsilon(\rho)$.
Then $\{\kappa_{\epsilon,\delta}\}_{0<\epsilon/\delta^{2}\leq\epsilon(\rho),\,\epsilon,\delta\in(0,1]}$ 
is a family of K\"ahler forms on $T^{*}{\bf P}^{1}$.
Moreover, the relation $\phi_{\epsilon, \delta}(z)=\delta^{2}\phi_{\epsilon/\delta^{2}}(z/\delta)$
implies that the automorphism of $T^{*}{\bf P}^{1}$ induced from the one $z\mapsto z/\delta$ on $V(\infty)$
yields an isometry of K\"ahler manifolds
$(T^{*}{\bf P}^{1},\kappa_{\epsilon, \delta})\cong(T^{*}{\bf P}^{1},\delta^{2}\kappa_{\epsilon/\delta^{2}})$
such that
\begin{equation}
\label{eqn:isometry:1}
(\widetilde{V}(2\delta),\kappa_{\epsilon, \delta})\cong(\widetilde{V}(2),\delta^{2}\kappa_{\epsilon/\delta^{2}}).
\end{equation}

\begin{lemma}
\label{lemma:comparison:2}
There exist constants $C_{1},C_{2}>0$ such that 
the following inequality of $(1,1)$-forms on $T^{*}{\bf P}^{1}$ holds for all $\epsilon,\delta\in(0,1]$ with $0<\epsilon/\delta^{2}\leq\epsilon(\rho)$:
$$
C_{1}\kappa_{\epsilon}\leq\kappa_{\epsilon,\delta}\leq C_{2}\kappa_{\epsilon}.
$$
\end{lemma}

\begin{pf}
{\em (Step 1) } 
By Lemma~\ref{lemma:comparison:2} (Step 1), 
we get $C_{1}\gamma_{\epsilon}^{\rm EH}\leq\kappa_{\epsilon}\leq C_{2}\gamma_{\epsilon}^{\rm EH}$ on $\widetilde{V}(2)$.
By \eqref{eqn:isometry:1} and the relation $\delta^{2}\gamma_{\epsilon/\delta^{2}}^{\rm EH}=\gamma_{\epsilon}^{\rm EH}$, 
this implies the inequality
$C_{1}\gamma_{\epsilon}^{\rm EH}\leq\kappa_{\epsilon,\delta}\leq C_{2}\gamma_{\epsilon}^{\rm EH}$ on $\widetilde{V}(2\delta)$.
Hence we get $C_{1}C_{2}^{-1}\kappa_{\epsilon}\leq\kappa_{\epsilon,\delta}\leq C_{2}C_{1}^{-1}\kappa_{\epsilon}$ on $\widetilde{V}(2\delta)$.
\par{\em (Step 2) } 
Next we compare $\kappa_{\epsilon,\delta}$ and $\kappa_{\epsilon}$ on $T^{*}{\bf P}^{1}\setminus\widetilde{V}(2\delta)$.
By definition, we have $\kappa_{\epsilon,\delta}=\gamma_{0}^{\rm EH}$ on $T^{*}{\bf P}^{1}\setminus\widetilde{V}(2\delta)$.
Let $H_{\epsilon}$ be the automorphism of $T^{*}{\bf P}^{1}$ induced from the automorphism $z\mapsto\sqrt{\epsilon}z$ of $V(\infty)={\bf C}^{2}/\pm1$.
Then $H_{\epsilon}$ is an isomorphism from 
$T^{*}{\bf P}^{1}\setminus\widetilde{V}(2\delta/\sqrt{\epsilon})$ to $T^{*}{\bf P}^{1}\setminus\widetilde{V}(2\delta)$ inducing the isometries
\begin{equation}
\label{eqn:isometry:2}
(T^{*}{\bf P}^{1}\setminus\widetilde{V}(2\delta),\gamma_{\epsilon}^{\rm EH})
\cong
(T^{*}{\bf P}^{1}\setminus\widetilde{V}(2\delta/\sqrt{\epsilon}),\epsilon\gamma^{\rm EH}),
\end{equation}
\begin{equation}
\label{eqn:isometry:3}
(T^{*}{\bf P}^{1}\setminus\widetilde{V}(2\delta),\gamma_{0}^{\rm EH})
\cong
(T^{*}{\bf P}^{1}\setminus\widetilde{V}(2\delta/\sqrt{\epsilon}),\epsilon\gamma_{0}^{\rm EH}).
\end{equation}
Since $\epsilon/\delta^{2}\leq\epsilon(\rho)$ and hence $\delta/\sqrt{\epsilon}>1/\sqrt{\epsilon(\rho)}$,
we have the inclusion 
$T^{*}{\bf P}^{1}\setminus\widetilde{V}(\frac{2\delta}{\sqrt{\epsilon}}) \subset T^{*}{\bf P}^{1}\setminus\widetilde{V}(2/\sqrt{\epsilon(\rho)})$.
By \eqref{eqn:5:estimate2:error}, there exist constants $C'_{1},C'_{2}>0$ such that
$C'_{1}\gamma^{\rm EH}\leq\gamma^{\rm EH}_{0}\leq C'_{2}\gamma^{\rm EH}$
on $T^{*}{\bf P}^{1}\setminus\widetilde{V}(2/\sqrt{\epsilon(\rho)})$.
This, together with \eqref{eqn:isometry:2}, \eqref{eqn:isometry:3}, yields the inequality
$C'_{1}\gamma^{\rm EH}_{\epsilon}\leq\gamma^{\rm EH}_{0}\leq C'_{2}\gamma^{\rm EH}_{\epsilon}$ on 
$T^{*}{\bf P}^{1}\setminus\widetilde{V}(2\delta)$ for all $\epsilon,\delta\in(0,1]$ with $0<\epsilon/\delta^{2}\leq\epsilon(\rho)$.
Since $\kappa_{\epsilon,\delta}=\gamma_{0}^{\rm EH}$ on $T^{*}{\bf P}^{1}\setminus\widetilde{V}(2\delta)$, we get
$C'_{1}\gamma^{\rm EH}_{\epsilon}\leq\kappa_{\epsilon,\delta}\leq C'_{2}\gamma^{\rm EH}_{\epsilon}$ on 
$T^{*}{\bf P}^{1}\setminus\widetilde{V}(2\delta)$.
By Lemma~\ref{lemma:comparison:1}, this implies the inequality 
$C''_{1}\kappa_{\epsilon}\leq\kappa_{\epsilon,\delta}\leq C''_{2}\kappa_{\epsilon}$ on 
$T^{*}{\bf P}^{1}\setminus\widetilde{V}(2\delta)$, where $C''_{1},C''_{2}>0$ are constants independent of 
$\epsilon,\delta\in(0,1]$ with $0<\epsilon/\delta^{2}\leq\epsilon(\rho)$.
This completes the proof.
\end{pf}

\subsubsection
{A two parameter family of K\"ahler metrics on $\widetilde{X}$}
\par
Modifying the construction \eqref{eqn:family:Kahler:metric:1}, we introduce a two parameter family of $\theta$-invariant K\"ahler forms 
on $\widetilde{X}$ by
\begin{equation}
\label{eqn:family:Kahler:metric:2}
\gamma_{\epsilon,\delta}
:=
\begin{cases}
\begin{array}{ll}
\kappa_{\epsilon,\delta}&\hbox{on }\bigcup_{{\frak p}\in{\rm Sing}(X)}\widetilde{V}(2)_{\frak p},
\\
\gamma_{0}&\hbox{on }X\setminus\bigcup_{{\frak p}\in{\rm Sing}(X)}\widetilde{V}(2)_{\frak p}
\end{array}
\end{cases}
\end{equation}
for $\epsilon,\delta\in(0,1]$ with $0<\epsilon/\delta^{2}\leq\epsilon(\rho)$.

\begin{lemma}
\label{lemma:comparison:3}
There exist constants $C_{1},C_{2}>0$ such that 
the following inequality of $(1,1)$-forms on $\widetilde{X}$ hold for all $\epsilon,\delta\in(0,1]$ with $0<\epsilon/\delta^{2}\leq\epsilon(\rho)$:
$$
C_{1}\gamma_{\epsilon}\leq\gamma_{\epsilon,\delta}\leq C_{2}\gamma_{\epsilon}.
$$
\end{lemma}

\begin{pf}
On $\bigcup_{{\frak p}\in{\rm Sing}(X)}\widetilde{V}(2)_{\frak p}$, the result follows from Lemma~\ref{lemma:comparison:2}.
On $X\setminus\bigcup_{{\frak p}\in{\rm Sing}(X)}\widetilde{V}(2)_{\frak p}$, the result is obvious since
$\gamma_{\epsilon,\delta}=\gamma_{\epsilon}=\gamma_{0}$ is independent of $\epsilon,\delta$ there.
\end{pf}

\begin{lemma}
\label{lemma:Ricci}
There exists a constant $C_{3}>0$ such that 
the following estimate holds for all $\epsilon,\delta\in(0,1]$ with $0<\epsilon/\delta^{2}\leq\epsilon(\rho)$:
$$
|{\rm Ric}(\gamma_{\epsilon,\delta})|_{\gamma_{\epsilon,\delta}}\leq C_{3}(\epsilon\delta^{-4}+1).
$$
\end{lemma}

\begin{pf}
Since $\gamma_{\epsilon,\delta}=\gamma_{0}$ on $X\setminus\bigcup_{{\frak p}\in{\rm Sing}(X)}\widetilde{V}(2)_{\frak p}$,
it suffices to prove the estimate on $\bigcup_{{\frak p}\in{\rm Sing}(X)}\widetilde{V}(2)_{\frak p}$.
Since $\gamma_{\epsilon,\delta}=i\partial\bar{\partial}\|z\|^{2}$ is a flat metric on
$\bigcup_{{\frak p}\in{\rm Sing}(X)}\widetilde{V}(2)_{\frak p}\setminus\widetilde{V}(2\delta)_{\frak p}$,
it suffices to prove the estimate on $\bigcup_{{\frak p}\in{\rm Sing}(X)}\widetilde{V}(2\delta)_{\frak p}$.
By \eqref{eqn:isometry:1}, we get on each $\widetilde{V}(2\delta)_{\frak p}$
$$
|{\rm Ric}(\gamma_{\epsilon,\delta})|_{\gamma_{\epsilon,\delta}}
=
|{\rm Ric}(\delta^{2}\gamma_{\epsilon/\delta^{2}})|_{\gamma_{\epsilon,\delta}}
=
\delta^{-2}|{\rm Ric}(\gamma_{\epsilon/\delta^{2}})|_{\gamma_{\epsilon/\delta^{2}}}
\leq
\delta^{-2}C(\epsilon/\delta^{2})
=
C\epsilon\delta^{-4},
$$
where we used Lemma~\ref{lemma:Ricci:1} (4) to get the inequality
$|{\rm Ric}(\gamma_{\epsilon/\delta^{2}})|_{\gamma_{\epsilon/\delta^{2}}}\leq C(\epsilon/\delta^{2})$ on $\widetilde{V}(2\delta)_{\frak p}$. 
This completes the proof.
\end{pf}

\par
Fix a nowhere vanishing holomorphic $2$-form
$$
\eta\in H^{0}(\widetilde{X},K_{\widetilde{X}})\setminus\{0\}.
$$
Since $(\varPi^{-1})^{*}(\eta|_{V(1)_{\frak p}})$ is a nowhere vanishing holomorphic $2$-form on $V(1)_{\frak p}\setminus\{0\}$,
there exists by the Hartogs extension theorem a nowhere vanishing holomorphic function $f_{\frak p}(z)$ on $B(1)$ such that
$$
(\varPi^{-1})^{*}(\eta|_{V(1)_{\frak p}})=f_{\frak p}(z)\,dz_{1}\wedge dz_{2}
$$
and $f_{\frak p}(-z)=f_{\frak p}(z)$.
Since $\gamma_{\epsilon,\delta}=\gamma_{\epsilon}$ on $\widetilde{V}(\delta)_{\frak p}$ and hence
$$
\left.(\varPi^{-1})^{*}(\gamma_{\epsilon,\delta}^{2}/2!)\right|_{V(\delta)_{\frak p}\setminus\{0\}}
=
(i\partial\bar{\partial}F_{\epsilon})^{2}/2!
=
(i)^{2}dz_{1}\wedge d\bar{z}_{1}\wedge dz_{2}\wedge d\bar{z}_{2},
$$
we get the equality of functions on $\widetilde{V}(\delta)_{\frak p}$
\begin{equation}
\label{eqn:Ricci:flatness}
\left.\frac{\eta\wedge\overline{\eta}}{\gamma_{\epsilon,\delta}^{2}/2!}\right|_{\widetilde{V}(\delta)_{\frak p}}
=
\pi^{*}|f_{\frak p}(z)|^{2}.
\end{equation}
In particular, we have the following:
\begin{itemize}
\item[(i)]
On each $\widetilde{V}(\epsilon)_{\frak p}$, the volume form of $\gamma_{\epsilon,\delta}$ is independent of $\epsilon\in(0,\epsilon(\rho)]$.
\item[(ii)]
$f_{\frak p}(0)$ is independent of $\delta\in(0,1]$ and the choice of the cut-off function $\rho$.
\end{itemize}
Since $\gamma_{\epsilon,\delta}$ converges to $\gamma_{0}$ outside $\bigcup_{{\frak p}\in{\rm Sing}(X)}\widetilde{V}(\delta)_{\frak p}$,
we get the continuity
\begin{equation}
\label{eqn:continuity:volume}
\lim_{\epsilon\to0}{\rm Vol}(\widetilde{X},\gamma_{\epsilon,\delta})={\rm Vol}(X,\gamma_{0}).
\end{equation}

\subsection
{\bf Ricci-flat K\"ahler form on the blowing-down of $\widetilde{X}^{\theta}$}
\par
Recall that 
$$
\pi \colon ( \widetilde{X}, \widetilde{X}^{\theta} ) \to ( X, {\rm Sing}\,X )
$$
is the blowing-down of the disjoint union of $(-2)$-curves $\widetilde{X}^{\theta} = \amalg_{{\mathfrak p} \in {\rm Sing}\,X} E_{\mathfrak p}$. 
Then ${\frak p}=\pi(E_{\frak p})$.
Under the identification $\psi_{\frak p}\colon(U_{\frak p},E_{\frak p})\cong(\widetilde{V}(1)_{\frak p},E)$,
$\pi\colon\widetilde{X}\to X$ is identified with the blowing-down $\varPi\colon T^{*}{\bf P}^{1}\to{\bf C}^{2}/\{\pm1\}$ on each $V(1)_{\frak p}$.
\par
By \cite{KobayashiTodorov87}, there exists a Ricci-flat orbifold K\"ahler form $\omega_{\eta}$ on $X$ such that
$$
\pi^{*}\omega_{\eta}^{2}/2!=\eta\wedge\overline{\eta}.
$$
By \eqref{eqn:Ricci:flatness}, we have
$$
\left.\pi^{*}\omega_{\eta}^{2}/\gamma_{\epsilon,\delta}^{2}\right|_{\widetilde{V}(\delta)_{\frak p}}=\varPi^{*}|f_{\frak p}(z)|^{2}.
$$
Since the right hand side is independent of $\epsilon\in(0,1)$, we get by putting $\epsilon\to0$
$$
\left.\omega_{\eta}^{2}/\gamma_{0}^{2}\right|_{\widetilde{V}(\delta)_{\frak p}}=|f_{\frak p}|^{2}.
$$
Hence we get the following relation by regarding $\eta$ as a nowhere vanishing holomorphic $2$-form on both $\widetilde{X}$ and $X$
$$
\left.\frac{\eta\wedge\overline{\eta}}{\gamma_{\epsilon,\delta}^{2}/2!}\right|_{E_{\frak p}}
=
|f_{\frak p}(0)|^{2}
=
\frac{\eta\wedge\overline{\eta}}{\gamma_{0}^{2}}({\frak p}).
$$

\section
{Behavior of some geometric quantities under the degeneration}
\par
In this section, we study the behavior of the second Chern form, the Bott-Chern term, and the analytic torsion of the fixed curves
when $\gamma_{\epsilon,\delta}$ converges to the orbifold metric $\gamma_{0}$.

\subsection
{\bf Behavior of the second Chern form as $\epsilon\to0$}
\par

\begin{proposition}
\label{prop:limit:2nd:Chern:form:1}
For any $\delta\in(0,1]$, one has
$$
\lim_{\epsilon\to0}\pi_{*}c_{2}(\widetilde{X},\gamma_{\epsilon,\delta})
=
c_{2}(X,\gamma_{0})+\frac{3}{2}\sum_{{\frak p}\in{\rm Sing}(X)}\delta_{\frak p}
$$
as currents on $X$, where $\delta_{\frak p}$ is the Dirac $\delta$-current supported at ${\frak p}$. In particular,
$$
\frac{1}{24}\int_{Y}c_{2}(Y,\gamma_{0})=\frac{1}{32}(16-k).
$$
\end{proposition}

\begin{pf}
Let $h\in C^{\infty}(X)$. By the definition of the K\"ahler form $\gamma_{\epsilon,\delta}$, we have
\begin{equation}
\label{eqn:5:2nd:Chern:form:1}
\begin{aligned}
\int_{\widetilde{X}}\pi^{*}h\cdot c_{2}(\widetilde{X},\gamma_{\epsilon,\delta})
&=
\int_{X\setminus\bigcup_{{\frak p}\in{\rm Sing}(X)}\widetilde{V}(\delta)_{\frak p}}h\cdot c_{2}(X,\gamma_{\epsilon,\delta})
+
\sum_{{\frak p}\in{\rm Sing}(X)}h({\frak p})\int_{\widetilde{V}(\delta)_{\frak p}}c_{2}(\widetilde{X},\gamma_{\epsilon,\delta})
\\
&\quad+
\sum_{{\frak p}\in{\rm Sing}(X)}\int_{\widetilde{V}(\delta)_{\frak p}}\pi^{*}\{h-h({\frak p})\}\cdot c_{2}(\widetilde{X},\gamma_{\epsilon,\delta}).
\end{aligned}
\end{equation}
\par
For $a>0$, let $T_{a}(z):=az$ be the homothety of ${\bf C}^{2}$ and
let $\widetilde{T}_{a}\colon T^{*}{\bf P}^{1}\to T^{*}{\bf P}^{1}$ be the biholomorphic map induced by $T_{a}$.
Then $\widetilde{T}_{\epsilon}$ induces an isometry of K\"ahler manifolds
$$
\widetilde{T}_{\epsilon}\colon(\widetilde{V}(\epsilon^{-2}),\epsilon^{2}\,\gamma^{\rm EH})
\cong
(\widetilde{V}(1),\gamma^{\rm EH}_{\epsilon^{2}}).
$$
Under the identification $T^{*}{\bf P}^{1}\setminus E\cong V(\infty)\setminus\{0\}$, we have the following estimates
$$
\left\|\gamma^{\rm EH}(z) - i\partial\bar{\partial}\|z\|^{2}\right\|\leq C(1+\|z\|)^{-4},
\qquad
\|c_{2}(T^{*}{\bf P}^{1},\gamma^{\rm EH})(z)\|\leq C(1+\|z\|)^{-6}
$$
for $\|z\|\gg1$ by \eqref{eqn:5:estimate1:error}, 
where $C>0$ is a constant and the norm is with respect to $\gamma^{\rm EH}$.
\par
Since there is a constant $C'>0$ with $\left|h|_{V(\delta)_{\frak p}}(z)-h({\frak p})\right|\leq C'\|z\|/(1+\|z\|)$ on $V(\delta)_{\frak p}$, 
we get
\begin{equation}
\label{eqn:5:estimate:2nd:Chern}
\begin{aligned}
\,&
\left|
\int_{\widetilde{V}(\delta)_{\frak p}}\{\pi^{*}h-h({\frak p})\}\cdot c_{2}(\widetilde{X},\gamma_{\epsilon,\delta})
\right|
=
\left|
\int_{\widetilde{V}(\delta)_{\frak p}}
\pi^{*}\{h|_{V(\delta)_{\frak p}}-h({\frak p})\}\cdot c_{2}(T^{*}{\bf P}^{1},\gamma^{\rm EH}_{\epsilon})
\right|
\\
&=
\left|
\int_{\widetilde{V}(\delta\sqrt{\epsilon}^{-1})}\widetilde{T}_{\epsilon}^{*}
\pi^{*}\{h|_{V(\delta)_{\frak p}}-h({\frak p})\}\cdot c_{2}(T^{*}{\bf P}^{1},\gamma^{\rm EH})
\right|
\\
&\leq
\int_{V(\delta\sqrt{\epsilon}^{-1})}C'\frac{\sqrt{\epsilon}\|z\|}{1+\sqrt{\epsilon}\|z\|}\cdot\frac{C}{1+\|z\|^{6}}\,\frac{(\gamma^{\rm EH})^{2}}{2!}
\leq
C''\sqrt{\epsilon}\to0\quad(\epsilon\to0),
\end{aligned}
\end{equation}
where $C''>0$ is a constant.
By \eqref{eqn:5:2nd:Chern:form:1} and \eqref{eqn:5:estimate:2nd:Chern}, we get
\begin{equation}
\label{eqn:5:2nd:Chern:form:2}
\begin{aligned}
\lim_{\epsilon\to0}\int_{\widetilde{X}}\pi^{*}h\cdot c_{2}(\widetilde{X},\gamma_{\epsilon,\delta})
&=
\int_{X\setminus\bigcup_{{\frak p}\in{\rm Sing}(X)}V(\delta)_{\frak p}}h\cdot c_{2}(X,\gamma_{0})
\\
&\quad
+
\sum_{{\frak p}\in{\rm Sing}(X)}h({\frak p})\int_{T^{*}{\bf P}^{1}}c_{2}(T^{*}{\bf P}^{1},\gamma^{\rm EH})
\\
&=
\int_{X}h\cdot c_{2}(X,\gamma_{0})
+
\sum_{{\frak p}\in{\rm Sing}(X)}h({\frak p})\int_{T^{*}{\bf P}^{1}}c_{2}(T^{*}{\bf P}^{1},\gamma^{\rm EH}),
\end{aligned}
\end{equation}
where we used the vanishing of $c_{2}(X,\gamma_{0})$ on $V(\delta)_{\frak p}$ to get the second equality.
Setting $h=1$ in \eqref{eqn:5:2nd:Chern:form:2} and comparing it with the formula \cite[p.396 l.5]{Kobayashi85}, we get
\begin{equation}
\label{eqn:5:GBC:instanton}
\int_{T^{*}{\bf P}^{1}}c_{2}(T^{*}{\bf P}^{1},\gamma^{\rm EH})=\chi({\bf P}^{1})-\frac{1}{|{\bf Z}_{2}|}=\frac{3}{2}.
\end{equation}
The first assertion follows from \eqref{eqn:5:2nd:Chern:form:2} and \eqref{eqn:5:GBC:instanton}.
\par
Since $\#{\rm Sing}(Y)=k$, we get by the first assertion
$$
2\int_{Y}c_{2}(Y,\gamma_{0})
=
\int_{X}c_{2}(X,\gamma_{0})
=
\int_{\widetilde{X}}c_{2}(\widetilde{X},\gamma_{\epsilon,\delta})-\frac{3}{2}\sum_{{\frak p}\in{\rm Sing}(X)}\int_{X}\delta_{\frak p}
=
24-\frac{3}{2}k.
$$
This proves the second assertion.
\end{pf}

\subsection
{\bf Behavior of the Bott-Chern terms as $\epsilon\to0$}
\par

\begin{proposition}
\label{prop:limit:2nd:Chern:form:2}
For any $\delta\in(0,1]$, one has
$$
\begin{aligned}
\,&
\lim_{\epsilon\to0}
\int_{\widetilde{X}}
\log
\left\{
\frac{\eta\wedge\overline{\eta}}{\gamma_{\epsilon,\delta}^{2}/2!}\cdot\frac{{\rm Vol}(\widetilde{X},\gamma_{\epsilon,\delta})}{\|\eta\|_{L^{2}}^{2}}
\right\}\,
c_{2}(\widetilde{X},\gamma_{\epsilon,\delta})
\\
&=
\int_{X}
\log\left\{\frac{\eta\wedge\overline{\eta}}{\gamma_{0}^{2}/2!}\cdot\frac{{\rm Vol}(X,\gamma_{0})}{\|\eta\|_{L^{2}}^{2}}\right\}\,c_{2}(X,\gamma_{0})
+
\frac{3}{2}\sum_{{\frak p}\in{\rm Sing}(X)}\log\left\{|f_{\frak p}(0)|^{2}\frac{{\rm Vol}(X,\gamma_{0})}{\|\eta\|_{L^{2}}^{2}}\right\}.
\end{aligned}
$$
\end{proposition}

\begin{pf}
Since $\gamma_{\epsilon,\delta}$ converges to $\gamma_{0}$ outside $\bigcup_{{\frak p}\in{\rm Sing}(X)}\widetilde{V}(\delta)_{\frak p}$
and since ${\rm Vol}(\widetilde{X},\gamma_{\epsilon,\delta})$ converges to ${\rm Vol}(X,\gamma_{0})$ as $\epsilon\to0$, 
we get the convergence
$$
\begin{aligned}
\,&
\left(
\int_{\widetilde{X}\setminus\bigcup_{{\frak p}\in{\rm Sing}(X)}\widetilde{V}(\delta)_{\frak p}}
+
\sum_{{\frak p}\in{\rm Sing}(X)}\int_{\widetilde{V}(\delta)_{\frak p}}
\right)
\log
\left\{
\frac{\eta\wedge\overline{\eta}}{\gamma_{\epsilon,\delta}^{2}/2!}\cdot\frac{{\rm Vol}(\widetilde{X},\gamma_{\epsilon,\delta})}{\|\eta\|_{L^{2}}^{2}}
\right\}
c_{2}(\widetilde{X},\gamma_{\epsilon,\delta})
\\
&\to
\int_{\widetilde{X}\setminus\bigcup_{{\frak p}\in{\rm Sing}(X)}\widetilde{V}(\delta)_{\frak p}}
\log
\left\{
\frac{\eta\wedge\overline{\eta}}{\gamma_{0}^{2}/2!}\cdot\frac{{\rm Vol}(X,\gamma_{0})}{\|\eta\|_{L^{2}}^{2}}
\right\}
c_{2}(\widetilde{X},\gamma_{0})
\\
&\quad+
\lim_{\epsilon\to0}
\sum_{{\frak p}\in{\rm Sing}(X)}\int_{\widetilde{V}(\delta)_{\frak p}}
\left\{
\log\pi^{*}|f_{\frak p}(z)|^{2}c_{2}(\widetilde{X},\gamma_{\epsilon,\delta})
+
\log\frac{{\rm Vol}(\widetilde{X},\gamma_{\epsilon,\delta})}{\|\eta\|_{L^{2}}^{2}}c_{2}(\widetilde{X},\gamma_{\epsilon,\delta})
\right\}
\\
&=
\int_{X\setminus\bigcup_{{\frak p}\in{\rm Sing}(X)}V(\delta)_{\frak p}}
\log
\left\{
\frac{\eta\wedge\overline{\eta}}{\gamma_{0}^{2}/2!}\cdot\frac{{\rm Vol}(X,\gamma_{0})}{\|\eta\|_{L^{2}}^{2}}
\right\}
c_{2}(X,\gamma_{0})
\\
&\quad+
\frac{3}{2}\sum_{{\frak p}\in{\rm Sing}(X)}\log
\left(|f_{\frak p}(0)|^{2}\frac{{\rm Vol}(X,\gamma_{0})}{\|\eta\|_{L^{2}}^{2}}\right)
\end{aligned}
$$
as $\epsilon\to0$, where the last equality follows from Proposition~\ref{prop:limit:2nd:Chern:form:1}.
Since $c_{2}(X,\gamma_{0})=0$ on $\bigcup_{{\frak p}\in{\rm Sing}(X)}V(\delta)_{\frak p}$, we get the result.
\end{pf}

\begin{corollary}
\label{cor:double:limit:torsion:K3}
For any $\delta\in(0,1]$, one has
$$
\begin{aligned}
\lim_{\delta\to0}\lim_{\epsilon\to0}\tau(\widetilde{X},\gamma_{\epsilon,\delta})
&=
\prod_{{\frak p}\in{\rm Sing}(X)}\left\{|f_{\frak p}(0)|^{2}\frac{{\rm Vol}(X,\gamma_{0})}{\|\eta\|_{L^{2}}^{2}}\right\}^{-\frac{1}{16}}
\\
&\quad\times
\exp
\left(
-\frac{1}{24}\int_{X}
\log\left\{\frac{\eta\wedge\overline{\eta}}{\gamma_{0}^{2}/2!}\cdot\frac{{\rm Vol}(X,\gamma_{0})}{\|\eta\|_{L^{2}}^{2}}\right\}\,
c_{2}(X,\gamma_{0})
\right).
\end{aligned}
$$
In particular, the limit $\lim_{\delta\to0}\lim_{\epsilon\to0}\tau(\widetilde{X},\gamma_{\epsilon,\delta})$ is independent of the choice of $\rho$.
\end{corollary}

\begin{pf}
By Theorem~\ref{thm:holomorphic:torsion:K3} and Proposition~\ref{prop:limit:2nd:Chern:form:2}, we get the desired equality.
The independence of the double limit $\lim_{\delta\to0}\lim_{\epsilon\to0}\tau(\widetilde{X},\gamma_{\epsilon,\delta})$ from $\rho$
is obvious, because the right hand side is independent of the choice of $\rho$.
\end{pf}

Define the Fubini-Study form on $E_{\frak p}$ by
$$
\omega_{\rm FS}(E_{\frak p}):=\left.\varPi^{*}\left(\frac{i}{2\pi}\partial\bar{\partial}\log\|z\|^{2}\right)\right|_{E_{\frak p}}
$$
Then for ${\frak p}\in{\rm Sing}(X)$, we have
$$
\gamma_{\epsilon,\delta}|_{E_{\frak p}}=\epsilon\,\omega_{\rm FS}(E_{\frak p})
$$
and an isomorphism of K\"ahler manifolds
$(E_{\frak p},\omega_{\rm FS}(E_{\frak p}))\cong({\bf P}^{1},\omega_{\rm FS})$.

\begin{proposition}
\label{prop:limit:anomaly}
For any $\delta\in(0,1]$, one has
$$
\lim_{\epsilon\to0}A_{M}(\widetilde{X},\theta,\gamma_{\epsilon,\delta})
=
\prod_{{\frak p}\in{\rm Sing}(X)}\left\{|f_{\frak p}(0)|^{2}\frac{{\rm Vol}(X,\gamma_{0})}{\|\eta\|_{L^{2}}^{2}}\right\}^{\frac{1}{4}}.
$$
\end{proposition}

\begin{pf}
Since $\gamma_{\epsilon,\delta}|_{E_{\frak p}}=\epsilon\,\omega_{\rm FS}(E_{\frak p})$ and 
since $\omega_{\rm FS}(E_{\frak p})$ is K\"ahler-Einstein, we get
$$
c_{1}(\widetilde{X}^{\theta},\gamma_{\epsilon,\delta}|_{\widetilde{X}^{\theta}})|_{E_{\frak p}}
=
\chi({\bf P}^{1})\,\omega_{\rm FS}(E_{\frak p})
=
2\,\omega_{\rm FS}(E_{\frak p}).
$$
Since $(\eta\wedge\overline{\eta})/(\gamma_{\epsilon,\delta}^{2}/2!)|_{E_{\frak p}}=|f_{\frak p}(0)|^{2}$
by \eqref{eqn:Ricci:flatness}, we get
$$
\begin{aligned}
\,&
A_{M}(\widetilde{X},\theta,\gamma_{\epsilon,\delta})
=
\exp
\left[
\frac{1}{8}\int_{\widetilde{X}^{\theta}}
\left.\log\left\{\frac{\eta\wedge\overline{\eta}}{\gamma_{\epsilon,\delta}^{2}/2!}
\cdot
\frac{{\rm Vol}(\widetilde{X},\gamma_{\epsilon,\delta})}{\|\eta\|_{L^{2}}^{2}}\right\}\right|_{\widetilde{X}^{\theta}}\,
c_{1}(\widetilde{X}^{\theta},\gamma_{\epsilon,\delta}|_{\widetilde{X}^{\theta}})
\right]
\\
&=
\exp
\left[
\frac{1}{4}\sum_{{\frak p}\in{\rm Sing}(X)}\int_{E_{\frak p}}
\left.\log\left\{
|f_{\frak p}(0)|^{2}\frac{{\rm Vol}(\widetilde{X},\gamma_{\epsilon,\delta})}{\|\eta\|_{L^{2}}^{2}}\right\}\right|_{E_{\frak p}}\,
\omega_{\rm FS}(E_{\frak p})
\right]
\\
&=
\exp
\left[
\frac{1}{4}\sum_{{\frak p}\in{\rm Sing}(X)}
\log\left\{
|f_{\frak p}(0)|^{2}\frac{{\rm Vol}(\widetilde{X},\gamma_{\epsilon,\delta})}{\|\eta\|_{L^{2}}^{2}}\right\}
\right]
\to
\prod_{{\frak p}\in{\rm Sing}(X)}\left\{|f_{\frak p}(0)|^{2}\frac{{\rm Vol}(X,\gamma_{0})}{\|\eta\|_{L^{2}}^{2}}\right\}^{\frac{1}{4}}
\end{aligned}
$$
as $\epsilon\to0$, where we used \eqref{eqn:continuity:volume} to get the last limit.
This completes the proof.
\end{pf}

\subsection
{\bf Behavior of the analytic torsion of the exceptional divisors}
\par

\begin{proposition}
\label{prop:torsion:fixed:curve}
For any $\delta\in(0,1]$ and ${\frak p}\in{\rm Sing}(X)$, the following equality holds for all $\epsilon\in(0,\delta^{2}\epsilon(\rho)]$
$$
\frac{{\rm Vol}(E_{\frak p},\gamma_{\epsilon,\delta}|_{E_{\frak p}})\tau(E_{\frak p},\gamma_{\epsilon,\delta}|_{E_{\frak p}})}
{{\rm Vol}({\bf P}^{1},\omega_{\rm FS})\tau({\bf P}^{1},\omega_{\rm FS})}
=
\epsilon^{1/3}.
$$
\end{proposition}

\begin{pf}
We recall a formula of Bost \cite[Prop.\,4.4]{Bost96}.
Let $(Z,g)$ be a compact K\"ahler manifold of dimension $d$ and let $\lambda>0$ be a constant.
By \cite[(4.2.4)]{Bost96}, we get
\begin{equation}
\label{eqn:5:formula:Bost}
\log\left(\frac{\tau(Z,\lambda g)}{\tau(Z,g)}\right)
=
\left(
-\sum_{i=0}^{d}(-1)^{i}(d-i)\,h^{0,i}(Z)+\int_{Z}{\rm Td}'(TZ)
\right)
\log\lambda,
\end{equation}
where the characteristic class ${\rm Td}'(E)$ is defined as follows (cf. \cite[Prop.\,4.4]{Bost96}).
If $\xi_{i}$ $(i=1,\ldots,r={\rm rk}(E))$ are the Chern roots of a vector bundle $E$, then
$$
{\rm Td}'(E)
:=
{\rm Td}(E)\cdot\sum_{i=1}^{r}\left(\frac{1}{\xi_{i}}-\frac{e^{-\xi_{i}}}{1-e^{-\xi_{i}}}\right).
$$
\par
Since
$$
{\rm Td}'(x)
=
\frac{x}{1-e^{-x}}\left(\frac{1}{x}-\frac{e^{-x}}{1-e^{-x}}\right)
=
\frac{1}{2}+\frac{1}{6}x+O(x^{2})
$$
and hence $\int_{{\bf P}^{1}}{\rm Td}'(T{\bf P}^{1})=1/3$, we get by \eqref{eqn:5:formula:Bost} applied to $(Z,g)=({\bf P}^{1},\omega_{\rm FS})$
\begin{equation}
\label{eqn:5:comparison:torsion:1}
\tau(E_{\frak p},\gamma_{\epsilon,\delta}|_{E_{\frak p}})/\tau({\bf P}^{1},\omega_{\rm FS})
=
\tau({\bf P}^{1},\epsilon\,\omega_{\rm FS})/\tau({\bf P}^{1},\omega_{\rm FS})
=
\epsilon^{-2/3}.
\end{equation}
Since
\begin{equation}
\label{eqn:5:comparison:volume}
{\rm Vol}(E_{\frak p},\gamma_{\epsilon,\delta}|_{E_{\frak p}})/{\rm Vol}({\bf P}^{1},\omega_{\rm FS})
=
{\rm Vol}({\bf P}^{1},\epsilon\,\omega_{\rm FS})/{\rm Vol}({\bf P}^{1},\omega_{\rm FS})
=
\epsilon,
\end{equation}
the result follows from \eqref{eqn:5:comparison:torsion:1} and \eqref{eqn:5:comparison:volume}.
\end{pf}

\section
{Spectrum and heat kernels under the degeneration}
\par
In this section, we prove a uniform lower bound of the $k$-th eigenvalue of the Laplacian 
and also a certain uniform exponential decay of the heat kernel for the degenerating family of metrics $\gamma_{\epsilon,\delta}$.

\subsection
{Uniformity of Sobolev inequality}
\par
In order to study the limit of the analytic torsions $\tau(\widetilde{X}, \gamma_{\epsilon,\delta})$ and 
$\tau_{{\mathbf Z}_{2}}(\widetilde{X}, \gamma_{\epsilon,\delta})(\theta)$ in the next section, we need to establish a uniform Sobolev inequality. 
First, we consider our model space  $(T^{*}{\bf P}^{1}, \gamma^{\rm EH}_{\epsilon})$, the Eguchi-Hanson instanton. 
Here $\gamma^{\rm EH}_{\epsilon}$ is the Ricci-flat K\"ahler metric constructed in Section 5.1 on $\widetilde{V}(\infty)=T^{*}{\bf P}^{1}$. 
Note that, for $0<\epsilon \leq 1$, under the identification $\Phi:\ ({\mathbf R}^4\!-\! B(\rho))/\{\pm 1\} \simeq \widetilde{V}(\infty)\! -\! K$ 
outside a compact neighborhood $K = \widetilde{V}(\rho) \subset \widetilde{V}(\infty)$ of the zero section of $T^{*}{\bf P}^{1}$ induced by
the identification $({\mathbf C}^2\!-\! B(\rho))/\{\pm 1\} = V(\infty) \!-\! V(\rho) = \widetilde{V}(\infty) \!-\! \widetilde{V}(\rho)$, one has 
$$
\Phi^*(\gamma^{\rm EH}_{\epsilon})_{ij} = \delta_{ij} + O(r^{-4})
$$ 
uniformly in $\epsilon$ by \eqref{eqn:5:estimate2:error}.

\begin{lemma} 
\label{lemma:Sobolev:ineq:1}
There is a constant $C$ such that for all  $0<\epsilon \leq 1$ the following holds.
\begin{itemize}
\item[(1)]  
For all $f\in C_0^{\infty}(\widetilde{V}(\infty))$,
$$  \| f \|_{L^4(\widetilde{V}(\infty),\gamma^{\rm EH}_{\epsilon})} \leq C \| df \|_{L^2(\widetilde{V}(\infty),\gamma^{\rm EH}_{\epsilon})}. $$
\item[(2)] 
Similarly, for all $\alpha \in A^{0,2}_0(\widetilde{V}(\infty))$,
$$  \| \alpha \|_{L^4(\widetilde{V}(\infty),\gamma^{\rm EH}_{\epsilon})} \leq C \| d\alpha \|_{L^2(\widetilde{V}(\infty),\gamma^{\rm EH}_{\epsilon})}. $$
\item[(3)] 
For all $\alpha \in A^{0,1}_0(\widetilde{V}(\infty))$,
$$  
\| \alpha \|^2_{L^4(\widetilde{V}(\infty),\gamma^{\rm EH}_{\epsilon})} 
\leq 
C^2 \left(\| \overline{\partial}\alpha \|^2_{L^2(\widetilde{V}(\infty),\gamma^{\rm EH}_{\epsilon})}
+
\| \overline{\partial}^*\alpha \|^2_{L^2(\widetilde{V}(\infty),\gamma^{\rm EH}_{\epsilon})}\right). 
$$
\end{itemize}
Here all norms are defined with respect to the metric $\gamma^{\rm EH}_{\epsilon}$.
\end{lemma}

\begin{pf}  
Since $(\widetilde{V}(\infty),\gamma^{\rm EH}_{\epsilon})\cong(\widetilde{V}(\infty),\epsilon\gamma^{\rm EH})$
by \eqref{eqn:scaling:property:EH} and since
the inequalities (1), (2), (3) above are invariant under the scaling of metrics $\gamma^{\rm EH}\mapsto\epsilon\gamma^{\rm EH}$,
it suffices to prove (1), (2), (3) for $\gamma^{\rm EH}$. In the rest to proof, all norms are defined with respect to $\gamma^{\rm EH}$.
Identifying a function in $C_0^{\infty}(\widetilde{V}(\infty)\! -\! K)$ with the corresponding $\pm 1$-invariant function on $\mathbf R^4$ 
with compact support via $\Phi$, we deduce from the Sobolev inequality for $\mathbf R^4$ that
$$ \|f\|_{L^4(\widetilde{V}(\infty))} \leq 2C\|df\|_{L^2(\widetilde{V}(\infty))}, \quad \forall f \in C_0^{\infty}(\widetilde{V}(\infty)\!-\! K),$$
where C is the Sobolev constant for $\mathbf R^4$. By an argument using partition of unity, there is a constant $C_K > 0$ such that
$$\|f\|_{L^4(\widetilde{V}(\infty))} \leq C_K(\|df\|_{L^2(\widetilde{V}(\infty))} + \|f\|_{L^2(K)}), 
\quad
\forall f \in C_0^{\infty}(\widetilde{V}(\infty)).$$
Assume that there is no constant $D > 0$ such that
$$ \|f\|_{L^2(K)}\leq D\|df\|_{L^2(\widetilde{V}(\infty))}, \quad \forall f \in C_0^{\infty}(\widetilde{V}(\infty)).$$
Then for any $n \in \mathbf N$, there is a function $f_n \in  C_0^{\infty}(\widetilde{V}(\infty))$ such that
$$\|f_n\|_{L^2(K)} = 1, \quad \|df_n\|_{L^2(\widetilde{V}(\infty))} \leq \frac{1}{n}.$$
Therefore, we have
$$ \|f_n\|_{L^4(\widetilde{V}(\infty))}\leq C_K(1 + 1/n) \leq 2C_K.$$
Passing to a subsequence if necessary, it follows that the sequence $f_n$ has a weak limit $f_{\infty} \in L^4(\widetilde{V}(\infty))$ 
with $df_{\infty}= 0$ as currents on $\widetilde{V}(\infty)$. This implies that in $L^4(\widetilde{V}(\infty))$, $f_{\infty} = 0$.
On the other hand, let $K'$ be a sufficiently big compact subset of $\widetilde{V}(\infty)$, whose open subset contains $K$. 
Now, for any compact subset $K' \subset \widetilde{V}(\infty)$, there is a constant $C_{K'} > 0$ such that
$$ \|f_n\|_{L^2(K')}\leq {\rm Vol}(K')^{1/2}\|f_n\|^{1/2}_{L^4(K')}\leq C_{K'} =\sqrt{2C_K{\rm Vol}(K')}.$$
Hence, by the Rellich lemma, we may assume (by passing to a subsequence if necessary again) 
that $f_n$ converges to $f_{\infty}$ strongly in $L^2(K')$. 
Since $K \subset K'$ and hence the convergence $f_n \rightarrow f_{\infty}$ in $L^2(K)$ is strong, 
we see that $ \|f_{\infty}\|_{L^2(K)} = \lim_{n \rightarrow \infty} \|f_n\|_{L^2(K)} = 1$.
This is a contradiction. Hence there exists a constant D such that $ \|f\|_{L^2(K)}\leq D\|df\|_{L^2(\widetilde{V}(\infty))}$.  
By setting $C = CK(1 + D$), we have
$$  \| f \|_{L^4(\widetilde{V}(\infty))} \leq C \| df \|_{L^2(\widetilde{V}(\infty))}. $$
This proves (1).

(2) is an immediate consequence of (1) and the isomorphism 
$C_0^{\infty}(\widetilde{V}(\infty))\ni f\mapsto f\overline{\eta}\in A^{0,2}_{0}(\widetilde{V}(\infty))$, 
which commutes with the operations involved. To see (3), let $\alpha \in A^{0,1}_0(\widetilde{V}(\infty))$. 
Then, by (1)
$$  
\| \alpha \|^2_{L^4(\widetilde{V}(\infty))} 
=
 \left( \int_{\widetilde{V}(\infty)} |\alpha|^4 dx\right)^{1/2} \leq C^2 \int_{\widetilde{V}(\infty)} |\, d|\alpha|\, |^2dx . 
$$
Using Kato's inequality, we have
$$ \int_{\widetilde{V}(\infty)} |\, d|\alpha|\, |^2dx \leq \int_{\widetilde{V}(\infty)} |\nabla\alpha|^2dx. $$
Now the Bochner formula \cite[(1.4.63)]{MM} gives 
$(\overline{\partial}\overline{\partial}^* + \overline{\partial}^*\overline{\partial})\alpha = \nabla^*\nabla \alpha$ 
since $(\widetilde{V}(\infty),\gamma^{\rm EH}_{\epsilon})$ is Ricci flat. Our result follows.
\end{pf}

\begin{lemma}
\label{lemma:Sobolev:ineq:2}
There is a constant $C$ such that for all  $\epsilon,\delta\in(0,1]$ with $\epsilon\delta^{-2}\leq\epsilon(\rho)$, 
and all $\alpha \in A^{0,q}_0(\widetilde{V}(\infty))$, $0\leq q\leq 2$,
$$  
\| \alpha \|^2_{L^4(\widetilde{V}(\infty),\kappa_{\epsilon,\delta})} 
\leq 
C^2 
\left(
\| \overline{\partial}\alpha \|^2_{L^{2}(\widetilde{V}(\infty),\kappa_{\epsilon,\delta})}
+
\| \overline{\partial}^*\alpha \|^2_{L^2(\widetilde{V}(\infty),\kappa_{\epsilon,\delta})}
\right),
$$
where the norms and $\bar{\partial}^{*}$ are defined with respect to the metric $\kappa_{\epsilon,\delta}$.
\end{lemma}

\begin{pf}
By Lemmas~\ref{lemma:comparison:1} and \ref{lemma:comparison:2}, there exist constants $C_{1},C_{2}>0$ such that
\begin{equation}
\label{eqn:quasi:isometry:1}
C_{1}\gamma^{\rm EH}_{\epsilon}\leq\kappa_{\epsilon,\delta}\leq C_{2}\gamma^{\rm EH}_{\epsilon}.
\end{equation}
for all $\epsilon,\delta\in(0,1]$ with $\epsilon\delta^{-2}\leq\epsilon(\rho)$.
Hence there is a constant $C_{3}>0$ such that
\begin{equation}
\label{eqn:quasi:isometry:2}
C_{3}^{-1}\| \alpha \|^2_{L^4(\widetilde{V}(\infty),\gamma^{\rm EH}_{\epsilon})} 
\leq
\| \alpha \|^2_{L^4(\widetilde{V}(\infty),\kappa_{\epsilon,\delta})} 
\leq
C_{3}\| \alpha \|^2_{L^4(\widetilde{V}(\infty),\gamma^{\rm EH}_{\epsilon})},
\end{equation}
\begin{equation}
\label{eqn:quasi:isometry:3}
C_{3}^{-1}\| \overline{\partial}\alpha \|^2_{L^{2}(\widetilde{V}(\infty),\gamma^{\rm EH}_{\epsilon})}
\leq
\| \overline{\partial}\alpha \|^2_{L^{2}(\widetilde{V}(\infty),\kappa_{\epsilon,\delta})}
\leq
C_{3}\| \overline{\partial}\alpha \|^2_{L^{2}(\widetilde{V}(\infty),\gamma^{\rm EH}_{\epsilon})}
\end{equation}
for all $\epsilon,\delta\in(0,1]$ with $\epsilon\delta^{-2}\leq\epsilon(\rho)$ and $\alpha\in A^{0,q}_{0}(\widetilde{V}(\infty))$.
\par
Let $\Lambda_{\epsilon.\delta}$ (resp. $\Lambda_{\epsilon}$) be the Lefschetz operator defined as the adjoint of the multiplication 
by $\kappa_{\epsilon,\delta}$ (resp. $\gamma^{\rm EH}_{\epsilon}$).
Since $\overline{\partial}^{*}=\pm i\Lambda_{\epsilon,\delta}\partial$ for $(0,q)$-forms by the K\"ahler identity, 
there exists by \eqref{eqn:quasi:isometry:1} a constant $C_{4}>0$ such that
\begin{equation}
\label{eqn:quasi:isometry:4}
C_{4}^{-1}\| \overline{\partial}^*\alpha \|^2_{L^2(\widetilde{V}(\infty),\gamma^{\rm EH}_{\epsilon})}
\leq
\| \overline{\partial}^*\alpha \|^2_{L^2(\widetilde{V}(\infty),\kappa_{\epsilon,\delta})}
\leq
C_{4}\| \overline{\partial}^*\alpha \|^2_{L^2(\widetilde{V}(\infty),\gamma^{\rm EH}_{\epsilon})}.
\end{equation}
By Lemma~\ref{lemma:Sobolev:ineq:1} (3) and \eqref{eqn:quasi:isometry:2}, \eqref{eqn:quasi:isometry:3}, \eqref{eqn:quasi:isometry:4},
we get the result.
\end{pf}

For the minimal resolution $\widetilde{X}$ and the family of K\"ahler metrics $\gamma_{\epsilon,\delta}$ 
constructed in Section 5.2 using the Eguchi-Hanson instanton, we have

\begin{proposition} 
\label{prop:Sobolev:ineq:3}
There is a constant $C$ such that for all  $\epsilon,\delta\in(0,1]$ with $\epsilon\delta^{-2}\leq\epsilon(\rho)$, 
and all $\alpha \in A^{0,q}(\widetilde{X})$, $0\leq q\leq 2$,
$$  
\| \alpha \|^2_{L^4(\widetilde{X},\gamma_{\epsilon,\delta})} 
\leq 
C^2 
\left(
\| \overline{\partial}\alpha \|^2_{L^2(\widetilde{X},\gamma_{\epsilon,\delta})}
+
\| \overline{\partial}^*\alpha \|^2_{L^2(\widetilde{X},\gamma_{\epsilon,\delta})} 
+
 \|\alpha \|^2_{L^2(\widetilde{X},\gamma_{\epsilon,\delta})}  
\right),
$$
where the norms are defined with respect to the metric $\gamma_{\epsilon,\delta}$.
\end{proposition}

\begin{pf}
Since $\gamma_{\epsilon,\delta}=\kappa_{\epsilon,\delta}$ on $\bigcup_{{\frak p}\in{\rm Sing}(X)}\widetilde{V}(\delta)_{\frak p}$,
the result follows from Lemma~\ref{lemma:Sobolev:ineq:2} and an easy partition of unity argument.
\end{pf}

\subsection
{A uniform lower bound of spectrum}
\par
Let $\square_{\epsilon,\delta}^{q}=(\bar{\partial}+\bar{\partial}^{*})^{2}$ (resp. $\square_{0}^{q}$)
be the Hodeg-Kodaira Laplacian of $(\widetilde{X},\gamma_{\epsilon,\delta})$ (resp. $(X,\gamma_{0})$) acting on $(0,q)$-forms.
Let $\lambda_{\epsilon,\delta}^{q}(k)$ (resp. $\lambda_{0}^{q}(k)$) be the $k$-th {\em non-zero} eigenvalue of the Laplacian 
$\square_{\epsilon,\delta}^{q}$ (resp. $\square_{0}^{q}$). Then the non-zero eigenvalues of $\square_{\epsilon,\delta}^{q}$ are given by
$$
0<\lambda_{\epsilon,\delta}^{q}(1)\leq\lambda_{\epsilon,\delta}^{q}(2)
\leq\cdots\leq\lambda_{\epsilon,\delta}^{q}(k)\leq\lambda_{\epsilon,\delta}^{q}(k+1)\leq\cdots
$$
and the set of corresponding eigenforms $\{\varphi_{k,\epsilon,\delta}^{q}\}_{k\in{\bf N}}$. We set $\lambda_{\epsilon,\delta}^{q}(0)=0$ and 
list the corresponding eigenforms $\varphi_{0,\epsilon,\delta}^{q}$ (here we abuse the notation as there would be 
$\dim H^0(\tilde X, \Omega^q_{\tilde X})$ many of them) so that $\{\varphi_{k,\epsilon,\delta}^{q}\}_{k=0}^{\infty}$ forms a complete 
orthonormal basis of $L^{0,q}_{\epsilon,\delta}(\widetilde{X})$,
the $L^{2}$-completion of $A^{0,q}(\widetilde{X})$ with respect to the norm associated to $\gamma_{\epsilon,\delta}$. Since
$$
K_{\epsilon,\delta}^{q}(t,x,y)=\sum_{k=0}^{\infty}e^{-t\lambda_{\epsilon,\delta}^{q}(k)}\varphi_{k,\epsilon,\delta}^{q}(x)\otimes\varphi_{k,\epsilon,\delta}^{q}(y)^{*},
$$
we get 
\begin{equation}
\label{eqn:CauchySchwartz}
\begin{aligned}
\left|
K_{\epsilon,\delta}^{q}(t,x,y)
\right|
&\leq
\sum_{k=0}^{\infty}e^{-t\lambda_{\epsilon,\delta}^{q}(k)}|\varphi_{k,\epsilon,\delta}^{q}(x)|\cdot|\varphi_{k,\epsilon,\delta}^{q}(y)|
\\
&\leq
\{\sum_{k=0}^{\infty}e^{-t\lambda_{\epsilon,\delta}^{q}(k)}|\varphi_{k,\epsilon,\delta}^{q}(x)|^{2}\}^{1/2}
\{\sum_{k=0}^{\infty}e^{-t\lambda_{\epsilon,\delta}^{q}(k)}|\varphi_{k,\epsilon,\delta}^{q}(y)|^{2}\}^{1/2}
\\
&=
\sqrt{{\rm tr}\,K_{\epsilon,\delta}^{q}(t,x,x)}\sqrt{{\rm tr}\,K_{\epsilon,\delta}^{q}(t,y,y)}.
\end{aligned}
\end{equation}

\begin{proposition}
\label{prop:heat:trace}
If $q=0,2$, then there are 
constants $A, C>0$ such that for all $\epsilon,\delta\in(0,1]$ with $\epsilon\delta^{-2}\leq\epsilon(\rho)$, 
and $x,y\in\widetilde{X}$, $t>0$, the following inequality holds:
\begin{equation}
\label{eqn:heat:kernel:bound:1}
0<
\left| K_{\epsilon,\delta}^{q}(t,x,y) \right|
\leq 
A e^{C(\epsilon\delta^{-4}+1)}(t^{-2}+1).
\end{equation}
Moreover, for all $(\epsilon,\delta)\in(0,1]$ with $\epsilon\delta^{-2}\leq\epsilon(\rho)$ 
and for all $t>0$, $q\geq0$, the following inequality holds:
\begin{equation}
\label{eqn:heat:kernel:bound:0}
{\rm Tr}\,e^{-t\,\square_{\epsilon,\delta}^{q}}\leq {\rm Vol}(\widetilde{X},\gamma_{\epsilon,\delta})Ae^{C(\epsilon\delta^{-4}+1)}(t^{-2}+1).
\end{equation}
\end{proposition}

\begin{pf}
{\em (Case 1) }
Let $q=0$.
By Proposition~\ref{prop:Sobolev:ineq:3}, the Sobolev constant is uniform for $\epsilon,\delta\in(0,1]$ with $\epsilon\delta^{-2}\leq\epsilon(\rho)$.
By \cite[Thms.\,2.1 and 2.16]{CarlenKusuokaStroock87}, there are constants $A>0$, $B\geq0$ such that for all $\epsilon,\delta\in(0,1]$ 
with $\epsilon\delta^{-2}\leq\epsilon(\rho)$, and $x,y\in\widetilde{X}$, $t>0$,
\begin{equation}
\label{eqn:heat:kernel:bound:weak}
0<
K_{\epsilon,\delta}^{q}(t,x,y)
\leq 
A\,e^{Bt}t^{-2}.
\end{equation}
Let $q=2$.  By Lemma~\ref{lemma:Ricci}, the Lichnerowicz formula and \cite[p.32 l.4-l.5]{HessSchraderUhlenbrock80}, we have
\begin{equation}
\label{eqn:heat:kernel:bound:weak:q=2}
| K_{\epsilon,\delta}^{q}(t,x,y) | \leq e^{t|{\rm Ric}\gamma_{\epsilon,\delta}|_{\infty}} | K_{\epsilon,\delta}^{0}(t,x,y) | \leq e^{C(\epsilon\delta^{-4}+1)t}Ae^{Bt}(t^{-2}+1).
\end{equation}
For $t\leq1$, we get \eqref{eqn:heat:kernel:bound:1} by \eqref{eqn:heat:kernel:bound:weak}, \eqref{eqn:heat:kernel:bound:weak:q=2}.
For $t\geq1$, since ${\rm tr}\,K_{\epsilon,\delta}^{q}(t,x,x)$ is a decreasing
function in $t$, we deduce \eqref{eqn:heat:kernel:bound:1} from \eqref{eqn:CauchySchwartz}, \eqref{eqn:heat:kernel:bound:weak}, \eqref{eqn:heat:kernel:bound:weak:q=2}
and the inequality
$$
\left|
K_{\epsilon,\delta}^{q}(t,x,y)
\right|
\leq
\sqrt{{\rm tr}\,K_{\epsilon,\delta}^{q}(1,x,x)}\sqrt{{\rm tr}\,K_{\epsilon,\delta}^{q}(1,y,y)}
\leq
2e^{C(\epsilon\delta^{-4}+1)}Ae^{B}.
$$
Since ${\rm Tr}\,e^{-t\,\square_{\epsilon,\delta}^{q}}=\int_{\widetilde{X}}{\rm tr}\,K_{\epsilon,\delta}^{q}(t,x,x)\,dx$, we get \eqref{eqn:heat:kernel:bound:0} 
from \eqref{eqn:heat:kernel:bound:1}.
\par{\em (Case 2) }
Let $q=1$. Since $\sum_{q}(-1)^{q}{\rm Tr}\,e^{-t\,\square_{\epsilon,\delta}^{q}}=0$ for all $t>0$, \eqref{eqn:heat:kernel:bound:0} for $q=1$ follows from
\eqref{eqn:heat:kernel:bound:0} for $q=0,2$. This completes the proof.
\end{pf}

Write $\lambda_{\epsilon}^{q}(k)$ for $\lambda_{\epsilon,1}^{q}(k)$.

\begin{lemma}
\label{lemma:bound:1st:eigenvalu:1}
There is a constant $\lambda>0$ such that for all $\epsilon\in(0,\epsilon(\rho)]$ and $q\geq0$,
\begin{equation}
\label{eqn:lower:bound:1st:eigenvalue}
\lambda_{\epsilon}^{q}(1)\geq\lambda>0.
\end{equation}
\end{lemma}

\begin{pf}
Since $\dim\widetilde{X}=2$ and hence $\lambda_{\epsilon}^{1}(1)=\lambda_{\epsilon}^{0}(1)$ or 
$\lambda_{\epsilon}^{1}(1)=\lambda_{\epsilon}^{2}(1)$, it suffices to prove \eqref{eqn:lower:bound:1st:eigenvalue} for $q=0,2$.
Assume that there is a sequence $\{\epsilon_{n}\}$ such that $\epsilon_{n}\to0$ and $\lambda_{\epsilon_{n}}^{q}(1)\to0$ as $n\to\infty$ for $q=0$ or $2$.
By the same argument as in \cite[p.434--p.436]{Yoshikawa96} using the uniformity of the Sobolev constant (cf. Proposition~\ref{prop:Sobolev:ineq:3}),
there is a holomorphic $q$-form $\psi$ on $X \setminus {\rm Sing}\,X$, which is possibly meromorphic on $\widetilde{X}$, with the following properties:
\begin{itemize}
\item[(i)]
The complex conjugation $\overline{\varphi_{1,\epsilon_{n}}^{q}}$ converges to $\psi$ 
on every compact subset of $X\setminus{\rm Sing}(X)$ as $n\to\infty$. 
\item[(ii)]
$\|\psi\|_{L^{2}}=1$ and $\pi^{*}\psi\perp H^{0}(\widetilde{X},\Omega^{q}_{\widetilde{X}})$ with respect to the degenerate K\"ahler metric $\pi^{*}\gamma_{0}$ on $\widetilde{X}$.
\end{itemize}
Since ${\rm Sing}\,X$ consists of isolated orbifold points, it follows from the Riemann extension theorem that $\psi$ extends to a holomorphic $q$-form on $X$
in the sense of orbifolds. When $q=0$, $\psi$ is a constant. When $q=2$, since $X$ has canonical singularities,
$\pi^{*}\psi$ is a holomorphic $2$-form on $\widetilde{X}$. In both cases,
the condition $\pi^{*}\psi\perp H^{0}(\widetilde{X},\Omega^{q}_{\widetilde{X}})$ implies $\psi =0$, which contradicts the other condition $\|\psi\|_{L^{2}}=1$.
This proves the result.
\end{pf}

\begin{lemma}
\label{lemma:bound:1st:eigenvalu:2}
There is a constant $\lambda'>0$ such that for all $\epsilon,\delta\in(0,1]$ with $\epsilon\delta^{-2}\leq\epsilon(\rho)]$ and $q\geq0$,
$$
\lambda_{\epsilon,\delta}^{q}(1)\geq\lambda'>0.
$$
\end{lemma}

\begin{pf}
Firstly we prove the inequality when $q=1$. Since $\widetilde{X}$ is a $K3$ surface and hence $\ker\square_{\epsilon,\delta}^{1}=0$,
we get by \eqref{eqn:lower:bound:1st:eigenvalue} 
\begin{equation}
\label{eqn:L2:estimate:1}
\lambda\,\|\alpha\|_{L^{2}(\widetilde{X},\gamma_{\epsilon})}^{2}
\leq
\|\bar{\partial}\alpha\|_{L^{2}(\widetilde{X},\gamma_{\epsilon})}^{2}+\|\bar{\partial}^{*}\alpha\|_{L^{2}(\widetilde{X},\gamma_{\epsilon})}^{2}
=
\|\partial\alpha\|_{L^{2}(\widetilde{X},\gamma_{\epsilon})}^{2}
\end{equation}
for all $\alpha\in A^{0,1}(\widetilde{X})$, where we used the coincidence of the $\bar{\partial}$-Laplacian and the $\partial$-Laplacian for K\"ahler manifolds 
to get the equality in \eqref{eqn:L2:estimate:1}. By Lemma~\ref{lemma:comparison:3}, there exist constants $C_{1}>0$ such that
for all $\alpha\in A^{0,1}(\widetilde{X})$,
$$
C_{1}^{-1}\|\alpha\|_{L^{2}(\widetilde{X},\gamma_{\epsilon})}^{2}
\leq
\|\alpha\|_{L^{2}(\widetilde{X},\gamma_{\epsilon,\delta})}^{2}
\leq
C_{1}\|\alpha\|_{L^{2}(\widetilde{X},\gamma_{\epsilon})}^{2},
$$
$$
C_{1}^{-1}\|\partial\alpha\|_{L^{2}(\widetilde{X},\gamma_{\epsilon})}^{2}
\leq
\|\partial\alpha\|_{L^{2}(\widetilde{X},\gamma_{\epsilon,\delta})}^{2}
\leq
C_{1}\|\partial\alpha\|_{L^{2}(\widetilde{X},\gamma_{\epsilon})}^{2}.
$$
Combining these inequalities and \eqref{eqn:L2:estimate:1}, we get for all $\alpha\in A^{0,1}(\widetilde{X})$
\begin{equation}
\label{eqn:L2:estimate:2}
C_{1}^{-1}\lambda\,\|\alpha\|_{L^{2}(\widetilde{X},\gamma_{\epsilon,\delta})}^{2}
\leq
C_{1}\|\partial\alpha\|_{L^{2}(\widetilde{X},\gamma_{\epsilon,\delta})}^{2}
=
C_{1}
\left(
\|\bar{\partial}\alpha\|_{L^{2}(\widetilde{X},\gamma_{\epsilon,\delta})}^{2}+\|\bar{\partial}^{*}\alpha\|_{L^{2}(\widetilde{X},\gamma_{\epsilon,\delta})}^{2}
\right).
\end{equation}
The result for $q=1$ follows from \eqref{eqn:L2:estimate:2}.
Since $\bar{\partial}\varphi_{\epsilon,\delta}^{0}(1)$ and $\bar{\partial}^{*}\varphi_{\epsilon,\delta}^{2}(1)$ are {\em non-zero} eigenforms of $\square_{\epsilon,\delta}^{1}$,
we get $\lambda'\leq\lambda_{\epsilon,\delta}^{1}(1)\leq\lambda_{\epsilon,\delta}^{0}(1)$ and 
$\lambda'\leq\lambda_{\epsilon,\delta}^{1}(1)\leq\lambda_{\epsilon,\delta}^{2}(1)$.
\end{pf}

\begin{theorem}
\label{thm:bound:eigenvalue}
There are 
constants $\Lambda, C>0$ such that for all $k\in{\bf N}$, $\epsilon,\delta\in(0,1]$ with $\epsilon\delta^{-2}\leq\epsilon(\rho)$ and $q\geq0$,
$$
\lambda_{\epsilon,\delta}^{q}(k)\geq\Lambda e^{-\frac{1}{2}C(\epsilon\delta^{-4}+1)}\,k^{1/2}.
$$
\end{theorem}

\begin{pf}
By Proposition~\ref{prop:heat:trace}, we get for all $\epsilon,\delta\in(0,1]$ with $\epsilon\delta^{-2}\leq\epsilon(\rho)$ and $t\in(0,1]$
$$
\sum_{i=1}^{k}e^{-t\lambda_{\epsilon,\delta}^{q}(i)}
\leq
h^{0,q}(\widetilde{X})+\sum_{i=1}^{\infty}e^{-t\lambda_{\epsilon,\delta}^{q}(i)}
=
{\rm Tr}\,e^{-t\square_{\epsilon,\delta}^{q}}
\leq
A'e^{C(\epsilon\delta^{-4}+1)} \,t^{-2},
$$
where $A'$ is a constant such that $A{\rm Vol}(\widetilde{X},\gamma_{\epsilon,\delta}) \leq A'$.
Since $\lambda'/\lambda_{\epsilon,\delta}^{q}(k)\leq1$ by Lemma~\ref{lemma:bound:1st:eigenvalu:2}, substituting $t:=\lambda'/\lambda_{\epsilon,\delta}^{q}(k)$
in this inequality and using $\lambda_{\epsilon,\delta}^{q}(i)/\lambda_{\epsilon,\delta}^{q}(k)\leq1$ for $i\leq k$, we get
$$
k\,e^{-\lambda'}
\leq
\sum_{i=1}^{k}e^{-\frac{\lambda'\lambda_{\epsilon,\delta}^{q}(i)}{\lambda_{\epsilon,\delta}^{q}(k)}}
\leq
A'e^{C(\epsilon\delta^{-4}+1)}\,\left(\frac{\lambda'}{\lambda_{\epsilon,\delta}^{q}(k)}\right)^{-2}.
$$
We get the result by setting $\Lambda:=(A')^{-1/2}\lambda' e^{-\lambda'/2}$.
\end{pf}

\begin{corollary}
\label{cor:uniform:exp:decay}
Let $C$ and $\Lambda$ be the same constants as in Theorem~\ref{thm:bound:eigenvalue} and
set $\Lambda(R):=\Lambda e^{-\frac{1}{2}CR}$ and $\Psi(R):=\sum_{k=1}^{\infty}e^{-\frac{1}{2}\Lambda(R) k^{1/2}}$. 
Then, for all $\epsilon,\delta\in(0,1]$ with $\epsilon\delta^{-2}\leq\epsilon(\rho)$ and $t\geq1$, 
the following inequality holds
$$
0
<
{\rm Tr}\,e^{-t\square_{\epsilon,\delta}^{q}}-h^{0,q}(\widetilde{X})
\leq 
\Psi(\epsilon\delta^{-4}+1)\,e^{-\frac{1}{2}\Lambda (1+\epsilon\delta^{-4})t}.
$$
\end{corollary}

\begin{pf}
Since $\lambda_{\epsilon,\delta}^{q}(k)\geq\frac{\Lambda(\epsilon\delta^{-4}+1)}{2}(k^{1/2}+1)$ by Theorem~\ref{thm:bound:eigenvalue}, we get 
$\sum_{k=1}^{\infty}e^{-t\lambda_{\epsilon,\delta}^{q}(k)} \leq 
e^{-t\Lambda(\epsilon\delta^{-4}+1)/2}\sum_{k=1}^{\infty}e^{-\frac{1}{2}t\Lambda(\epsilon\delta^{-4}+1) k^{1/2}} 
\leq \Psi(\epsilon\delta^{-4}+1)\,e^{-t\Lambda(\epsilon\delta^{-4}+1)/2}$
for $t\geq1$.
\end{pf}

We also need an estimate for the heat kernel $K_{\epsilon,\delta,\infty}^{q}(t,x,y)$ of $(\widetilde{V}(\infty),\kappa_{\epsilon,\delta})$.

\begin{proposition}
\label{prop:heat:trace:noncompact}
There are  
constants $A', C'>0$ such that for all $\epsilon,\delta\in(0,1]$ with $\epsilon\delta^{-2}\leq\epsilon(\rho)$, $x\in\widetilde{V}(\infty)$, 
$t>0$ and $q\geq0$, the following inequality holds:
$$
\left|K_{\epsilon,\delta,\infty}^{q}(t,x,y)\right|\leq A'e^{C'(\epsilon\delta^{-4}+1)} (t^{-2}+1).
$$
\end{proposition}

\begin{pf}
When $q=0$, the result follows from Lemma~\ref{lemma:Sobolev:ineq:2} and \cite[Thms.\,2.1 and 2.16]{CarlenKusuokaStroock87}.
Let $q>0$. Since $|{\rm Ric}(\gamma_{\epsilon,\delta})|\leq C(\epsilon\delta^{-4}+1)$ by Lemma~\ref{lemma:Ricci}, 
we deduce from \cite[p.32 l.4-l.5]{HessSchraderUhlenbrock80} and the Lichnerowicz formula for $\square_{\epsilon,\delta}^{q}$ that
$$
0<
\left|
K_{\epsilon,\delta,\infty}^{q}(t,x,y)
\right|
\leq
e^{C(\epsilon\delta^{-4}+1)t}K_{\epsilon,\delta,\infty}^{0}(t,x,y)
\leq 
A\,e^{C(\epsilon\delta^{-4}+1)t}t^{-2}.
$$
This proves the result for $t\leq1$. 
Since \eqref{eqn:CauchySchwartz} remains valid for $K_{\epsilon,\delta,\infty}^{q}(t,x,y)$ by the fact that $K_{\epsilon,\delta,\infty}^{q}(t,x,y)$
is obtained as the limit $R\to\infty$ of the Dirichlet heat kernel of $\widetilde{V}(R)$, 
the result for $t\geq1$ also follows.
This completes the proof.
\end{pf}

\section
{Behavior of (equivariant) analytic torsion}
\par
In the previous sections, the additional parameter $\delta$ is pretty harmless and the results still hold in its presence. 
This parameter will play more essential role in this section. Indeed, we shall prove the following:

\begin{theorem} \label{loat}
There exist constants $C_{0}(k),C_{1}(k)>0$ depending only on $k=\#{\rm Sing}(Y)$ such that
$$
\lim_{\delta\to0}\lim_{\epsilon\to0}\tau(\widetilde{X},\gamma_{\epsilon, \delta})=C_{0}(k)\cdot\tau(X,\gamma_{0}),
$$
$$
\lim_{\delta\to0}\lim_{\epsilon\to0}\epsilon^{k/3}\tau_{{\bf Z}_{2}}(\widetilde{X},\gamma_{\epsilon, \delta})(\theta)=C_{1}(k)\cdot\tau_{{\bf Z}_{2}}(X,\gamma_{0})(\iota).
$$
\end{theorem}

\subsection
{Existence of limit}
\par
By Corollary~\ref{cor:double:limit:torsion:K3}, the first limit exists and is independent of the choice of a cut-off function $\rho$.
For the second limit we have

\begin{proposition}
\label{prop:limit:equivariant:torsion}
For any $\delta\in(0,1]$, the number
$$
\epsilon^{k/3}\tau_{{\bf Z}_{2}}(\widetilde{X},\gamma_{\epsilon,\delta})(\theta){\rm Vol}(\widetilde{X},\gamma_{\epsilon,\delta})
$$
is independent of $\epsilon,\delta\in(0,1]$ with $0<\epsilon\delta^{-2}\leq\epsilon(\rho)$. 
In particular, for any $\delta\in(0,1]$, the following limit exists as $\epsilon\to0$:
$$
\lim_{\epsilon\to0}\epsilon^{k/3}\tau_{{\bf Z}_{2}}(\widetilde{X},\gamma_{\epsilon,\delta})(\theta)
$$
and the limit is independent of $\delta\in(0,1]$ and the choice of a cut-off function $\rho$.
\end{proposition}

\begin{pf}
{\em (Step 1) }
Let $g_{0},g_{1}$ be $\theta$-invariant K\"ahler metrics on $\widetilde{X}$.
Let $\widetilde{\rm Td}_{\theta}(T\widetilde{X};g_{0},g_{1})^{(1,1)}$ be the Bott-Chern class such that
$$
-dd^{c}\widetilde{\rm Td}_{\theta}(T\widetilde{X};g_{0},g_{1})
=
{\rm Td}_{\theta}(T\widetilde{X},g_{0})
-
{\rm Td}_{\theta}(T\widetilde{X},g_{1}).
$$
By Bismut \cite[Th.\,2.5]{Bismut95},
\begin{equation}
\label{eqn:anomaly:equivariant:torsion}
\log
\left(
\frac{\tau_{{\bf Z}_{2}}(\widetilde{X},g_{0})(\theta){\rm Vol}(\widetilde{X},g_{0})}{\tau_{{\bf Z}_{2}}(\widetilde{X},g_{1})(\theta){\rm Vol}(\widetilde{X},g_{1})}
\right)
=
\int_{\widetilde{X}^{\theta}}
\widetilde{\rm Td}_{\theta}(T\widetilde{X};g_{0},g_{1}).
\end{equation}
Since
$$
{\rm Td}_{\theta}(T\widetilde{X};g_{0},g_{1})^{(1,1)}
=
\frac{1}{8}
{c_{1}(T\widetilde{X})|_{\widetilde{X}^{\theta}}c_{1}(T\widetilde{X}^{\theta})}(g_{0},g_{1})
-
\frac{1}{12}
{c_{1}(T\widetilde{X}^{\theta})^{2}}(g_{0},g_{1})
$$
by \cite[Prop.\,5.3]{Yoshikawa04}, we have the following equality of Bott-Chern classes:
$$
\begin{aligned}
\widetilde{\rm Td}_{\theta}(T\widetilde{X};g_{0},g_{1})^{(1,1)}
&=
\frac{1}{8}
\widetilde{c_{1}(T\widetilde{X})|_{\widetilde{X}^{\theta}}c_{1}(T\widetilde{X}^{\theta})}(g_{0},g_{1})
-
\frac{1}{12}
\widetilde{c_{1}(T\widetilde{X}^{\theta})^{2}}(g_{0},g_{1})
\\
&=
\frac{1}{8}
\widetilde{c_{1}}(T\widetilde{X};g_{0},g_{1})|_{\widetilde{X}^{\theta}}c_{1}(T\widetilde{X}^{\theta},g_{1})
+
\frac{1}{8}
c_{1}(T\widetilde{X},g_{0})|_{\widetilde{X}^{\theta}}\widetilde{c_{1}}(T\widetilde{X}^{\theta};g_{0},g_{1})
\\
&\quad
-
\frac{1}{12}
\widetilde{c_{1}}(T\widetilde{X}^{\theta};g_{0},g_{1})\{c_{1}(T\widetilde{X}^{\theta},g_{0})+c_{1}(T\widetilde{X}^{\theta},g_{1})\},
\end{aligned}
$$
where \cite[Eq.\,(1.3.1.2)]{GilletSoule90} is used to get the second equality.
For a holomorphic line bundle $L$ and Hermitian metrics $h_{0},h_{1}$ on $L$, we have
$$
\widetilde{c_{1}}(L;h_{0},h_{1})=\log(h_{0}/h_{1})
$$
by \cite[Eq.\,(1.2.5.1)]{GilletSoule90}. (Our sign convention is different from the one in Gillet-Soul\'e \cite{GilletSoule90}.
Our $\widetilde{c_{1}}(L;h_{0},h_{1})$ is $-\widetilde{c_{1}}(L;h_{0},h_{1})$ in \cite{GilletSoule90}.)
Hence
\begin{equation}
\label{eqn:anomaly:equivariant:torsion:2}
\begin{aligned}
\widetilde{\rm Td}_{\theta}(T\widetilde{X};g_{0},g_{1})^{(1,1)}
&\equiv
\frac{1}{8}
\left.
\log\left(\frac{\det g_{0}}{\det g_{1}}\right)
\right|_{\widetilde{X}^{\theta}}
c_{1}(\widetilde{X}^{\theta},g_{1})
+
\frac{1}{8}
\left.
c_{1}(\widetilde{X},g_{0})
\log\left(\frac{g_{0}}{g_{1}}
\right|_{\widetilde{X}^{\theta}}
\right)
\\
&\quad
-
\frac{1}{12}
\left.
\log\left(\frac{g_{0}}{g_{1}}
\right|_{\widetilde{X}^{\theta}}
\right)
\{c_{1}(\widetilde{X}^{\theta},g_{0})+c_{1}(\widetilde{X}^{\theta},g_{1})\}
\\
&\qquad
\mod
{\rm Im}\,\partial+{\rm Im}\,\bar{\partial}.
\end{aligned}
\end{equation}
\par{\em (Step 2) }
We set $g_{0}=\gamma_{\epsilon,\delta}$ and $g_{1}=\gamma_{\epsilon(\rho)}$ in Step 1.
Since $g_{0}=\gamma_{\epsilon,\delta}$ is Ricci-flat on a neighborhood of $\widetilde{X}^{\theta}$, we have
$$
c_{1}(\widetilde{X},\gamma_{\epsilon,\delta})|_{\widetilde{X}^{\theta}}=0.
$$
Since the volume form of EH instanton $i\partial\bar{\partial}F_{\epsilon}$ is the standard Euclidean volume form
$$
\frac{(i\partial\bar{\partial}F_{\epsilon})^{2}}{2!}=i^{2}dz_{1}\wedge d\bar{z}_{1}\wedge dz_{2}\wedge d\bar{z}_{2}
$$
and since $\gamma_{\epsilon,\delta}=i\partial\bar{\partial}F_{\epsilon}$ on $\widetilde{V}(\delta)_{\frak p}$,
we get
$$
\left.
\left(
\frac{\det\gamma_{\epsilon,\delta}}{\det\gamma_{\epsilon(\rho)}}\right)
\right|_{\widetilde{X}^{\theta}}
=
\left.
\left(
\frac{\gamma_{\epsilon,\delta}^{2}/2!}{\gamma_{\epsilon(\rho)}^{2}/2!}\right)
\right|_{\widetilde{X}^{\theta}}
=
1.
$$
If $E_{i}\cong{\bf P}^{1}$ is a component of $\widetilde{X}^{\theta}$, 
then $(E_{i},\gamma_{\epsilon,\delta}|_{E_{i}})\cong({\bf P}^{1},\epsilon\,\omega_{\rm FS})$
and $(E_{i},\gamma_{\epsilon(\rho)}|_{E_{i}})\cong({\bf P}^{1},\epsilon(\rho)\omega_{\rm FS})$. Hence
$$
\left.\frac{\gamma_{\epsilon,\delta}}{\gamma_{\epsilon(\rho)}}\right|_{\widetilde{X}^{\theta}}=\frac{\epsilon}{\epsilon(\rho)}.
$$
All together, we get
\begin{equation}
\label{eqn:Bott-Chern:equivariant:Todd}
\begin{aligned}
\,&
\int_{\widetilde{X}^{\theta}}
\widetilde{\rm Td}_{\theta}(T\widetilde{X};\gamma_{\epsilon,\delta},\gamma_{\epsilon(\rho)})^{(1,1)}
=
-\int_{\widetilde{X}^{\theta}}
\frac{1}{12}
\left.
\log\left(\frac{\gamma_{\epsilon,\delta}}{\gamma_{\epsilon(\rho)}}
\right|_{\widetilde{X}^{\theta}}
\right)
\{c_{1}(\widetilde{X}^{\theta},\gamma_{\epsilon,\delta})+c_{1}(\widetilde{X}^{\theta},\gamma_{\epsilon(\rho)})\}
\\
&=
-\frac{\log(\epsilon/\epsilon(\rho))}{6}
\int_{\widetilde{X}^{\theta}}c_{1}(\widetilde{X}^{\theta})
=
-\frac{\log(\epsilon/\epsilon(\rho))}{6}\chi(\widetilde{X}^{\theta})
=
-\frac{k}{3}\,\log\frac{\epsilon}{\epsilon(\rho)},
\end{aligned}
\end{equation}
where we used the fact $\widetilde{X}^{\theta}=E_{1}\amalg\cdots\amalg E_{k}$, $k=\#{\rm Sing}\,X$, $E_{i}\cong{\bf P}^{1}$.
This, together with \eqref{eqn:anomaly:equivariant:torsion}, yields that
$$
\epsilon^{k/3}\tau_{{\bf Z}_{2}}(\widetilde{X},\gamma_{\epsilon,\delta})(\theta){\rm Vol}(\widetilde{X},\gamma_{\epsilon,\delta})
=
\epsilon(\rho)^{k/3}{\rm Vol}(\widetilde{X},\gamma_{\epsilon(\rho)})\,
\tau_{{\bf Z}_{2}}(\widetilde{X},\gamma_{\epsilon(\rho)})(\theta)
$$
is independent of $\epsilon,\delta\in(0,1]$ with $0<\epsilon\delta^{-2}\leq\epsilon(\rho)$.
\par{\em (Step 3) }
Let $\chi$ be another cut-off function to glue Eguchi-Hanson instanton to the initial K\"ahler form $\gamma_{0}$ on $X$
(cf. Sections~5.2.1 and 5.2.3). Then there exists $\epsilon(\chi)\in(0,1)$ such that the function
$\phi'_{\epsilon,\delta}(z):=\|z\|^{2}+\chi_{\delta}(\|z\|)E(z,\epsilon)$ on $V(\infty)\setminus\{0\}$ is a potential of a K\"ahler form
on $T^{*}{\bf P}^{1}=\widetilde{V}(\infty)$ for any $\epsilon,\delta\in(0,1]$ with $0<\epsilon\delta^{-2}\leq\epsilon(\chi)$.
Let $\gamma'_{\epsilon,\delta}$ be the families of K\"ahler forms on $\widetilde{X}$ constructed in the same way as in
\eqref{eqn:family:Kahler:metric:2} using $\kappa'_{\epsilon,\delta}:=i\partial\bar{\partial}\phi'_{\epsilon,\delta}$ instead of $\kappa_{\epsilon,\delta}$.
By Step 2, we get
$$
\epsilon^{k/3}\tau_{{\bf Z}_{2}}(\widetilde{X},\gamma'_{\epsilon,\delta})(\theta){\rm Vol}(\widetilde{X},\gamma'_{\epsilon,\delta})
=
\epsilon(\chi)^{k/3}{\rm Vol}(\widetilde{X},\gamma'_{\epsilon(\chi)})\,
\tau_{{\bf Z}_{2}}(\widetilde{X},\gamma'_{\epsilon(\chi)})(\theta)
$$
for any $\epsilon,\delta\in(0,1]$ with $0<\epsilon\delta^{-2}\leq\epsilon(\chi)$.
To prove the independence of the limit 
$\lim_{\epsilon\to0}\epsilon^{k/3}\tau_{{\bf Z}_{2}}(\widetilde{X},\gamma_{\epsilon,\delta})(\theta){\rm Vol}(\widetilde{X},\gamma_{\epsilon,\delta})$
from the choice of $\rho$, we must prove
\begin{equation}
\label{eqn:independence:cut-off:fcn}
\epsilon(\rho)^{k/3}{\rm Vol}(\widetilde{X},\gamma_{\epsilon(\chi)})\,
\tau_{{\bf Z}_{2}}(\widetilde{X},\gamma_{\epsilon(\chi)})(\theta)
=
\epsilon(\chi)^{k/3}{\rm Vol}(\widetilde{X},\gamma'_{\epsilon(\chi)})\,
\tau_{{\bf Z}_{2}}(\widetilde{X},\gamma'_{\epsilon(\chi)})(\theta).
\end{equation}
We set $g_{0}=\gamma_{\epsilon(\rho)}$ and $g_{1}=\gamma'_{\epsilon(\chi)}$ in \eqref{eqn:anomaly:equivariant:torsion:2}.
By the same computation as in \eqref{eqn:Bott-Chern:equivariant:Todd}, we get
$$
\int_{\widetilde{X}^{\theta}}
\widetilde{\rm Td}_{\theta}(T\widetilde{X};\gamma_{\epsilon(\rho)},\gamma'_{\epsilon(\chi)})^{(1,1)}
=
-\frac{k}{3}\log\frac{\epsilon(\chi)}{\epsilon(\rho)}.
$$
This, together with \eqref{eqn:anomaly:equivariant:torsion}, yields \eqref{eqn:independence:cut-off:fcn}. This completes the proof.
\end{pf}

\subsection
{A comparison of heat kernels}
\par
Recall that 
$K_{\epsilon,\delta}^{q}(t, x, y)$ denote the heat kernel of the Hodge-Kodaira Laplacian 
$\square_{q}^{\epsilon,\delta}$ for the K\"ahler metric 
$\gamma_{\epsilon,\delta}$ on $\widetilde{X}$, 
and $K_{0}^q(t, x, y)$ the heat kernel of the Hodge-Kodaira Laplacian $\square_{q}^{0}$ for the K\"ahler metric $\gamma_{0}$ on $X$.  
For $0<r\leq 4$ let 
$$
\widetilde{V}_r:=\bigcup_{{\frak p}\in{\rm Sing}(X)} \widetilde{V}(r)_{\frak p},
\qquad
\widetilde{X}_r:=\widetilde{X}-\widetilde{V}_r.
$$
Define $\widetilde{V}_{\infty}$ to be $\widetilde{V}_4$ extended by $k$ copies of the infinite cone $(\mathbf C^2 - B(4))/\{\pm 1\}$. 
The metric $\gamma_{\epsilon, \delta}|_{\widetilde{V}_4}$ similarly extends to a K\"ahler metric $\gamma_{\epsilon, \delta}^{\infty}$ on $\widetilde{V}_{\infty}$. 
We denote by $K_{\epsilon, \delta, \infty}^q(t, x, y)$ the corresponding heat kernel on $\widetilde{V}_{\infty}$. 
Similarly we have the corresponding $X_r, V_r, V_{\infty}$ on $X$, with $X_r$ identified with $\widetilde{X}_r$. 
Note that $V_{\infty}$ is just $k$ copies of the infinite cone.

We first established some uniform estimates on the heat kernel $K_{\epsilon,\delta}^{q}(t, x, y)$, improving on Proposition \ref{prop:heat:trace} when the points are in specific regions.

\begin{theorem} 
\label{thm:heat:kernel:estimates:0}	
There are  constants $A, C$ depending only on the Sobolev constant and dimension such that, for all  $\epsilon,\delta\in(0,1]$ with $\epsilon\delta^{-2}\leq\epsilon(\rho)$, and $0\leq q\leq 2$, we have

$$  | K_{\epsilon}^q(t, x, z) | \leq Ae^{C(1+ \epsilon \delta^{-4})} \delta^{-4} e^{-\frac{\delta^2}{32t}},  \ \ \ \forall x \in \widetilde{X}_{3\delta}, \ \  z\in  \widetilde{V}_{2\delta}, \ \ t>0. $$ 
Similarly we have,  $\forall x \in \widetilde{X}_{3\delta}, \  z\in  \widetilde{V}_{2\delta}, \  t>0$,
$$ |d K_{\epsilon}^q(t, x, z)| \leq Ae^{C(1+ \epsilon \delta^{-4})} \delta^{-5}  e^{-\frac{\delta^2}{32t}}, \ \ | d^*_{\epsilon, \delta} K_{\epsilon}^q(t, x, z) | \leq Ae^{C(1+ \epsilon \delta^{-4})} \delta^{-5}  e^{-\frac{\delta^2}{32t}}, $$
Here $d, \ d^*_{\epsilon, \delta}$ could act either on $x$ or $z$ variable. 
Finally, for $0<r<2\delta$,  $x \in \widetilde{X}_{3\delta}, \ \  z\in  \widetilde{V}_{2\delta, r} = \widetilde{V}_{2\delta}-\widetilde{V}_{r}$, 
and $i\in \mathbf N$, 
$$
| \nabla^i K_{\epsilon, \delta}^q(t-s, x, z) | \leq  C(i, \delta, r)  e^{-\frac{\delta^2}{32t}}, 
$$
for a constant $C(i, \delta, r)$ depending on $i, \delta, r$. Here $\nabla^i$ denotes the $i$-th covariant derivative with respect to the metric  $\gamma_{\epsilon,\delta}$, acting on either variable.
\end{theorem}

\begin{pf}
Throughout the proof we fix $ x \in \widetilde{X}_{3\delta}, \ \  z\in  \widetilde{V}_{2\delta}, \ \ t>0.$	
Since the Ricci curvature of $\gamma_{\epsilon, \delta}$ is  bounded by Lemma~\ref{lemma:Ricci},  
the Sobolev estimate together with the Moser iteration technique combined with the finite propagation speed argument as in
Cheeger-Gromov-Taylor \cite{CGT} gives the uniform estimate
$$  | K_{\epsilon}^q(t, x, z) | \leq Ae^{C(1+ \epsilon \delta^{-4})} \delta^{-4} e^{-\frac{\delta^2}{32t}}. $$
Indeed the finite propagation speed technique gives us the $L^2$ estimate
$$
 \| K_{\epsilon}^q(t, \cdot, \cdot) \|_{L^2(B_{\delta/4}(x)\times B_{\delta/4}(z) )} \leq c e^{-\frac{\delta^2}{16t}}
$$
for some uniform constant $c$. Now Moser iteration as in \cite{CGT}[pp.16-26], together with semi-group domination \cite{HessSchraderUhlenbrock80} yields
the desired estimate.

For the estimate on  $ d K_{\epsilon}^q(t, x, z), \  d^*_{\epsilon, \delta} K_{\epsilon}^q(t, x, z)$, let $\eta(r)$ be a smooth cut-off function which is identically $1$ for $|r|\leq \delta/8$ and identically $0$ for $|r|\geq \delta/4$, and $|\eta'| \leq \frac{16}{\delta}$. We will continue to denote by 
$\eta$ its composition with a distance function (either $d(x, \cdot)$ or $d(z, \cdot)$). Note then 
\begin{eqnarray*}
 \|(d+d^*_{\epsilon, \delta} )_z[\eta K_{\epsilon}^q(t, \cdot, \cdot)] \|^2_{L^2(B_{\delta/4}(x)\times B_{\delta/4}(z) )}=
 \|(d)_z[\eta K_{\epsilon}^q(t, \cdot, \cdot)] \|^2_{L^2(B_{\delta/4}(x)\times B_{\delta/4}(z) )} \\
  + 
 \|(d^*_{\epsilon, \delta} )_z[\eta K_{\epsilon}^q(t, \cdot, \cdot)] \|^2_{L^2(B_{\delta/4}(x)\times B_{\delta/4}(z) )},
\end{eqnarray*}
from which we deduce 
\begin{eqnarray*}
\|(d)_z[ K_{\epsilon}^q(t, \cdot, \cdot)] \|_{L^2(B_{\delta/8}(x)\times B_{\delta/8}(z) )} & \leq &
	\|(d+d^*_{\epsilon, \delta} )_z[\eta K_{\epsilon}^q(t, \cdot, \cdot)] \|_{L^2(B_{\delta/4}(x)\times B_{\delta/4}(z) )}
	 \\
& &	+ 
\frac{16}{\delta}	\| K_{\epsilon}^q(t, \cdot, \cdot) \|_{L^2(B_{\delta/4}(x)\times B_{\delta/4}(z) )} \\
& \leq &
\|(d+d^*_{\epsilon, \delta} )_z[ K_{\epsilon}^q(t, \cdot, \cdot)] \|_{L^2(B_{\delta/4}(x)\times B_{\delta/4}(z) )}
\\
& &	+ 
\frac{32}{\delta}	\| K_{\epsilon}^q(t, \cdot, \cdot) \|_{L^2(B_{\delta/4}(x)\times B_{\delta/4}(z) )} 
\end{eqnarray*}

Now the same finite propagation speed technique gives
$$
 \|(d+d^*_{\epsilon, \delta} )_z[ K_{\epsilon}^q(t, \cdot, \cdot)] \|_{L^2(B_{\delta/4}(x)\times B_{\delta/4}(z) )} \leq c' e^{-\frac{\delta^2}{16t}}, 
$$
which in turn gives
$$
\|(d)_z[ K_{\epsilon}^q(t, \cdot, \cdot)] \|_{L^2(B_{\delta/8}(x)\times B_{\delta/8}(z) )} \leq (c'+ c\frac{32}{\delta}) e^{-\frac{\delta^2}{32t}}, 
$$
The same method as above then yields
$$
|(d)_z[ K_{\epsilon}^q(t, x, z)] | \leq Ae^{C(1+ \epsilon \delta^{-4})}  \delta^{-5}  e^{-\frac{\delta^2}{32t}}.
$$
The others can be proven in exactly the same way.

Finally, for $0<r<2\delta$,  we note that the  curvature tensor and its derivatives of  $\gamma_{\epsilon,\delta}$ are bounded in $\widetilde{V}_{2\delta, r} = \widetilde{V}_{2\delta}-\widetilde{V}_{r}$ by a constant depending on $\delta, r$. Moreover the injectivity radius  of $\gamma_{\epsilon,\delta}$ in  $\widetilde{V}_{2\delta, r}$  is bounded away from zero by a constant depending on $\delta, r$. Hence, by the elliptic estimate combined with the argument as before,  we have,  for $x \in \widetilde{X}_{3\delta}, \ \  z\in  \widetilde{V}_{2\delta, r} = \widetilde{V}_{2\delta}-\widetilde{V}_{r}$,
	 and $i\in \mathbf N$, 
	$$
	| \nabla^i K_{\epsilon, \delta}^q(t-s, x, z) | \leq  C(i, \delta, r)  e^{-\frac{\delta^2}{32t}}, 
	$$
	for a constant $C(i, \delta, r)$ depending on $i, \delta, r$.
\end{pf}


\begin{theorem} \label{thm:heat:kernel:estimates:1}
There are  constants $A, C$ depending only on the Sobolev constant and dimension such that, for all  $\epsilon,\delta\in(0,1]$ with $\epsilon\delta^{-2}\leq\epsilon(\rho)$, and $0\leq q\leq 2$, we have
	$$  
	| K_{\epsilon, \delta}^q(t, x, y)- K_{0}^q(t, x, y) | \leq
Ae^{C(1+ \epsilon \delta^{-4})}	\delta^{-9} e^{-\frac{\delta^2}{16t}} {\rm vol}(\partial \widetilde{X}_{2\delta}),  
	\ \ \ \forall x, y \in \widetilde{X}_{3\delta},  \ \ t>0.  
	$$ 
Furthermore, $\forall x, y \in \widetilde{X}_{3\delta},  \ \ t>0$, we have the pointwise (although not necessarily uniform) convergence as $\epsilon \rightarrow 0$,
$$
 K_{\epsilon, \delta}^q(t, x, y)- K_{0}^q(t, x, y) \longrightarrow 0.
$$
\end{theorem}

\begin{pf}

For $0<r\leq 4$,  we apply the Duhamel principle \cite[(3.9)]{C} to $K_{\epsilon, \delta}^q(t, x, y)- K_{0}^q(t, x, y)$ on $\widetilde{X}_r$ to obtain
\begin{eqnarray*}
	K_{\epsilon, \delta}^q(t, x, y)- K_{0}^q(t, x, y) & = & -\int_0^t\int_{\widetilde{X}_r} \left[ (\partial_t  + \square_0^q) K_{\epsilon, \delta}^q(t-s, x, z) \right] \wedge * K_{0}^q(s, z, y) \\
	& &  + \int_0^t\int_{\partial \widetilde{X}_r}  * d  K_{\epsilon, \delta}^q(t-s, x, z)  \wedge K_{0}^q(s, z, y)  \\
	& & + (-1)^{4q+1} \int_0^t\int_{\partial \widetilde{X}_r}  K_{\epsilon, \delta}^q(t-s, x, z)  \wedge * d  K_{0}^q(s, z, y) \\
	& & + (-1)^{4q+1} \int_0^t\int_{\partial \widetilde{X}_r}  *  K_{\epsilon, \delta}^q(t-s, x, z)  \wedge d^*  K_{0}^q(s, z, y) \\
	& &  + \int_0^t\int_{\partial \widetilde{X}_r}  d^*  K_{\epsilon, \delta}^q(t-s, x, z)  \wedge *  K_{0}^q(s, z, y) .
\end{eqnarray*}

Now fix $x, y \in \widetilde{X}_{3\delta}$. First we let $r=2\delta$. Then the first term on the right hand side goes away and we are left with only boundary terms. By Theorem \ref{thm:heat:kernel:estimates:0}, and noticing that similar estimates hold for the orbifold heat kernel
\begin{equation}
\label{eqn:heat:kernel:estimate:2}
| K_{0}^q(t, x, z) | \leq C\delta^{-4} e^{-\frac{\delta^2}{32t}},  \ \ \ \forall x \in \widetilde{X}_{3\delta}, \ \  z\in  \widetilde{V}_{2\delta}, \ \ t>0,
\end{equation}
as well as its derivatives, we deduce then that
$$
	| K_{\epsilon, \delta}^q(t, x, y)- K_{0}^q(t, x, y) | \leq
Ae^{C(1+ \epsilon \delta^{-4})}	\delta^{-9} e^{-\frac{\delta^2}{16t}} {\rm vol}(\partial \widetilde{X}_{2\delta}). 
$$

To prove the pointwise convergence, we let $r<2\delta$, and denote $\widetilde{V}_{2\delta, r}=\widetilde{V}_{2\delta}-\widetilde{V}_r$. Then the Duhamel principle becomes
\begin{eqnarray*}
	K_{\epsilon, \delta}^q(t, x, y)- K_{0}^q(t, x, y) & = & -\int_0^t\int_{\widetilde{V}_{2\delta}, r} \left[ (\partial_t  + \square_0^q) K_{\epsilon, \delta}^q(t-s, x, z) \right] \wedge * K_{0}^q(s, z, y) \\
	& &  + \int_0^t\int_{\partial \widetilde{X}_r}  * d  K_{\epsilon, \delta}^q(t-s, x, z)  \wedge K_{0}^q(s, z, y)  \\
	& & + (-1)^{4q+1} \int_0^t\int_{\partial \widetilde{X}_r}  K_{\epsilon, \delta}^q(t-s, x, z)  \wedge * d  K_{0}^q(s, z, y) \\
	& & + (-1)^{4q+1} \int_0^t\int_{\partial \widetilde{X}_r}  *  K_{\epsilon, \delta}^q(t-s, x, z)  \wedge d^*  K_{0}^q(s, z, y) \\
	& &  + \int_0^t\int_{\partial \widetilde{X}_r}  d^*  K_{\epsilon, \delta}^q(t-s, x, z)  \wedge *  K_{0}^q(s, z, y) .
\end{eqnarray*}

Since $ (\partial_t  + \square_0^q) K_{\epsilon}^q(t-s, x, z)= (\square_0^q-\square_{\epsilon}^q) K_{\epsilon}^q(t-s, x, z)$, 
and $\gamma_{\epsilon, \delta}= i \partial \overline{\partial} \phi_{\epsilon, \delta}$, 
$\phi_{\epsilon, \delta}(z)=\|z\|^2 + \rho_{\delta}(z)E(z, \epsilon)$ on $\widetilde{V}_{2\delta, r}$, 
by \eqref{eqn:5:estimate:error}, \eqref{eqn:5:estimate2:error}, \eqref{eqn:5:estimate3:error}, we have, for $x \in \widetilde{X}_{3\delta}, \ \  z\in  \widetilde{V}_{2\delta, r}$, 
$$
| (\partial_t  + \square_0^q) K_{\epsilon, \delta}^q(t-s, x, z) | \leq \epsilon C(\delta, r)  e^{-\frac{\delta^2}{32t}}, 
$$
for a constant $C(\delta, r)$ depending on $\delta, r$ but not on $\epsilon$. 

Combining with the uniform estimates in Theorem \ref{thm:heat:kernel:estimates:1}, 
we obtain, for $x, y \in \widetilde{X}_{3\delta}$,
$$ | K_{\epsilon, \delta}^q(t, x, y)- K_{0}^q(t, x, y) | \leq \epsilon t  C'(\delta, r)  e^{-\frac{\delta^2}{16t}} + C''(\delta)t  e^{-\frac{\delta^2}{16t}} {\rm vol}(\partial \widetilde{X}_r). $$

Now for any $\eta>0$, we take $r$ sufficiently small so that $ C''(\delta)t  e^{-\frac{\delta^2}{16t}} {\rm vol}(\partial \widetilde{X}_r)< \frac{\eta}{2}$. Then we take $\epsilon$ sufficiently small such that $\epsilon t  C'(\delta, r)  e^{-\frac{\delta^2}{16t}}< \frac{\eta}{2}$. Hence 
$$ | K_{\epsilon, \delta}^q(t, x, y)- K_{0}^q(t, x, y) | < \eta. $$
This proves the pointwise convergence.
\end{pf}

\begin{remark}
Since we have the Ricci curvature lower bound, the pointwise convergence of the heat kernels  should also be a consequence of 
some general spectral convergence results due to Cheeger-Colding \cite{CC} for the case $q=0$, Honda \cite{Honda17} for the case $q=1$, 
and Bei \cite{Bei19} for $q=n=2$. See also \cite{Di}.
\end{remark}

Our next task is to compare the heat kernel $K_{0}^q(t, x, y)$ for $(X, \gamma_{0})$ with the heat kernel $K_{0, \infty}^q(t, x, y)$ of $V_{\infty}$ when $x, y \in V_{3\delta}$.

\begin{theorem} \label{thm:heat:kernel:estimates:2}
	There is a constant $C$ depending only on the Sobolev constant and dimension such that, for $\delta\leq 1$, 
	$$  
	| K_{0}^q(t, x, y)- K_{0, \infty}^q(t, x, y) | \leq
	C e^{-\frac{1}{16t}} {\rm vol}(\partial V_4),  
	\ \ \ \forall x, y \in V_{3\delta},  \ \ t>0.  
	$$ 
\end{theorem}

\begin{pf}
The Duhamel principle \cite[(3.9)]{C} applied to $K_{0}^q(t, x, y)- K_{0, \infty}^q(t, x, y)$ on $V_4$ gives us
\begin{eqnarray*}
	K_{0}^q(t, x, y)-  K_{0, \infty}^q(t, x, y) & = &
 \int_0^t\int_{\partial V_4}  * d  K_{0}^q(s, x, z)  \wedge K_{0, \infty}^q(t-s, z, y)  \\
	& & + (-1)^{4q+1} \int_0^t\int_{\partial V_4}  K_{0}^q(s, x, z)  \wedge * d K_{0, \infty}^q(t-s, z, y)  \\
	& & + (-1)^{4q+1} \int_0^t\int_{\partial V_4}  *  K_{0}^q(s, x, z)  \wedge d^* K_{0, \infty}^q(t-s, z, y)  \\
	& &  + \int_0^t\int_{\partial V_4}  d^*  K_{0}^q(s, x, z)  \wedge *  K_{0, \infty}^q(t-s, z, y)  .
\end{eqnarray*}

Thus we obtain, for $x, y\in V_{3\delta}, \delta \leq 1$,  using the estimate \ref{eqn:heat:kernel:estimate:2}, except with the $\delta$ there replaced by a fixed constant, say $1/4$, as well as a similar estimate for $K_{0, \infty}^q(t, x, y) $, 
$$
 |	K_{0}^q(t, x, y)-  K_{0, \infty}^q(t, x, y) | \leq C e^{-\frac{1}{16t}}  {\rm vol}(\partial V_4) .
$$

\end{pf}

Our final task here is to compare the heat kernel $K_{\epsilon,\delta}^{q}(t, x, y)$ with
$K_{\epsilon, \delta, \infty}^q(t, x, y)$,  the heat kernel on $\widetilde{V}_{\infty}$, when $x, y \in \widetilde{V}_{3\delta}$.

\begin{theorem} \label{thm:heat:kernel:estimates:3}
There are  constants $A, C$ depending only on the Sobolev constant and dimension such that, for all  $\epsilon,\delta\in(0,1]$ with $\epsilon\delta^{-2}\leq\epsilon(\rho)$, and $0\leq q\leq 2$, we have
$$  
	| K_{\epsilon,\delta}^{q}(t, x, y)- K_{\epsilon,\delta, \infty}^{q}(t, x, y) | \leq
Ae^{C(1+ \epsilon \delta^{-4})} e^{-\frac{1}{16t}} {\rm vol}(\partial V_4),  
	\ \ \ \forall x, y \in \widetilde{V}_{3\delta},  \ \ t>0.  
$$ 
\end{theorem}

\begin{pf}
The proof follows the same line as above. We apply the Duhamel principle to $K_{\epsilon,\delta}^{q}(t, x, y)- K_{\epsilon,\delta, \infty}^{q}(t, x, y)$ on $\widetilde{V}_{4}$ and use the heat kernel estimate in Theorem \ref{thm:heat:kernel:estimates:0} as well as the analogous estimate for $K_{\epsilon,\delta, \infty}^{q}(t, x, y)$ to obtain the desired estimate.
\end{pf}

%

\subsection
{Partial analytic torsion}
\par
Recall that in Section 4.1, for a compact K\"ahler orbifold $(Z, \gamma)$ of diemsnion $n$,
$$
\zeta_{q}(s)=\sum_{\lambda\in\sigma(\square_{q})\setminus\{0\}}\lambda^{-s}\,\dim E(\lambda;\square_{q}) 
=
\frac{1}{\Gamma(s)} \int_0^{\infty} t^{s-1} {\rm Tr}(e^{-t\square_q }P^{\perp}_q) \, dt
$$
with $P^{\perp}_q$ the orthogonal projection onto the orthogonal compliment of ${\rm ker}\, \square_q$, and (the logarithm of) the analytic torsion
$$
\ln \tau(Z,\gamma)=-\sum_{q=0}^{n}(-1)^{q}q\,\zeta'_{q}(0)= -\zeta'_{T}(0),
$$
where
$$
\zeta_{T}(s)=\frac{1}{\Gamma(s)} \int_0^{\infty} t^{s-1} {\rm Tr}_s(Ne^{-t\square}P^{\perp}) \, dt.
$$
Here $\square$ denotes the Hodge-Kodaira Laplacian on $A^{0, *}(Z)$, $P^{\perp}$ the orthogonal projection onto the orthogonal compliment of ${\rm ker}\, \square$,  
${\rm Tr}_s$ the supertrace on $A^{0, *}(Z)$, i.e., the alternating sum of the traces on each degree, and $N$ the so called number operator 
which simply multiply a differential form by its degree.


By Lidskii theorem
\begin{eqnarray*}
{\rm Tr}_s(Ne^{-t\square}P^{\perp}) & = & \int_Z {\rm tr}_s(N K(t, x, x)P^{\perp}(x, x)) \, dx \\
& = & \sum_{q=0}^{n}(-1)^{q}q\,\int_Z {\rm tr}(K_q(t, x, x)P^{\perp}(x, x)) \, dx
\end{eqnarray*}
where $K(t, x, y)$, $K_q(t, x, y)$ denotes the heat kernel of $\square$, $\square_q$, respectively, $P^{\perp}(x, x)$ the Schwartz kernel of $P^{\perp}$, 
and ${\rm tr}_s$ (abusing notation) also the pointwise supertrace.

At this point it is convenient to introduce what is called ``partial analytic torsion" in \cite{D}. For a domain $D\subset Z$, we define
$$
\zeta_{T}^{D, Z}(s)=\frac{1}{\Gamma(s)} \int_0^{\infty} t^{s-1} \int_D {\rm tr}_s(N K(t, x, x)P^{\perp}(x, x)) \, dx \, dt
$$
and
\begin{equation} \label{pat}
\ln \tau(D, Z,\gamma)= -\left(\zeta_T^{D, Z}\right)'(0).
\end{equation}

Clearly
\begin{equation} \label{aopat}
\ln \tau(Z,\gamma)=\ln \tau(D, Z,\gamma)+ \ln \tau(Z-D, Z,\gamma).
\end{equation}

Similarly we can define the equivariant version $\tau_{{\bf Z}_{2}}(D, Z, \gamma)(\theta)$ for $\theta$-invariant domain $D\subset Z$. 
That is, we define 
$$
\zeta_{T, \theta}^{D, Z}(s)=\frac{1}{\Gamma(s)} \int_0^{\infty} t^{s-1} \int_D {\rm tr}_s(N K(t, x,\theta x)P^{\perp}(x, \theta x)) \, dx \, dt
$$
and
\begin{equation}   \label{peat}
\ln \tau_{{\bf Z}_{2}}(D, Z, \gamma)(\theta)= -\left(\zeta_{T, \theta}^{D, Z}\right)'(0).
\end{equation}
Then the discussion applies to the equivariant version as well.

\subsection
{Limit of partial analytic torsion I}
\par

\begin{theorem} \label{lopat1} 
For $0<\delta \leq 1$, we have
$$
\lim_{\epsilon \rightarrow 0} \ln \tau(\widetilde{X}_{3\delta}, \widetilde{X}, \gamma_{\epsilon, \delta})= \ln \tau(X_{3\delta}, X, \gamma_{0}),
$$
and
$$
\lim_{\epsilon \rightarrow 0} \ln \tau_{{\bf Z}_{2}}(\widetilde{X}_{3\delta}, \widetilde{X}, \gamma_{\epsilon, \delta})(\theta)
=
 \ln \tau_{{\bf Z}_{2}}(X_{3\delta}, X, \gamma_{0})(\iota).
$$
\end{theorem}

\begin{pf}
\par{\em (Step 1) }
Let
$$
 {\rm tr}_s(N K_{\epsilon, \delta}(t, x, x)) \sim \sum_{i=0}^{\infty} a_i^{\epsilon, \delta}(x)\, t^{i-2}
$$
be the pointwise small time asymptotic expansion, and write
\begin{eqnarray*}
\zeta_{T}^{\widetilde{X}_{3\delta}, \widetilde{X}}(s) & = & \frac{1}{\Gamma(s)} 
\left[ \int_1^{\infty} t^{s-1} \int_{\widetilde{X}_{3\delta}} {\rm tr}_s(N K_{\epsilon, \delta}(t, x, x)P^{\perp}_{\epsilon, \delta}(x, x)) \, dx \, dt \right. \\
& & +   \int_0^{1} t^{s-1} \int_{\widetilde{X}_{3\delta}} [ {\rm tr}_s(N K_{\epsilon, \delta}(t, x, x))- \sum_{i=0}^{2} a_i^{\epsilon, \delta}(x)\, t^{i-2} ] \, dx \, dt \\
& & + \left. \sum_{i=0}^{1} \int_{\widetilde{X}_{3\delta}} \frac{a_i^{\epsilon, \delta}(x)}{s+ i-2} \, dx  
+ \frac{1}{s} \int_{\widetilde{X}_{3\delta}}[ a_2^{\epsilon, \delta}(x) -  {\rm tr}_s(N P_{\epsilon, \delta}(x, x))] \, dx \right] ,
\end{eqnarray*}
where $P_{\epsilon, \delta}(x, x)$ is the Schwartz kernel of  $P_{\epsilon, \delta}$, the orthogonal projection onto ${\rm ker}\, \square_{\epsilon, \delta}$. 
We obtain
\begin{eqnarray*}
\ln \tau(\widetilde{X}_{3\delta}, \widetilde{X}, \gamma_{\epsilon, \delta}) & = & - \int_1^{\infty} t^{-1} 
\int_{\widetilde{X}_{3\delta}} {\rm tr}_s(N K_{\epsilon, \delta}(t, x, x)P^{\perp}_{\epsilon, \delta}(x, x)) \, dx \, dt \\
& & -  \int_0^{1} t^{-1} \int_{\widetilde{X}_{3\delta}}\left[ {\rm tr}_s(N K_{\epsilon, \delta}(t, x, x))- \sum_{i=0}^{2} a_i^{\epsilon, \delta}(x)\, t^{i-2}\right] \, dx \, dt \\
& & - \sum_{i=0}^{1} \int_{\widetilde{X}_{3\delta}} 
\frac{a_i^{\epsilon, \delta}(x)}{ i-2} \, dx + \Gamma'(1) \int_{\widetilde{X}_{3\delta}}\left[ a_2^{\epsilon, \delta}(x) -  {\rm tr}_s(N P_{\epsilon, \delta}(x, x))\right] \, dx
\end{eqnarray*}
and similarly for $\ln \tau(X_{3\delta}, X, \gamma_{0})$.
Since the asymptotic expansion depends only on the local data, we have $a_i^{\epsilon, \delta}(x)=a_i^{0}(x)$ on $\widetilde{X}_{3\delta}$. 
Hence
\begin{equation}
\label{eqn:difference:partial:torsion}
\begin{aligned}
\,&
\ln \tau(\widetilde{X}_{3\delta}, \widetilde{X}, \gamma_{\epsilon, \delta}) - \ln \tau(X_{3\delta}, X, \gamma_{0})
\\
&=
- \int_1^{\infty} t^{-1} 
\int_{\widetilde{X}_{3\delta}} {\rm tr}_s \left[ N K_{\epsilon, \delta}(t, x, x)P^{\perp}_{\epsilon, \delta}(x, x) - N K_{0}(t, x, x)P^{\perp}_{0}(x, x) \right] \, dx \, dt
\\
&\quad
-  \int_0^{1} t^{-1} \int_{\widetilde{X}_{3\delta}}\left[ {\rm tr}_s(N K_{\epsilon, \delta}(t, x, x) -N K_{0}(t, x, x)) \right] \, dx \, dt 
\\
&\quad
-  \Gamma'(1) \int_{\widetilde{X}_{3\delta}} {\rm tr}_s(N P_{\epsilon, \delta}(x, x) - N P_{0}(x, x))\, dx.
\end{aligned}
\end{equation}
We estimate each term of the right hand side.
\par{\em (Step 2) }
Let $\Lambda>0$ be the same constant as in Corollary~\ref{cor:uniform:exp:decay}.
By Corollary~\ref{cor:uniform:exp:decay}, 
$$
\begin{aligned}
\int_{\widetilde{X}_{3\delta}}\left|{\rm tr}_s(N K_{\epsilon, \delta}(t, x, x)P^{\perp}_{\epsilon, \delta}(x, x))\right|\,dx
&\leq
\sum_{q\geq0}q\,({\rm Tr}\,e^{-t\square_{\epsilon,\delta}}-h^{0,q}(\widetilde{X}))
\\
&\leq
\Psi(\epsilon\delta^{-4}+1)\,\exp[-\frac{1}{2}t\Lambda(\epsilon\delta^{-4}+1)]
\end{aligned}
$$
for all $\epsilon,\delta \in (0,1]$ with $\epsilon\delta^{-2}\leq\epsilon(\rho)$ and $t\geq1$,
where $\Psi(R)$ and $\Lambda(R)$ were defined in Corollary~\ref{cor:uniform:exp:decay}.
Hence for any $\nu >0$, there is $T'=T'(\nu)>0$ depending only on $\nu$ such that for all $\epsilon,\delta \in (0,1]$ with 
$\epsilon<\min\{\epsilon(\rho)\delta^{2},\delta^{4}\}$, 
and $T>T'$,
\begin{equation}
\label{eqn:estimate:[T,infty]}
\int_T^{\infty} t^{-1} \int_{\widetilde{X}_{3\delta}} \left|{\rm tr}_s(N K_{\epsilon, \delta}(t, x, x)P^{\perp}_{\epsilon, \delta}(x, x)) \right|\, dx \, dt 
\leq
\Psi(2)\int_{T}^{\infty}e^{-\Lambda(2) t/2}\frac{dt}{t}
< \nu
\end{equation}
and similarly for the same term involving $K_0$.
By Theorem \ref{thm:heat:kernel:estimates:1} and Lebesgue dominated convergence theorem, there exists $\epsilon_0>0$ such that 
\begin{equation}
\label{eqn:difference:heat:kernel}
\left| \int_1^{T} t^{-1} \int_{\widetilde{X}_{3\delta}} {\rm tr}_s[N( K_{\epsilon, \delta}(t, x, x)) - K_{0}(t, x, x))] \, dx \, dt  \right| 
< \nu,
\end{equation}
whenever $\epsilon < \epsilon_0$. Similarly,
\begin{equation}
\label{eqn:difference:heat:kernel:[0,1]}
\left| \int_0^{1} t^{-1} \int_{\widetilde{X}_{3\delta}} {\rm tr}_s[N( K_{\epsilon, \delta}(t, x, x)) - K_{0}(t, x, x))] \, dx \, dt  \right| 
< \nu,
\end{equation}
whenever $\epsilon < \epsilon_0$.

On the other hand
\begin{equation}
\label{eqn:heat:ker:harmonic:proj}
{\rm tr}_s(N K_{\epsilon, \delta}(t, x, x)P^{\perp}_{\epsilon, \delta}(x, x)) ={\rm tr}_s(N K_{\epsilon, \delta}(t, x, x)) - {\rm tr}_s(N P_{\epsilon, \delta}(x, x))
\end{equation}
and similarly for $K_0$. 
Recall $\ker\square_{0}=\ker\square_{\epsilon,\delta}={\bf C}\cdot1\oplus{\bf C}\cdot\bar{\eta}$. 
For $x\in\widetilde{X}_{3\delta}$, we get
\begin{equation}
\label{eqn:harmonic:projection}
{\rm tr}_{s}[N(P_{\epsilon,\delta}(x,x)-P_{0}(x,x))]
=
\frac{2}{\|\eta\|_{L^{2}}^{2}}
\left(
\frac{\eta\wedge\bar{\eta}}{\gamma_{\epsilon,\delta}^{2}/2!}(x)-\frac{\eta\wedge\bar{\eta}}{\gamma_{0}^{2}/2!}(x)
\right)
=0,
\end{equation}
because $\gamma_{\epsilon,\delta}=\gamma_{0}$ on $\widetilde{X}_{3\delta}$.
It follows from \eqref{eqn:difference:heat:kernel}, \eqref{eqn:heat:ker:harmonic:proj}, \eqref{eqn:harmonic:projection} that
\begin{equation}
\label{eqn:integral:[1,T]}
 \left| \int_1^{T} t^{-1} \int_{\widetilde{X}_{3\delta}} {\rm tr}_s(N[ K_{\epsilon, \delta}(t, x, x)P^{\perp}_{\epsilon, \delta}(x, x)-K_0(t, x, x)P^{\perp}_0(x, x)]) \, dx \, dt \right| 
< \nu  .
\end{equation}
%
Substituting \eqref{eqn:estimate:[T,infty]}, \eqref{eqn:difference:heat:kernel}, \eqref{eqn:difference:heat:kernel:[0,1]}, \eqref{eqn:harmonic:projection}, \eqref{eqn:integral:[1,T]} 
into \eqref{eqn:difference:partial:torsion}, we get
$$
\left| \ln \tau(\widetilde{X}_{3\delta}, \widetilde{X}, \gamma_{\epsilon, \delta}) -\ln \tau(X_{3\delta}, X, \gamma_{0}) \right| < 3\nu,
$$
whenever $\epsilon < \epsilon_0$.
Since $\nu>0$ can be chosen arbitrarily small, this finishes the proof of the first formula. 

To prove the  result about the equivariant torsion,  we follow the  same line of argument, except with a simplification, since $\theta$ has no fixed points in  $\widetilde{X}_{3\delta}$. Indeed,  
\begin{eqnarray*}
	\ln \tau_{{\bf Z}_{2}}(\widetilde{X}_{3\delta}, \widetilde{X}, \gamma_{\epsilon, \delta})(\theta) & = & - \int_1^{\infty} t^{-1} 
	\int_{\widetilde{X}_{3\delta}} {\rm tr}_s(N K_{\epsilon, \delta}(t, x, \theta x)P^{\perp}_{\epsilon, \delta}(x, \theta x)) \, dx \, dt \\
	& & -  \int_0^{1} t^{-1} \int_{\widetilde{X}_{3\delta}} {\rm tr}_s(N K_{\epsilon, \delta}(t, x, \theta x))\, dx \, dt \\
	& & - \Gamma'(1) \int_{\widetilde{X}_{3\delta}} {\rm tr}_s(N P_{\epsilon, \delta}(x, \theta x)\, dx
\end{eqnarray*}
and 
\begin{equation}
\label{eqn:difference:equivariant:partial:torsion}
\begin{aligned}
\,&
	\ln \tau_{{\bf Z}_{2}}(\widetilde{X}_{3\delta}, \widetilde{X}, \gamma_{\epsilon, \delta})(\theta) -\ln \tau_{{\bf Z}_{2}}(X_{3\delta}, X, \gamma_{0})(\iota)
\\
&=
- \int_1^{\infty} t^{-1} 
\int_{\widetilde{X}_{3\delta}} {\rm tr}_s \left[ N K_{\epsilon, \delta}(t, x, \theta x)P^{\perp}_{\epsilon, \delta}(x, \theta x) - N K_{0}(t, x, \theta x)P^{\perp}_{0}(x,\theta  x) \right] \, dx \, dt
\\
&\quad
-  \int_0^{1} t^{-1} \int_{\widetilde{X}_{3\delta}}\left[ {\rm tr}_s(N K_{\epsilon, \delta}(t, x,\theta  x) -N K_{0}(t, x,\theta  x)) \right] \, dx \, dt 
\\
&\quad
-  \Gamma'(1) \int_{\widetilde{X}_{3\delta}} {\rm tr}_s(N P_{\epsilon, \delta}(x, \theta  x) - N P_{0}(x,\theta  x))\, dx.
\end{aligned}
\end{equation}
Now we proceed as before.
\end{pf}

\subsection
{Limit of partial analytic torsion II}
\par
To relate $\ln \tau(X_{3\delta}, X, \gamma_{0})$ to $\ln \tau(X, \gamma_{0})$, by (\ref{aopat}), it suffices to show

\begin{theorem} \label{lopat3} We have
$$
\lim_{\delta \rightarrow 0} \ln \tau(V_{3\delta}, X, \gamma_0)= 0, \ \ \ \lim_{\delta \rightarrow 0} \ln \tau_{{\bf Z}_{2}}(V_{3\delta}, X, \gamma_0)(\iota)=0.
$$
\end{theorem}

\begin{remark}
This is closely related to \cite{DY} where analytic torsions on orbifolds defined from conical singularity pointview are shown to be the same as the ones defined from orbifold singularity pointview. 
\end{remark}

\begin{pf} 
Again the proof for both formulas work the same so we only present the first one. Moreover, the argument works for any orbifold singularity but 
we will work with the cyclic quotient singularity of type $\frac{1}{4}(1,1)$ in our situation. 
First of all, by the same kind of argument as above, using Theorem \ref{thm:heat:kernel:estimates:2} and ${\rm vol}(V_{3\delta}) \rightarrow 0$ as $\delta \rightarrow 0$, one has
$$
\lim_{\delta \rightarrow 0} \ln \tau(V_{3\delta}, X, \gamma_0)= \lim_{\delta \rightarrow 0} \ln \tau(V_{3\delta}, V_{\infty}, \kappa_0).
$$
Now the right hand side can be explicitly computed since the heat kernel is explicitly known. Indeed, as $V_{\infty}$ is just $k$ copies of 
$\mathbf C^2/\mathbf Z_2$, the (orbifold) heat kernel of $( V_{\infty}, \kappa_0)$ on the $(0,q)$ forms 
is $k\binom{n}{q}$ ($n=2$ in our case) copies of
$$
K_0(t, x,x')= \frac{1}{(4\pi t)^{n/2}}\left( e^{-\frac{|x-x'|^2}{4t}} + e^{-\frac{|x+x'|^2}{4t}} \right).
$$
In terms of the polar coordinates $x=(r, y), y\in S^{2n-1}$,
$$
K_0(t, x, x)= \frac{1}{(4\pi t)^{n/2}}(1+ e^{-r^2/t}).
$$
Thus
$$
\int_{V_{3\delta}} K_0(t, x, x) dx = c_n \delta^n t^{-n/2} + d_n \int_0^{\frac{3\delta}{t^{1/2}}} \xi^{n-1}e^{-\xi^2} d\xi,
$$
where $c_n=\frac{3^n\omega_n}{n(4\pi)^{n/2}}$, $d_n=\frac{\omega_n}{(4\pi)^{n/2}}$, $\omega_n={\rm vol}(S^{n-1})$.

The second term has different asymptotic behaviors for $t\rightarrow 0$ and $t\rightarrow \infty$. Since $$
\int_0^{\frac{3\delta}{t^{1/2}}} \xi^{n-1}e^{-\xi^2} d\xi= \int_0^{\infty} \xi^{n-1}e^{-\xi^2} d\xi -\int^\infty_{\frac{3\delta}{t^{1/2}}} \xi^{n-1}e^{-\xi^2} d\xi,
$$
by some elementary inequality, it is a constant $d_n'=d_n \int_0^{\infty} \xi^{n-1}e^{-\xi^2} d\xi$ plus 
an exponentially decaying term as $t\rightarrow 0$ (or one could just invoke the known asymptotic for the complimentary error function for large argument). 
On the other hand, it is also straightforward to see that as $t \rightarrow \infty$, the second term is $O(t^{-\frac{n-1}{2}})$.

Set
\begin{eqnarray*}
\zeta_\delta(s) & = &\frac{1}{\Gamma(s)}\int_0^\infty t^{s-1} \left(\int_{V_{3\delta}} K_0(t, x, x) dx\right) dt \\ 
&  = & \frac{1}{\Gamma(s)}\left[ \int_0^1 t^{s-1} \left(\int_{V_{3\delta}} K_0(t, x, x) dx\right) dt  + \int_1^\infty t^{s-1} \left(\int_{V_{3\delta}} K_0(t, x, x) dx\right) dt \right],
\end{eqnarray*}
where the first term is defined through analytic continuation from a region where the real part of $s$ is sufficiently large, 
whereas the second term defined through analytic continuation from a region where the real part of $s$ is sufficiently negative. Therefore
\begin{eqnarray*}
 \frac{1}{\Gamma(s)}\int_0^1 t^{s-1} \left(\int_{V_{3\delta}} K_0(t, x, x) dx\right) dt & = & \frac{1}{\Gamma(s)} \frac{c_n\delta^n}{s-n/2} + \frac{d_n'}{\Gamma(s+1)}  \\
& & -  \frac{d_n}{\Gamma(s)}\int_0^1 t^{s-1}\int^\infty_{\frac{3\delta}{t^{1/2}}} \xi^{n-1}e^{-\xi^2} d\xi dt  ,
\end{eqnarray*}
and
\begin{eqnarray*}
\frac{1}{\Gamma(s)}\int_1^\infty t^{s-1} \left(\int_{V_{3\delta}} K_0(t, x, x) dx\right) dt & = &  - \frac{1}{\Gamma(s)} \frac{c_n\delta^n}{s-n/2} \\
& & + \frac{d_n}{\Gamma(s)}\int_0^1 t^{s-1}\int_0^{\frac{3\delta}{t^{1/2}}} \xi^{n-1}e^{-\xi^2} d\xi dt.
\end{eqnarray*}

Thus,
$$
\zeta_\delta'(0)
= 
-d_n'\Gamma'(1) -d_n \int_0^1 t^{-1}\int^\infty_{\frac{3\delta}{t^{1/2}}}\xi^{n-1}e^{-\xi^2} d\xi dt  + d_n \int_0^1 t^{-1}\int^{\frac{3\delta}{t^{1/2}}}_0 \xi^{n-1}e^{-\xi^2} d\xi dt .
$$
By a simple change of integration we arrive at
$$
\zeta_\delta'(0)= -d_n'\Gamma'(1) -d_n \int_{3\delta}^\infty 2\ln \frac{3\delta}{\xi} \xi^{n-1}e^{-\xi^2} d\xi + d_n \int_0^{3\delta} 2\ln \frac{3\delta}{\xi} \xi^{n-1}e^{-\xi^2} d\xi.
$$

This has a logarithmic divergence ($2d_n' \ln 3\delta$) as $\delta \rightarrow 0$, but
$$
\ln \tau(V_{3\delta}, V_{\infty}, \gamma_0)= - k\zeta_\delta'(0)\sum_{q=0}^n (-1)^q q \binom{n}{q} =0
$$
by combinatorial formula since $n\geq 2$ (in fact equal to $2$ in this case).
The proof for the partial equivariant torsion is almost the same. We just need to insert the action of the involution $\theta$ into the heat kernel, which will result in only the $d_n$ terms similar to the above formulas. 
\end{pf}

\begin{corollary}
\label{cor:double:limit:partial:torsion}
We have
$$
\lim_{\delta \rightarrow 0}\lim_{\epsilon \rightarrow 0} \ln \tau(\widetilde{X}_{3\delta}, \widetilde{X}, \gamma_{\epsilon,\delta})= \ln \tau(X, \gamma_0), \ \ \ 
$$
$$
\lim_{\delta \rightarrow 0}\lim_{\epsilon \rightarrow 0} \ln \tau_{{\bf Z}_{2}}(\widetilde{X}_{3\delta}, \widetilde{X}, \gamma_{\epsilon,\delta})(\iota)
=\ln \tau_{{\bf Z}_{2}}(X, \gamma_0)(\iota).
$$
\end{corollary}

\begin{pf}
Since $\ln \tau(X, \gamma_0)=\ln \tau(X_{3\delta},X, \gamma_0)+\ln \tau(V_{3\delta},X,\gamma_0)$ and
$\ln \tau_{{\bf Z}_{2}}(X, \gamma_0)(\iota)=\ln \tau_{{\bf Z}_{2}}(X_{3\delta},X, \gamma_0)(\iota)+\ln \tau_{{\bf Z}_{2}}(V_{3\delta},X,\gamma_0)(\iota)$
by \eqref{aopat}, we get by Theorem~\ref{lopat3} 
$$
\lim_{\delta\to0}\ln \tau(X_{3\delta},X, \gamma_0)=\ln \tau(X, \gamma_0),
\qquad
\lim_{\delta\to0}\ln \tau_{{\bf Z}_{2}}(X_{3\delta},X, \gamma_0)(\iota)=\ln \tau_{{\bf Z}_{2}}(X, \gamma_0)(\iota),
$$ 
which, together with Theorem~\ref{lopat1}, yields the result. 
\end{pf}

\subsection
{Limit of partial analytic torsion III}
\par
On the other hand, we have

\begin{theorem} \label{lopat2}
The following equalities hold:
$$
\lim_{\delta \rightarrow 0}\lim_{\epsilon \rightarrow 0} \ln \tau(\widetilde{V}_{3\delta}, \widetilde{X}, \gamma_{\epsilon, \delta})
=
k\,\ln C_0^{\rm EH}(\rho), 
$$
$$
\lim_{\delta \rightarrow 0}\lim_{\epsilon \rightarrow 0} 
\ln \left[\epsilon^{k/3}\tau_{{\bf Z}_{2}}(\widetilde{V}_{3\delta}, \widetilde{X}, \gamma_{\epsilon, \delta})(\theta)\right]=k\,\ln C_1^{\rm EH}(\rho),
$$
where the constants $C_0^{\rm EH}(\rho)$, $C_1^{\rm EH}(\rho)$
depend only on 
the cut-off function $\rho$.
\end{theorem}

At this stage, the constants $C_0^{\rm EH}(\rho)$, $C_1^{\rm EH}(\rho)$ may depend on $\rho$. The fact that they are  independent of 
$\rho$ will be postponed to the next subsection.

\subsubsection
{An integral expression of $\tau(\widetilde{V}_{3\delta}, \widetilde{X}, \gamma_{\epsilon, \delta})$ and
$\tau_{{\bf Z}_{2}}(\widetilde{V}_{3\delta}, \widetilde{X}, \gamma_{\epsilon, \delta})(\theta)$}
\par
For the proof of Theorem~\ref{lopat2}, as before, we compute
\begin{eqnarray*}
\ln \tau(\widetilde{V}_{3\delta}, \widetilde{X}, \gamma_{\epsilon, \delta}) 
& = & - \int_1^{\infty} t^{-1} \int_{\widetilde{V}_{3\delta}} {\rm tr}_s(N K_{\epsilon, \delta}(t, x, x)P^{\perp}_{\epsilon, \delta}(x, x)) \, dx \, dt \\
& & -  \int_0^{1} t^{-1} \int_{\widetilde{V}_{3\delta}}\left[ {\rm tr}_s(N K_{\epsilon, \delta}(t, x, x))- \sum_{i=0}^{2} a_i^{\epsilon, \delta}(x)\, t^{i-2}\right] \, dx \, dt \\
& & - \sum_{i=0}^{1} \int_{\widetilde{V}_{3\delta}} \frac{a_i^{\epsilon, \delta}(x)}{ i-2} \, dx 
+ \Gamma'(1) \int_{\widetilde{V}_{3\delta}}\left[ a_2^{\epsilon, \delta}(x) -  {\rm tr}_s(N P_{\epsilon, \delta}(x, x))\right] \, dx
\end{eqnarray*}
and
\begin{eqnarray*}
\ln \tau_{{\bf Z}_{2}}(\widetilde{V}_{3\delta}, \widetilde{X}, \gamma_{\epsilon, \delta}) 
& = & - \int_1^{\infty} t^{-1} \int_{\widetilde{V}_{3\delta}} {\rm tr}_s(N K_{\epsilon, \delta}(t, x, \theta(x))P^{\perp}_{\epsilon, \delta}(x, \theta(x))) \, dx \, dt \\
& & -  \int_0^{1} \frac{dt}{t} \left[ \int_{\widetilde{V}_{3\delta}}{\rm tr}_s(N K_{\epsilon, \delta}(t, x, \theta(x)))\,dx - \sum_{i=0}^{1}t^{i-1} \int_{E}b_i(z)\, dz\right] \\
& & +  \int_{E} b_0(z) \, dz 
+ \Gamma'(1) \left[ \int_{E}b_1(z)\,dz - \int_{\widetilde{V}_{3\delta}} {\rm tr}_s(N P_{\epsilon, \delta}(x, \theta(x)))\, dx\right].
\end{eqnarray*}
We study the behavior of each term in the right hand side as $\epsilon\to0$ and $\delta\to0$. For this, we set
$$
\begin{aligned}
I(\epsilon,\delta;\rho)
&:=
- \int_{0}^{1}\frac{dt}{t}\int_{\widetilde{V}(3\delta)}\left[ {\rm tr}_{s}(N K_{\epsilon, \delta,\infty}(t, x, x))- \sum_{i=0}^{2} a_{i}^{\epsilon, \delta}(x)\, t^{i-2}\right] \, dx
\\
&\qquad
- \sum_{i=0}^{1} \int_{\widetilde{V}(3\delta)} \frac{a_{i}^{\epsilon, \delta}(x)}{ i-2} \, dx 
+ \Gamma'(1) \int_{\widetilde{V}(3\delta)}a_{2}^{\epsilon, \delta}(x) \, dx,
\end{aligned}
$$
$$
\begin{aligned}
J(\epsilon,\delta;\rho)
&:=
- \int_{0}^{1}\frac{dt}{t}
\left[ 
\int_{\widetilde{V}(3\delta)}{\rm tr}_{s}(N K_{\epsilon, \delta,\infty}(t, x, \theta(x)))\, dx - \sum_{i=0}^{1}t^{i-1}\int_{E}b_{i}^{\epsilon}(z)\,dz
\right] 
\\
&\qquad
+ \int_{E} \frac{b_{0}^{\epsilon}(z)}{2} \, dz + \Gamma'(1) \int_{E}b_{1}^{\epsilon}(z) \, dz.
\end{aligned}
$$
Since $K_{\epsilon, \delta,\infty}(t, x, y)=\oplus_{q}K_{\epsilon, \delta,\infty}^{q}(t, x, y)$ is the heat kernel of $(T^{*}{\bf P}^{1},\kappa_{\epsilon,\delta})$,
$I(\epsilon,\delta;\rho)$ and $J(\epsilon,\delta;\rho)$ depend only on $\epsilon,\delta\in(0,1]$ with $\epsilon\delta^{-2}\leq\epsilon(\rho)$ and the cut-off function $\rho$.
Since $\widetilde{V}_{3\delta}$ is a $k$-copies of $\widetilde{V}(3\delta)$, we have
$$
\begin{aligned}
\ln \tau(\widetilde{V}_{3\delta}, \widetilde{X}, \gamma_{\epsilon, \delta})
&=
-\int_1^{\infty} \frac{dt}{t} \int_{\widetilde{V}_{3\delta}} {\rm tr}_s(N K_{\epsilon, \delta}(t, x, x)P^{\perp}_{\epsilon, \delta}(x, x)) \, dx
\\
&\quad
-\int_{0}^{1}\frac{dt}{t}\int_{\widetilde{V}_{3\delta}}{\rm tr}_{s}\{N K_{\epsilon,\delta}(t,x,x)-N K_{\epsilon,\delta,\infty}(t,x,x)\}dx
\\
&\quad
-\Gamma'(1)\int_{\widetilde{V}_{3\delta}}{\rm tr}_{s}(N P_{\epsilon,\delta}(x,x))\,dx+k\cdot I(\epsilon,\delta;\rho)
\end{aligned}
$$
and similarly
$$
\begin{aligned}
\ln \tau_{{\bf Z}_{2}}(\widetilde{V}_{3\delta}, \widetilde{X}, \gamma_{\epsilon, \delta})
&=
-\int_1^{\infty} \frac{dt}{t} \int_{\widetilde{V}_{3\delta}} {\rm tr}_s(N K_{\epsilon, \delta}(t, x, \theta(x))P^{\perp}_{\epsilon, \delta}(x, \theta(x))) \, dx
\\
&\quad
-\int_{0}^{1}\frac{dt}{t}\int_{\widetilde{V}_{3\delta}}{\rm tr}_{s}\{N K_{\epsilon,\delta}(t,x,\theta(x))-N K_{\epsilon,\delta,\infty}(t,x,\theta(x))\}dx
\\
&\quad
-\Gamma'(1)\int_{\widetilde{V}_{3\delta}}{\rm tr}_{s}\{N P_{\epsilon,\delta}(x,\theta(x))\}\,dx+k\cdot J(\epsilon,\delta;\rho).
\end{aligned}
$$

\subsubsection
{Limit of the first integral}
\par

\begin{proposition}
\label{prop:first:term}
The following equality holds:
$$
\lim_{\delta\to0}\lim_{\epsilon\to0}
\int_1^{\infty} \frac{dt}{t} \int_{\widetilde{V}_{3\delta}} {\rm tr}_s(N K_{\epsilon, \delta}(t, x, x)P^{\perp}_{\epsilon, \delta}(x, x)) \, dx=0.
$$
The same is true for the first integral in the expression of $\ln \tau_{{\bf Z}_{2}}(\widetilde{V}_{3\delta}, \widetilde{X}, \gamma_{\epsilon, \delta})$.
\end{proposition}

\begin{pf}
Let $\nu>0$ be arbitrary.
As in the proof of Theorem~\ref{lopat1} Step 2,
there is $T=T(\nu)>0$ depending only on $\nu$ such that
\begin{equation}
\label{eqn:integral:V3delta:[T,infty]}
\int_T^{\infty} t^{-1} \int_{\widetilde{V}_{3\delta}} \left|{\rm tr}_s(N K_{\epsilon, \delta}(t, x, x)P^{\perp}_{\epsilon, \delta}(x, x)) \right|\, dx \, dt 
< \nu
\end{equation}
for all $\epsilon,\delta \in (0,1]$ with $\epsilon \leq \min\{\epsilon(\rho)\delta^{2}, \delta^{4}\}$, which will be assumed throughout the proof.
By Theorem \ref{thm:heat:kernel:estimates:3},
\begin{equation}
\label{eqn:diff:heat:kernel:near:exc}
\left| \int_1^{T} t^{-1} \int_{\widetilde{V}_{3\delta}} {\rm tr}_s[N\{ K_{\epsilon, \delta}(t, x, x)) - K_{\epsilon, \delta,\infty}(t, x, x)\}] \, dx \, dt  \right| 
\leq 
C(T)\, {\rm vol}(\widetilde{V}_{3\delta}),
\end{equation}
where $C(T)$ is a constant depending only on $T$. 
By \eqref{eqn:harmonic:projection}, we get
\begin{equation}
\label{eqn:harmonic:proj:near:exc:div}
\int_{\widetilde{V}_{3\delta}}{\rm tr}_{s}[N(P_{\epsilon,\delta}(x,x)]\,dx
=
\int_{\widetilde{V}_{3\delta}}\frac{2\eta\wedge\bar{\eta}}{\|\eta\|_{L^{2}}^{2}}
\leq
2\frac{\|\eta\wedge\bar{\eta}/\gamma_{0}^{2}\|_{L^{\infty}}}{\|\eta\|_{L^{2}}}{\rm Vol}(\widetilde{V}_{3\delta}).
\end{equation}
By \eqref{eqn:heat:ker:harmonic:proj}, \eqref{eqn:diff:heat:kernel:near:exc}, \eqref{eqn:harmonic:proj:near:exc:div}, we get
\begin{equation}
\label{eqn:integral:[1,T]:near:exc:div}
\begin{aligned}
\,&
 \left| 
 \int_1^{T} t^{-1} \int_{\widetilde{V}_{3\delta}} {\rm tr}_s[N\{ K_{\epsilon, \delta}(t, x, x)P^{\perp}_{\epsilon, \delta}(x, x)-K_{\epsilon, \delta,\infty}(t, x,x)\}] \, dx \, dt 
 \right| 
 \\
& \leq  
 \{C(T)+2\frac{\|\eta\wedge\bar{\eta}/\gamma_{0}^{2}\|_{L^{\infty}}}{\|\eta\|_{L^{2}}}\log T\}\, {\rm vol}(\widetilde{V}_{3\delta}).
 \end{aligned}
\end{equation}
By Proposition~\ref{prop:heat:trace:noncompact}, there is a constant $A>0$ such that
\begin{equation}
\label{eqn:integral:[1,T]:noncompact}
\int_1^{T} t^{-1} \int_{\widetilde{V}_{3\delta}} \left| K_{\epsilon, \delta,\infty}(t, x,x) \right|\, dx \, dt
\leq 
Ae^{C(\epsilon\delta^{-4}+1)T}\,\log T\cdot{\rm vol}(\widetilde{V}_{3\delta})
\end{equation}
for all $\epsilon,\delta\in(0,1]$ with $\epsilon/\delta^{-2}\leq\epsilon(\rho)$.
By \eqref{eqn:integral:[1,T]:near:exc:div}, \eqref{eqn:integral:[1,T]:noncompact}, we get
\begin{equation}
\label{eqn:integral:V3delta:[1,T]}
\left| 
\int_1^{T} t^{-1} \int_{\widetilde{V}_{3\delta}} {\rm tr}_s[N K_{\epsilon, \delta}(t, x, x)P^{\perp}_{\epsilon, \delta}(x, x)] \, dx \, dt 
\right| 
\leq  
\widetilde{C}(T)\, {\rm vol}(\widetilde{V}_{3\delta}),
\end{equation}
where $\widetilde{C}(T)=C(T)+
(2\frac{\|\eta\wedge\bar{\eta}/\gamma_{0}^{2}\|_{L^{\infty}}}{\|\eta\|_{L^{2}}}+Ae^{2CT})\log T$.
Since $\nu>0$ can be chosen arbitrarily small, by taking into account that ${\rm vol}(\widetilde{V}_{3\delta}) $ goes to zero as $\delta \rightarrow 0$, 
the result follows from \eqref{eqn:integral:V3delta:[T,infty]}, \eqref{eqn:integral:V3delta:[1,T]}.
\end{pf}

\subsubsection
{Limit of the second integral}
\par

\begin{proposition}
\label{prop:second:term}
The following equality holds:
$$
\lim_{\delta\to0}\lim_{\epsilon\to0}
\int_{0}^{1}\frac{dt}{t}\int_{\widetilde{V}_{3\delta}}{\rm tr}_{s}\{N K_{\epsilon,\delta}(t,x,x)-N K_{\epsilon,\delta,\infty}(t,x,x)\}dx=0.
$$
The same is true for the second integral in the expression of $\ln \tau_{{\bf Z}_{2}}(\widetilde{V}_{3\delta}, \widetilde{X}, \gamma_{\epsilon, \delta})$.
\end{proposition}

\begin{pf}
The proof is the same as above, using the estimate of Theorem~\ref{thm:heat:kernel:estimates:3}.  
Indeed, we have, for all $\epsilon,\delta \in (0,1]$ with $\epsilon \leq \min\{\epsilon(\rho)\delta^{2}, \delta^{4}\}$, there is a constant $C>0$ such that 
$$
|{\rm tr}_{s}\{N K_{\epsilon,\delta}(t,x,x)-N K_{\epsilon,\delta,\infty}(t,x,x)\}|\leq C\,t
$$ 
for all $(x,t)\in \widetilde{V}_{3\delta}\times(0,1]$. Hence
\begin{equation}
\label{eqn:estimate:integral:5}
\left|
\int_{0}^{1}\frac{dt}{t}\int_{\widetilde{V}_{3\delta}}{\rm tr}_{s}\{N K_{\epsilon,\delta}(t,x,x)-N K_{\epsilon,\delta,\infty}(t,x,x)\}dx
\right|
\leq
C\,{\rm Vol}(\widetilde{V}_{3\delta},\gamma_{\epsilon, \delta}).
\end{equation}
By 
 the fact that 
 $$\lim_{\delta\to0}\lim_{\epsilon\to0}{\rm Vol}(\widetilde{V}_{3\delta},\gamma_{\epsilon, \delta})= \lim_{\delta\to0}{\rm Vol}(\widetilde{V}_{3\delta},\gamma_{0})=0,
 $$
  we get the result.
\end{pf}

\subsubsection
{Proof of Theorem~\ref{lopat2} }
\par
By \eqref{eqn:harmonic:projection}, we get  
$$
\lim_{\delta\to0}\lim_{\epsilon\to0}\int_{\widetilde{V}_{3\delta}}{\rm tr}_{s}(N P_{\epsilon,\delta}(x,x))\,dx
=
\lim_{\delta\to0}\lim_{\epsilon\to0}\int_{\widetilde{V}_{3\delta}}{\rm tr}_{s}(N P_{\epsilon,\delta}(x,\theta(x)))\,dx
=0.
$$
From Propositions~\ref{prop:first:term} and \ref{prop:second:term}, it follows that
$$
\lim_{\delta\to0}\lim_{\epsilon\to0}
\ln \tau(\widetilde{V}_{3\delta}, \widetilde{X}, \gamma_{\epsilon, \delta})
=
k\,\lim_{\delta\to0}\lim_{\epsilon\to0}I(\epsilon,\delta;\rho),
$$
$$
\lim_{\delta\to0}\lim_{\epsilon\to0}
\ln\left[\epsilon^{k/3} \tau_{{\bf Z}_{2}}(\widetilde{V}_{3\delta}, \widetilde{X}, \gamma_{\epsilon, \delta})\right]
=
k\,\lim_{\delta\to0}\lim_{\epsilon\to0}\left[J(\epsilon,\delta;\rho)+\frac{1}{3}\ln\epsilon\right].
$$
Since the right hand side depend only on the choice of $\rho$, we get the result by setting
$$
\ln C_{0}^{\rm EH}(\rho):=\lim_{\delta\to0}\lim_{\epsilon\to0}I(\epsilon,\delta;\rho),
\qquad
\ln C_{1}^{\rm EH}(\rho):=\lim_{\delta\to0}\lim_{\epsilon\to0}\left[J(\epsilon,\delta;\rho)+\frac{1}{3}\ln\epsilon\right].
$$
This completes the proof, provided that these double limits exist. This will be addressed in what follows.
\qed

\begin{remark} 
$C_0^{\rm EH}(\rho)$, respectively $C_1^{\rm EH}(\rho)$, is  renormalized (resp. equivariant) analytic torsion for the asymptotically conical space 
$\widetilde{V}(\infty)=(T^{*}{\bf P}^{1},\gamma^{\rm EH})$.
\end{remark}

\subsection
{Proof of Theorem \ref{loat}}
\par
Since
$$
\ln \tau(\widetilde{X}, \gamma_{\epsilon, \delta})
=
\ln \tau(\widetilde{X}_{3\delta}, \widetilde{X}, \gamma_{\epsilon, \delta})
+
\ln \tau(\widetilde{V}_{3\delta}, \widetilde{X}, \gamma_{\epsilon, \delta})
$$
and
$$
\ln \tau_{{\bf Z}_{2}}(\widetilde{X}, \gamma_{\epsilon, \delta})
=
\ln \tau_{{\bf Z}_{2}}(\widetilde{X}_{3\delta}, \widetilde{X}, \gamma_{\epsilon, \delta})
+
\ln \tau_{{\bf Z}_{2}}(\widetilde{V}_{3\delta}, \widetilde{X}, \gamma_{\epsilon, \delta})
$$
by the definition of partial (equivariant) analytic torsion, we get by Corollary~\ref{cor:double:limit:partial:torsion} and Theorem \ref{lopat2}
\begin{equation}
\label{eqn:limit:torsion}
\lim_{\delta\to0}\lim_{\epsilon\to0}
\ln \tau(\widetilde{X}, \gamma_{\epsilon, \delta})
=
\ln \tau(X, \gamma_{0})
+
k\,\ln C_{0}^{\rm EH}(\rho),
\end{equation}
\begin{equation}
\label{eqn:limit:equiv:torsion}
\lim_{\delta\to0}\lim_{\epsilon\to0}
\ln\left[\epsilon^{k/3} \tau_{{\bf Z}_{2}}(\widetilde{X}, \gamma_{\epsilon, \delta})\right]
=
\ln \tau_{{\bf Z}_{2}}(X, \gamma_{0})
+
k\,\ln C_{1}^{\rm EH}(\rho).
\end{equation}
As the double limits on the left hand side of (\ref{eqn:limit:torsion}), (\ref{eqn:limit:equiv:torsion}) exist by virtue of Corollary~\ref{cor:double:limit:torsion:K3} and Proposition~\ref{prop:limit:equivariant:torsion}, so do the double limits in defining $\ln C_{0}^{\rm EH}(\rho)$ and $\ln C_{1}^{\rm EH}(\rho)$. 

On the other hand, again by Corollary~\ref{cor:double:limit:torsion:K3} and Proposition~\ref{prop:limit:equivariant:torsion}, the double limits $\lim_{\delta\to0}\lim_{\epsilon\to0}\ln \tau(\widetilde{X}, \gamma_{\epsilon, \delta})$
and
$\lim_{\delta\to0}\lim_{\epsilon\to0}\ln\left[\epsilon^{k/3} \tau_{{\bf Z}_{2}}(\widetilde{X}, \gamma_{\epsilon, \delta})\right]$ 
are independent of the choice of $\rho$.
Hence $C_{0}^{\rm EH}(\rho)$ and $C_{1}^{\rm EH}(\rho)$ in \eqref{eqn:limit:torsion}, \eqref{eqn:limit:equiv:torsion} are in fact independent of $\rho$.
This completes the proof of Theorem \ref{loat}.
\qed


\section
{A holomorphic torsion invariant of log-Enriques surfaces}
\par
In this section, we introduce a holomorphic torsion invariant of log-Enriques surfaces and give its explicit formula as a function on the moduli space.

\subsection
{A construction of invariant} \label{inv}
\par

\begin{theorem}
\label{thm:equiv:torsion:K3}
There is a constant $C(k)$ depending only on $k=\#{\rm Sing}(Y)$ with
$$
\begin{aligned}
\tau_{M}(\widetilde{X},\theta)
&=
C(k){\rm Vol}(Y,\gamma_{0})^{\frac{4-k}{4}}\,\tau(Y,\gamma_{0})^{2}
\times
\prod_{{\frak p}\in{\rm Sing}(X)}\left\{|f_{\frak p}(0)|^{2}\frac{{\rm Vol}(Y,\gamma_{0})}{\|\eta\|_{L^{2}(Y)}^{2}}\right\}^{\frac{5}{16}}
\\
&\quad\times
\exp
\left(
\frac{1}{12}\int_{Y}
\log\left\{\frac{\eta\wedge\bar{\eta}}{\gamma_{0}^{2}/2!}\cdot\frac{{\rm Vol}(Y,\gamma_{0})}{\|\eta\|_{L^{2}(Y)}^{2}}\right\}\,
c_{2}(Y,\gamma_{0})
\right).
\end{aligned}
$$
\end{theorem}

\begin{pf}
Since $M_{k}^{\perp}\cong\Lambda_{k}(2)$, we have $\frac{14-r(H^{2}(\widetilde{X},{\bf Z})^{+})}{4}=\frac{4-k}{4}$. 
By its independence of the choice of $\theta$-invariant K\"ahler metric on $\widetilde{X}$, $\tau_{M}(\widetilde{X},\theta)$ is given by
$$
\begin{aligned}
\,&
\lim_{\delta\to0}\lim_{\epsilon\to0}
\tau(\widetilde{X},\gamma_{\epsilon,\delta})\tau_{{\bf Z}_{2}}(\widetilde{X},\gamma_{\epsilon,\delta})(\theta)\,
{\rm Vol}(\widetilde{X},\gamma_{\epsilon,\delta})^{\frac{4-k}{4}}
{\rm Vol}(\widetilde{X}^{\theta},\gamma_{\epsilon,\delta}|_{\widetilde{X}^{\theta}})
\tau(\widetilde{X}^{\theta},\gamma_{\epsilon,\delta}|_{\widetilde{X}^{\theta}})
\\
&\quad\times
A_{M}(\widetilde{X},\theta,\gamma_{\epsilon,\delta})\,
\exp
\left(
\frac{1}{24}\int_{\widetilde{X}}
\log\left\{
\frac{\eta\wedge\bar{\eta}}{\gamma_{\epsilon,\delta}^{2}/2!}
\cdot
\frac{{\rm Vol}(\widetilde{X},\gamma_{\epsilon,\delta})}{\|\eta\|_{L^{2}(\widetilde{X})}^{2}}
\right\}\,
c_{2}(\widetilde{X},\gamma_{\epsilon,\delta})
\right)
\\
&=
\lim_{\delta\to0}\lim_{\epsilon\to0}
\{\epsilon^{\frac{k}{3}}\tau(\widetilde{X},\gamma_{\epsilon,\delta})\tau_{{\bf Z}_{2}}(\widetilde{X},\gamma_{\epsilon,\delta})(\theta)\}
\\
&\quad\times
\lim_{\delta\to0}\lim_{\epsilon\to0}
\prod_{{\frak p}\in{\rm Sing}(X)}\epsilon^{-\frac{1}{3}}{\rm Vol}(E_{\frak p},\gamma_{\epsilon,\delta}|_{E_{\frak p}})
\tau(E_{\frak p},\gamma_{\epsilon,\delta}|_{E_{\frak p}})
\times
\lim_{\delta\to0}\lim_{\epsilon\to0}A_{M}(\widetilde{X},\theta,\gamma_{\epsilon,\delta})
\\
&\quad\times
\lim_{\delta\to0}\lim_{\epsilon\to0}
\exp
\left(
\frac{1}{24}\int_{\widetilde{X}}
\log
\left\{
\frac{\eta\wedge\bar{\eta}}{\gamma_{\epsilon,\delta}^{2}/2!}
\cdot
\frac{{\rm Vol}(\widetilde{X},\gamma_{\epsilon,\delta})}{\|\eta\|_{L^{2}(\widetilde{X})}^{2}}
\right\}\,
c_{2}(\widetilde{X},\gamma_{\epsilon,\delta})
\right).
\end{aligned}
$$
By Propositions~\ref{prop:limit:2nd:Chern:form:2}, \ref{prop:limit:anomaly}, \ref{prop:torsion:fixed:curve},
Corollary~\ref{cor:double:limit:torsion:K3} and Theorem~\ref{loat}, we get
$$
\begin{aligned}
\tau_{M}(\widetilde{X},\theta)
&=
(C_{0}^{\rm EH}C_{1}^{\rm EH})^{k}\,\tau(X,\gamma_{0})\tau_{{\bf Z}_{2}}(X,\gamma_{0})(\iota)\,
\{2{\rm Vol}(Y,\gamma_{0})\}^{\frac{4-k}{4}}
\\
&\quad\times
\{{\rm Vol}({\bf P}^{1},\omega_{\rm FS})\tau({\bf P}^{1},\omega_{\rm FS})\}^{k}
\times
\prod_{{\frak p}\in{\rm Sing}(X)}\left\{|f_{\frak p}(0)|^{2}\frac{{\rm Vol}(X,\gamma_{0})}{\|\eta\|_{L^{2}(X)}^{2}}\right\}^{\frac{1}{4}+\frac{1}{16}}
\\
&\quad\times
\exp
\left[
\frac{1}{24}\int_{X}
\log\left\{\frac{\eta\wedge\bar{\eta}}{\gamma_{0}^{2}/2!}\cdot\frac{{\rm Vol}(X,\gamma_{0})}{\|\eta\|_{L^{2}(X)}^{2}}\right\}\,
c_{2}(X,\gamma_{0})
\right].
\end{aligned}
$$
Since
$$
\tau(Y,\gamma_{0})^{2}=\tau(X,\gamma_{0})\tau_{{\bf Z}_{2}}(X,\gamma_{0})(\iota),
\qquad
{\rm Vol}(X,\gamma_{0})/\|\eta\|_{L^{2}(X)}^{2}={\rm Vol}(Y,\gamma_{0})/\|\eta\|_{L^{2}(Y)}^{2}
$$
and since $X$ is a double covering of $Y$, we get the result by setting
\begin{equation}
\label{eqn:constant}
C(k)=2\{2^{-\frac{1}{4}}C_{0}^{\rm EH}C_{1}^{\rm EH}{\rm Vol}({\bf P}^{1},\omega_{\rm FS})\tau({\bf P}^{1},\omega_{\rm FS})\}^{k}.
\end{equation}
This completes the proof.
\end{pf}

\begin{theorem}
\label{thm:anomaly:log:Enriques}
Let $\gamma$ be a K\"ahler form on $Y$ in the sense of orbifolds. Then the following equality holds:
$$
\frac{\tau(Y,\gamma){\rm Vol}(Y,\gamma)}{\tau(Y,\omega_{\eta}){\rm Vol}(Y,\omega_{\eta})}
=
\left\{
\prod_{{\frak p}\in{\rm Sing}(Y)}\left(\frac{\omega_{\eta}^{2}}{\gamma^{2}}\right)({\frak p})
\right\}^{-\frac{5}{32}}
\exp
\left\{
-\frac{1}{24}\int_{Y}
\log\left(\frac{\omega_{\eta}^{2}}{\gamma^{2}}\right)\,c_{2}(Y,\gamma)
\right\}.
$$
\end{theorem}

\begin{pf}
Let ${\frak p}\in{\rm Sing}(Y)$ and let $({\mathcal U}_{\frak p},0)\subset({\bf C}^{2},0)$ be an open subset which uniformizes the germ $(Y,{\frak p})$.
We have an isomorphism $(Y,{\frak p})\cong({\bf C}^{2}/\Gamma_{\frak p},0)$ of germs, where $\Gamma_{\frak p}={\bf Z}/4{\bf Z}=\langle i\rangle$,
such that $\omega_{\eta}$ and $\gamma$ lift to K\"ahler metrics on ${\mathcal U}_{\frak p}$.
Following Ma \cite{Ma05}, we define $Y^{\Sigma}$ as the union $Y^{\Sigma}:=Y^{i}\amalg Y^{i^{2}}\amalg Y^{i^{3}}$, where
$Y^{i^{\nu}}=\{{\frak p}^{i^{\nu}}\}_{{\frak p}\in{\rm Sing}(Y)}$ and the germ $(Y^{i^{\nu}},{\frak p}^{i^{\nu}})$ is equipped with orbifold structure
$(Y^{i^{\nu}},{\frak p}^{i^{\nu}})\cong({\bf C}^{2}/\langle i^{\nu}\rangle,0)$.
\par
Recall that the characteristic class ${\rm Td}^{\Sigma}(TY)$ supported on the singular locus of $Y$ appears in the Riemann-Roch theorem for orbifolds,
for which we refer the reader to e.g. \cite{Ma05}.
By the anomaly formula for Quillen metrics for orbifolds \cite{Ma05}, we get
\begin{equation}
\label{eqn:5:anomaly:0}
\begin{aligned}
\,&
\log
\left(
\frac{\tau(Y,\gamma){\rm Vol}(Y,\gamma)}{\tau(Y,\omega_{\eta}){\rm Vol}(Y,\omega_{\eta})}
\right)
=
\frac{1}{4}\int_{Y^{\Sigma}}\widetilde{\rm Td}^{\Sigma}(TY;\gamma,\omega_{\eta})
+
\frac{1}{24}\int_{Y}\widetilde{c_{1}c_{2}}(TY;\gamma,\omega_{\eta})
\\
&=
\frac{1}{4}\sum_{{\frak p}\in{\rm Sing}(Y)}\sum_{\nu=1}^{3}
\left[
\left(\widetilde{\frac{\rm Td}{\rm e}}\right)_{\nu/2}(TU_{\frak p};\gamma,\omega_{\eta})
\right]^{(0,0)}({\frak p})
+
\frac{1}{24}\int_{Y}\widetilde{c_{1}c_{2}}(TY;\gamma,\omega_{\eta})^{(2,2)}.
\end{aligned}
\end{equation}
Here, for $\theta\in{\bf R}$ and a square matrix $A$, we define
$\left({\frac{\rm Td}{\rm e}}\right)_{\theta}(A):=\det\left(\frac{I}{I-e^{\pi i\theta}A}\right)$
and $(\widetilde{{\rm Td}/{\rm e}})_{\theta}$ is the Bott-Chern secondary class associated to $\left({\rm Td}/{\rm e}\right)_{\theta}(A)$ 
such that for any holomorphic vector bundle $E$ and Hermitian metrics $h$, $h'$ on $E$
$$
-dd^{c}\left(\widetilde{\frac{\rm Td}{\rm e}}\right)_{\theta}(E;h,h')
=
\left(\frac{\rm Td}{\rm e}\right)_{\theta}\left(-\frac{1}{2\pi i}R(E,h)\right)
-
\left(\frac{\rm Td}{\rm e}\right)_{\theta}\left(-\frac{1}{2\pi i}R(E,h')\right).
$$
Similarly, $\widetilde{c_{1}c_{2}}$ is the Bott-Chern secondary class associated to the invariant polynomial $c_{1}(A)c_{2}(A)$
such that for any holomorphic vector bundle $E$ and Hermitian metrics $h$, $h'$ on $E$
$$
-dd^{c}\widetilde{c_{1}c_{2}}(E;h,h')=c_{1}(E,h)c_{2}(E,h)-c_{1}(E,h')c_{2}(E,h').
$$
\par
For $A={\rm diag}(\lambda_{1},\lambda_{2})$, we have
$$
\left(\frac{{\rm Td}}{{\rm e}}\right)_{\frac{\nu}{2}}(A)=\frac{1}{(1-i^{-\nu})^{2}}\{1-\frac{i^{-\nu}}{1-i^{-\nu}}c_{1}(A)+O(2)\}.
$$
Thus we get
\begin{equation}
\label{eqn:5:anomaly:1}
\begin{aligned}
\sum_{\nu=1}^{3}
\left[
\left(\widetilde{\frac{\rm Td}{\rm e}}\right)_{\frac{\nu}{2}}(TU_{\frak p};\gamma,\omega_{\eta})
\right]^{(0,0)}
({\frak p})
&=
-\sum_{\nu=1}^{3}\frac{i^{-\nu}}{(1-i^{-\nu})^{3}}\,\widetilde{c_{1}}(TU_{\frak p};\gamma,\omega_{\eta})({\frak p})
\\
&=
\frac{5}{8}\widetilde{c_{1}}(TU_{\frak p};\gamma,\omega_{\eta})({\frak p})
=
-\frac{5}{8}\log\left(\omega_{\eta}^{2}/\gamma^{2}\right)({\frak p}).
\end{aligned}
\end{equation}
\par
On the other hand, by the same computations as in \eqref{eqn:4:BottChern}, we get
\begin{equation}
\label{eqn:5:anomaly:2}
\widetilde{c_{1}c_{2}}(TY;\gamma,\omega_{\eta})^{(2,2)}=-\log\left(\omega_{\eta}^{2}/\gamma^{2}\right)\,c_{2}(TY,\gamma).
\end{equation}
Substituting \eqref{eqn:5:anomaly:1} and \eqref{eqn:5:anomaly:2} into \eqref{eqn:5:anomaly:0}, we get the result.
\end{pf}

\begin{theorem}
\label{thm:comparison:torsion:logEnriques}
For every Ricci-flat log-Enriques surface $(Y,\omega)$, one has
$$
{\rm Vol}(Y,\omega)^{\frac{4-k}{8}}\tau(Y,\omega)=C(k)^{-1}\,\tau_{M}(\widetilde{X},\theta)^{\frac{1}{2}},
$$
where $C(k)$ is the same constant as in Theorem~\ref{thm:equiv:torsion:K3}.
\end{theorem}

\begin{pf}
We put $\gamma=\gamma_{0}$ in Theorem~\ref{thm:anomaly:log:Enriques}.
Then we get by Theorem~\ref{thm:equiv:torsion:K3}
$$
\begin{aligned}
\,&
\tau(Y,\omega_{\eta}){\rm Vol}(Y,\omega_{\eta})
=
\tau(Y,\gamma_{0}){\rm Vol}(Y,\gamma_{0})^{\frac{4-k}{8}}{\rm Vol}(Y,\gamma_{0})^{\frac{4+k}{8}}
\\
&\qquad\qquad\qquad\qquad\quad\times
\{
\prod_{{\frak p}\in{\rm Sing}(X)}\left(\frac{\omega_{\eta}^{2}}{\gamma_{0}^{2}}\right)({\frak p})
\}^{\frac{5}{32}}
\exp
\left[
\frac{1}{24}\int_{Y}
\log\left(\frac{\omega_{\eta}^{2}}{\gamma_{0}^{2}}\right)\,c_{2}(Y,\gamma_{0})
\right]
\\
&=
C(k)^{-1}\tau_{M}(\widetilde{X},\theta)^{\frac{1}{2}}{\rm Vol}(Y,\gamma_{0})^{\frac{4+k}{8}}
\times
\prod_{{\frak p}\in{\rm Sing}(X)}\left\{|f_{\frak p}(0)|^{2}\frac{{\rm Vol}(Y,\gamma_{0})}{\|\eta\|_{L^{2}(Y)}^{2}}\right\}^{-\frac{5}{32}}
\\
&\quad\times
\exp
\left[
-\frac{1}{24}\int_{Y}
\log\left\{\frac{\eta\wedge\bar{\eta}}{\gamma_{0}^{2}/2!}\cdot\frac{{\rm Vol}(Y,\gamma_{0})}{\|\eta\|_{L^{2}(Y)}^{2}}\right\}\,
c_{2}(Y,\gamma_{0})
\right]
\\
&\quad\times
\{
\prod_{{\frak p}\in{\rm Sing}(X)}\left(\frac{\omega_{\eta}^{2}}{\gamma_{0}^{2}}\right)({\frak p})
\}^{\frac{5}{32}}
\exp
\left[
\frac{1}{24}\int_{Y}
\log\left(\frac{\eta\wedge\bar{\eta}}{\gamma_{0}^{2}/2!}\right)\,c_{2}(Y,\gamma_{0})
\right].
\end{aligned}
$$
Since $|f_{\frak p}(0)|^{2}=[\eta\wedge\bar{\eta}/(\gamma_{0}^{2}/2!)]({\frak p})=[\omega_{\eta}^{2}/\gamma_{0}^{2}]({\frak p})$, we get
$$
\begin{aligned}
\,&
\tau(Y,\omega_{\eta}){\rm Vol}(Y,\omega_{\eta})
\\
&=
C(k)^{-1}\tau_{M}(\widetilde{X},\theta)^{\frac{1}{2}}{\rm Vol}(Y,\gamma_{0})^{\frac{4+k}{8}}
\times
\prod_{{\frak p}\in{\rm Sing}(X)}\left(|f_{\frak p}(0)|^{2}\frac{{\rm Vol}(Y,\gamma_{0})}{{\rm Vol}(Y,\omega_{\eta})}\right)^{-\frac{5}{32}}
\\
&\quad\times
\{
\prod_{{\frak p}\in{\rm Sing}(X)}|f_{\frak p}(0)|^{2}
\}^{\frac{5}{32}}
\exp
\left[
-\frac{1}{24}\int_{Y}
\log\left\{\frac{{\rm Vol}(Y,\gamma_{0})}{{\rm Vol}(Y,\omega_{\eta})}\right\}\,
c_{2}(Y,\gamma_{0})
\right]
\\
&=
C(k)^{-1}\tau_{M}(\widetilde{X},\theta)^{\frac{1}{2}}{\rm Vol}(Y,\gamma_{0})^{\frac{4+k}{8}}
\left(\frac{{\rm Vol}(Y,\gamma_{0})}{{\rm Vol}(Y,\omega_{\eta})}\right)^{-\frac{5}{32}k}
\exp\left[
-\frac{16-k}{32}\log\left(\frac{{\rm Vol}(Y,\gamma_{0})}{{\rm Vol}(Y,\omega_{\eta})}\right)
\right]
\\
&=
C(k)^{-1}\tau_{M}(\widetilde{X},\theta)^{\frac{1}{2}}{\rm Vol}(Y,\omega_{\eta})^{\frac{4+k}{8}},
\end{aligned}
$$
where we used the second assertion of Proposition~\ref{prop:limit:2nd:Chern:form:1} to get the second equality.
This proves the result.
\end{pf}

\begin{theorem}
\label{thm:torsion:logEnriques}
Let $\gamma$ be a K\"ahler form on $Y$ in the sense of orbifolds and let $\varXi\in H^{0}(Y,K_{Y}^{\otimes2})\setminus\{0\}$
be a nowhere vanishing bicanonical form on $Y$. Then 
$$
\begin{aligned}
\tau_{k}(Y)
&:=
\tau(Y,\gamma){\rm Vol}(Y,\gamma)\,\|\varXi\|_{L^{1}(Y)}^{-\frac{4+k}{8}}
\left\{
\prod_{{\frak p}\in{\rm Sing}(Y)}\left(\frac{\gamma^{2}/2!}{|\varXi|}\right)({\frak p})
\right\}^{\frac{5}{32}}
\\
&\quad\times
\exp
\left[
\frac{1}{24}\int_{Y}
\log\left(\frac{|\varXi|}{\gamma^{2}/2!}\right)\,c_{2}(Y,\gamma)
\right]
\end{aligned}
$$
is independent of the choices of $\gamma$ and $\varXi$, where $|\varXi|:=\sqrt{\varXi\otimes\overline{\varXi}}$
is the Ricci-flat volume form on $Y$ induced by $\varXi$.
In fact,
$$
\tau_{k}(Y)=C(k)^{-1}\tau_{M}(\widetilde{X},\theta)^{\frac{1}{2}}.
$$
\end{theorem}

\begin{pf}
Let $\omega$ be a Ricci-flat K\"ahler form on $Y$ in the sense of orbifolds such that $\omega^{2}/2!=|\varXi|$.
Since ${\rm Vol}(Y,\omega)=\|\varXi\|_{L^{1}(Y)}$, we get by Theorem~\ref{thm:comparison:torsion:logEnriques}
\begin{equation}
\label{eqn:5:torsion:Y:1}
{\rm Vol}(Y,\omega)\tau(Y,\omega)
=
{\rm Vol}(Y,\omega)^{\frac{4+k}{8}}{\rm Vol}(Y,\omega)^{\frac{4-k}{8}}\tau(Y,\omega)
=
C(k)^{-1}\,\|\varXi\|_{L^{1}(Y)}^{\frac{4+k}{8}}\,\tau_{M}(\widetilde{X},\theta)^{\frac{1}{2}}.
\end{equation}
\par
Let $\xi\in H^{0}(\widetilde{X},K_{\widetilde{X}})$ be a nowhere vanishing holomorphic $2$-form on $\widetilde{X}$ such that
$(p\circ\pi)^{*}\varXi=\xi^{\otimes2}$. Since $\omega=\omega_{\xi}$, i.e., $\omega^{2}/2!=\xi\wedge\overline{\xi}=|\varXi|$,
we get by Theorem~\ref{thm:anomaly:log:Enriques}
\begin{equation}
\label{eqn:5:comparison:torsion:Y}
\frac{\tau(Y,\gamma){\rm Vol}(Y,\gamma)}{\tau(Y,\omega){\rm Vol}(Y,\omega)}
=
\left\{
\prod_{{\frak p}\in{\rm Sing}(Y)}\left(\frac{|\varXi|}{\gamma^{2}/2!}\right)({\frak p})
\right\}^{-\frac{5}{32}}
\exp
\left[
-\frac{1}{24}\int_{Y}
\log\left(\frac{|\varXi|}{\gamma^{2}/2!}\right)\,c_{2}(Y,\gamma)
\right].
\end{equation}
Comparing \eqref{eqn:5:torsion:Y:1} and \eqref{eqn:5:comparison:torsion:Y}, we get
\begin{equation}
\label{eqn:5:torsion:Y:2}
\begin{aligned}
\tau(Y,\gamma){\rm Vol}(Y,\gamma)
&=
C(k)^{-1}\tau_{M}(\widetilde{X},\theta)^{\frac{1}{2}}\,\|\varXi\|_{L^{1}(Y)}^{\frac{4+k}{8}}
\left\{
\prod_{{\frak p}\in{\rm Sing}(Y)}\left(\frac{|\varXi|}{\gamma^{2}/2!}\right)({\frak p})
\right\}^{-\frac{5}{32}}
\\
&\quad\times
\exp
\left[
-\frac{1}{24}\int_{Y}
\log\left(\frac{|\varXi|}{\gamma^{2}/2!}\right)\,c_{2}(Y,\gamma)
\right].
\end{aligned}
\end{equation}
From \eqref{eqn:5:torsion:Y:2}, we get $\tau_{k}(Y)=C(k)^{-1}\tau_{M}(\widetilde{X},\theta)^{1/2}$.
Since the right hand side is independent of the choices of $\gamma$ and $\varXi$, so is $\tau_{k}(Y)$.
This completes the proof.
\end{pf}

\subsection
{Del Pezzo surfaces and an explicit formula for the invariant $\tau_{k}$}
\par
In this subsection, we give an explicit formula for $\tau_{k}$ as an automorphic function on the K\"ahler moduli of Del Pezzo surfaces.
Let $1\leq k\leq9$. We define the unimodular Lorentzian lattices $L_{k}$  and ${\Bbb U}(-1)$ as
$$
L_{k}:=\begin{pmatrix}1&0\\0&-I_{9-k}\end{pmatrix}
\quad(k\not=2),
\qquad
L_{2} := \begin{pmatrix}1&0\\0&-1\end{pmatrix}\hbox{ or }\begin{pmatrix}0&1\\1&0\end{pmatrix},
$$
$$
{\Bbb U}(-1):=\begin{pmatrix}
0&-1
\\
-1&0
\end{pmatrix}.
$$
We fix an isometry of lattices $\Lambda_{k}\cong{\Bbb U}(-1)\oplus L_{k}$ and identify $\Lambda_{k}$ with ${\Bbb U}(-1)\oplus L_{k}$.
\par
Let $V$ be a Del Pezzo surface of degree $k$, i.e., 
$$
k = \deg V = \int_{V} c_{1}(V)^{2}.
$$
Then $V = {\rm Bl}_{9-k}( {\mathbf P}^{2} )$ is the blowing-up of ${\mathbf P}^{2}$ at $9-k$ points in general position when $k\not=8$. 
When $k=8$, $V\cong\Sigma_{0}$ or $\Sigma_{1}$, where
$\Sigma_{n} = {\mathbf P}( {\mathcal O}_{{\mathbf P}^{1}} \oplus {\mathcal O}_{{\mathbf P}^{1}}(n) )$ is the Hirzebruch surface. 
Notice that $\Sigma_{0}={\mathbf P}^{1} \times {\mathbf P}^{1}$ and $\Sigma_{1}={\rm Bl}_{1}( {\mathbf P}^{2} )$.
When $V\not\cong\Sigma_{0}$, $H^{2}(V, {\mathbf Z})$ endowed with the cup product pairing is isometric to $L_{k}$
by identifying $H, E_{1},\ldots,E_{9-k}$ with the standard basis of $L_{k}$, where $H \in H^{2}(V, {\mathbf Z})$ is the class obtained from
the hyperplane class of $H^{2}({\mathbf P}^{2}, {\mathbf Z})$ and $E_{i}$ $(i=1,\ldots,9-k)$ are the classes of exceptional divisors. 
Similarly, $H(V, {\mathbf Z})$ endowed with the Mukai pairing is isometric to $\Lambda_{k}$. In what follows, we identify $L_{k}$ (resp. $\Lambda_{k}$)
with $H^{2}(V, {\mathbf Z})$ (resp. $H(V, {\mathbf Z})$) in this way.
\par
Recall that the type IV domain $\Omega_{k}$ associated with $\Lambda_{k}$ was defined in Section 3.
We identify $\Omega_{H(V, {\mathbf Z})}$ with the tube domain 
$H^{2}(V, {\mathbf Z})\otimes{\bf R}+i\,{\mathcal C}_{H^{2}(V, {\mathbf Z})} \subset H^{2}(V, {\mathbf C})$ via the map
\begin{equation}
\label{eqn:7:tube:domain}
H^{2}(V, {\mathbf Z})\otimes{\bf R}+i\,{\mathcal C}_{H^{2}(V, {\mathbf Z})}
\ni y\to[\exp(y)]:=\left[\left(1, y, y^{2}/2\right)\right]\in
\Omega_{H(V, {\mathbf Z})},
\end{equation}
where ${\mathcal C}_{H^{2}(V, {\mathbf Z})}:=\{v\in H^{2}(V, {\mathbf R}); \, v^{2}>0 \}$ is the positive cone of $H^{2}(V, {\mathbf R})$.
\par
Let ${\mathcal K}_{V} \subset {\mathcal C}_{H^{2}(V, {\mathbf Z})}$ be the K\"ahler cone of $V$, i.e., the cone of $H^{2}(V, {\mathbf R})$ 
consisting of K\"ahler classes on $V$. 
Let ${\rm Eff}(V) \subset H^{2}(V, {\mathbf R})$ be the effective cone of $V$, i.e., the dual cone of the K\"ahler cone ${\mathcal K}_{V}$.

\begin{definition}
\label{def:Borcherds:product}
Define the infinite product $\Phi_{V}(z)$ on $H^{2}(V, {\mathbf Z})\otimes{\bf R}+i\,{\mathcal K}_{V}$ by
$$
\begin{aligned}
\Phi_{V}(z)
&:=
e^{\pi i\langle c_{1}(V), z\rangle}
\prod_{\alpha\in {\rm Eff}(V)}
(1-e^{2\pi i\langle\alpha,z\rangle})^{c_{k}^{(0)}(\alpha^{2})}
\\
&\qquad\times
\prod_{\beta\in {\rm Eff}(V),\,\,
\beta/2\equiv c_{1}(V)/2\mod H^{2}(V, {\mathbf Z})}
(1-e^{\pi i\langle\beta,z\rangle})^{c_{k}^{(1)}(\beta^{2}/4)},
\end{aligned}
$$
where $\{c_{k}^{(0)}(l)\}_{l\in{\bf Z}}$, $\{c_{k}^{(1)}(l)\}_{l\in{\bf Z}+k/4}$ are defined by the generating functions
$$
\sum_{l\in{\bf Z}}c_{k}^{(0)}(l)\,q^{l}
=
\frac{\eta(2\tau)^{8}\theta_{{\Bbb A}_{1}}(\tau)^{k}}{\eta(\tau)^{8}\eta(4\tau)^{8}},
\qquad
\sum_{l\in\frac{k}{4}+{\bf Z}}c_{k}^{(1)}(l)\,q^{l}
=
-8\,\frac{\eta(4\tau)^{8}\theta_{{\Bbb A}_{1}+1/2}(\tau)^{k}}{\eta(2\tau)^{16}}.
$$
Here $\theta_{{\Bbb A}_{1}+\epsilon/2}(\tau):=\sum_{n\in{\bf Z}}q^{(n+\epsilon/2)^{2}}$ and $\eta(\tau):=q^{1/24}\prod_{n>0}(1-q^{n})$.
\end{definition}

Let ${\mathcal C}_{H^{2}(V, {\mathbf Z})}^{+}$ be the connected component of ${\mathcal C}_{H^{2}(V, {\mathbf Z})}$ that contains ${\mathcal K}_{V}$
and let $\Omega_{H(V, {\mathbf Z})}^{+}$ be the component of $\Omega_{H(V, {\mathbf Z})}$ corresponding to 
$H^{2}(V, {\mathbf R}) + i\,{\mathcal C}_{H^{2}(V, {\mathbf Z})}^{+}$ via the isomorphism \eqref{eqn:7:tube:domain}.
By Borcherds \cite[Th.\,13.3]{Borcherds98} (cf. \cite{Yoshikawa09}), $\Phi_{V}(z)$ converges absolutely for those 
$z\in H^{2}(V, {\mathbf R}) + i\,{\mathcal K}_{V}$ with $\Im z\gg0$ and extends to an automorphic form on $\Omega_{H(V, {\mathbf Z})}^{+}$ 
for $O^{+}(H(V, {\mathbf Z}))$ of weight $\deg V+4$ with zero divisor
${\rm div}(\Phi_{V})=\sum_{d\in H(V, {\mathbf Z}),\,d^{2}=-1}d^{\perp}$
under the identification $H^{2}(V, {\mathbf R}) + i\,{\mathcal C}_{H^{2}(V, {\mathbf Z})}^{+} \cong \Omega_{H(V, {\mathbf Z})}^{+}$.
\par
Recently, an explicit Fourier series expansion of $\Phi_{V}(z)$ is discovered by Gritsenko \cite[Cor.\,5.1]{Gritsenko18}.
It is also remarkable that $\Phi_{V}$ is the denominator function of a generalized Kac-Moody algebra, whose real and imaginary simple roots
are explicitly given by the Fourier series expansion of $\Phi_{V}$ \cite[\S 6.2, Th.\,6.1 Eq.\,(6.1), (6.10)]{GritsenkoNikulin18}.
In this sense, the series of Borcherds products $\Phi_{V}$ associated to Del Pezzo surfaces is quite analogous to
the Borcherds $\Phi$-function of rank $10$.
\par
We define the Petersson norm of $\Phi_{V}(z)$ by
$$
\|\Phi_{V}(z)\|^{2}:=\langle\Im z,\Im z\rangle^{4+\deg V}|\Phi_{V}(z)|^{2},
$$
where $z\in H^{2}(V, {\mathbf R}) + i\,{\mathcal C}_{H^{2}(V, {\mathbf Z})}^{+}$.
Then $\|\Phi_{V}\|^{2}$ is an $O^{+}(H(V, {\mathbf Z}))$-invariant $C^{\infty}$ function on $\Omega_{H(V, {\mathbf Z})}^{+}$.
Hence $\|\Phi_{V}\|^{2}$ is identified with a $C^{\infty}$ function on ${\mathcal M}_{\deg V}$ in the sense of orbifolds.

\begin{theorem}
Let $1\leq k\leq9$. There exists a constant $\widetilde{C}(k)>0$ depending only on $k$ such that
for every $2$-elementary $K3$ surface $(\widetilde{X},\theta)$ of type $M_{k}:=\Lambda_{k}(2)^{\perp}$,
$$
\tau_{M_{k}}(\widetilde{X},\theta)=\widetilde{C}(k)\,\| \Phi_{V}(\overline{\varpi}(\widetilde{X},\theta)) \|^{-1/2},
$$
where $k = \deg V$.
\end{theorem}

\begin{pf}
See \cite[Th.\,4.2 (1)]{Yoshikawa09} and \cite[Th.\,0.1]{Yoshikawa13}.
\end{pf}

\begin{theorem}
\label{thm:torsion:Bor:prod}
Let $1\leq k\leq9$. Then there exists a constant $C_{k}>0$ depending only on $k$ such that
for every good log-Enriques surface $Y$ with $\#{\rm Sing}(Y) = \deg V$,
$$
\tau_{\deg V}(Y) = C_{\deg V} \left\| \Phi_{V}(\overline{\varpi}(Y)) \right\|^{-1/4}.
$$
\end{theorem}

\begin{pf}
We set $k = \deg V$. When $k=2$, we define $V=\Sigma_{0}$ when $Y$ is of even type and $V=\Sigma_{1}$ when $Y$ is of odd type.
Let $(\widetilde{X},\theta)$ be the $2$-elementary $K3$ surface of type $M_{k}$ associated to $Y$.
By the definition of the period of $Y$, we have $\overline{\varpi}(Y)=\overline{\varpi}(\widetilde{X},\theta)$.
Hence
\begin{equation}
\label{eqn:7:comparison:automorphic:form}
\left\|\Phi_{V}\left(\overline{\varpi}(Y)\right)\right\|=\|\Phi_{V}(\overline{\varpi}(\widetilde{X},\theta))\|.
\end{equation}
By Theorems 5.11 and 7.2 and \eqref{eqn:7:comparison:automorphic:form}, we get
\begin{equation}
\label{eqn:7:torsion:Y}
\begin{aligned}
\tau_{k}(Y)
&=
C(k)^{-1}\tau_{M_{k}}(\widetilde{X},\theta)^{1/2}
=
C(k)^{-1}\widetilde{C}(k)\,\|\Phi_{V}(\overline{\varpi}(\widetilde{X},\theta))\|^{-1/4}
\\
&=
C(k)^{-1}\widetilde{C}(k)\,\left\|\Phi_{V}\left(\overline{\varpi}(Y)\right)\right\|^{-1/4}.
\end{aligned}
\end{equation}
Setting $C_{k}:=C(k)^{-1}\widetilde{C}(k)$ in \eqref{eqn:7:torsion:Y}, we get the result.
\end{pf}

\subsection
{The quasi-pullback of  $\Phi_{V}$}
\par
Let $\pi\colon \widetilde{V}:={\rm Bl}_{p}(V)\to V$ be the blow-up of $V$ at $p$ and let $E:=\pi^{-1}(p)$ be the exceptional curve of $\pi$. 
Then we have a map of cohomologies $\pi^{*}\colon H(V,{\mathbf Z})\to H(\widetilde{V},{\mathbf Z})$, which induces the canonical identification
$$
H(V,{\mathbf Z})\cong \pi^{*}H(V,{\mathbf Z}) = \{ [x]\in H(\widetilde{V},{\mathbf Z});\,\langle[E],x\rangle=0 \}.
$$
Since $[E]$ is a norm $(-1)$-vector of $H^{2}(\widetilde{V}, {\mathbf Z})$, this implies that $\mathcal{KM}(V)$ is identified with a component of 
the Heegner divisor of norm $(-1)$-vectors of $\mathcal{KM}(\widetilde{V})$. 
Since $O(H(\widetilde{V}, {\mathbf Z}))$ acts transitively on the norm $(-1)$-vectors of $H(\widetilde{V}, {\mathbf Z})$ except
the case $\deg \widetilde{V} = 7$, i.e., $H(\widetilde{V}, {\mathbf Z}) \cong {\mathbb U}^{\oplus2} \oplus \langle -1 \rangle$, 
$\mathcal{KM}(V)$ coincides with the Heegner divisor of norm $(-1)$-vectors of $\mathcal{KM}(\widetilde{V})$ when $\deg V \not=7$. 
When $\deg V=7$, the Heegner divisor of norm $(-1)$-vectors of $\mathcal{KM}(\widetilde{V})$ consists of two components;
one is given by $\mathcal{KM}(\varSigma_{0})$ and the other is given by $\mathcal{KM}(\varSigma_{1})$, where 
$\varSigma_{n} = {\mathbf P}( {\mathcal O}_{{\mathbf P}^{1}} \oplus {\mathcal O}_{{\mathbf P}^{1}}(n) )$ is the Hirzebruch surface.
In the following theorem, we use the convention that {\em a Del Pezzo surface of degree $0$ is an Enriques surface}.

\begin{theorem}
$\Phi_{V}$ is the quasi-pullback of $\Phi_{\widetilde{V}}$ to $\mathcal{KM}(V) = [E]^{\perp}$, up to a constant. 
Namely, in the infinite product expression in Definition~\ref{def:Borcherds:product}, the following equality holds
$$
\Phi_{V} = {\rm Const}.\, \left. \frac{\Phi_{\widetilde{V}}(\cdot)}{\langle \cdot, [E]\rangle} \right|_{[E]^{\perp}},
$$
where $\langle z, [E] \rangle$ is the linear form on $H^{2}(\widetilde{V}, {\mathbf C})$ defined by the norm $(-1)$-vector $[E]$.
\end{theorem}

\begin{pf}
The result is a special case of \cite[Th.\,1.1]{Ma19}. See also \cite[Example 3.17]{Ma19}.
\end{pf}

This theorem can be summarized as the following diagrams:
$$
\begin{array}{lllllllllll}
\mathcal{KM}({\rm Enr})
&\supset
&\mathcal{KM}({\rm dP}_{1})
&\supset
&\cdots
&\supset
&\mathcal{KM}({\rm dP}_{7})
&\supset
&\mathcal{KM}(\Sigma_{1})
&\supset
&\mathcal{KM}({\mathbf P}^{2})
\\
\Phi_{\rm Enr}
&
\to
&
\Phi_{{\rm dP}_{1}}
&
\to
&
\cdots
&
\to
&
\Phi_{{\rm dP}_{7}}
&
\to
&
\Phi_{\Sigma_{1}}
&
\to
&
\Phi_{{\mathbf P}^{2}}
\\
\eta_{1^{-8}2^{8}4^{-8}}
&
\to
&
\eta_{1^{-8}2^{8}4^{-8}} \theta
&
\to
&
\cdots
&
\to
&
\eta_{1^{-8}2^{8}4^{-8}} \theta^{7}
&
\to
&
\eta_{1^{-8}2^{8}4^{-8}} \theta^{8}
&
\to
&
\eta_{1^{-8}2^{8}4^{-8}} \theta^{9}
\end{array}
$$
and
$$
\begin{array}{lllllllll}
\mathcal{KM}({\rm dP}_{7})
&\supset
&\mathcal{KM}(\Sigma_{0})
\\
\Phi_{{\rm dP}_{7}}
&
\to
&
\Phi_{\Sigma_{0}}
\\
\eta_{1^{-8}2^{8}4^{-8}} \theta^{7}
&
\to
&
\eta_{1^{-8}2^{8}4^{-8}} \theta^{8}
\end{array}
$$
where the inclusion implies the embedding as the discriminant divisor, the arrow in the second line implies the quasi-pullback (up to a constant),
and the arrow in the third line describe the change of elliptic modular form for $\Gamma_{0}(4)$ corresponding to $\Phi_{V}$. 
We remark that there are no inclusions of $\mathcal{KM}({\mathbf P}^{2})$ into $\mathcal{KM}(\Sigma_{0})$.

\section
{The invariant $\tau_{k}$ and the BCOV invariant}
\par

\subsection
{The BCOV invariant of log-Enriques surfaces}
\par
In this subsection, we prove that the invariant $\tau_{k}$ is viewed as the BCOV invariant of good log-Enriques surfaces.
Recall that for a compact connected K\"ahler orbifold $(V, \gamma)$, 
the BCOV torsion $T_{\rm BCOV}(V,\gamma)$ is defined as
$$
T_{\rm BCOV}(V,\gamma) := \exp( -\sum_{p,q\geq 0} (-1)^{p+q}pq\,\zeta_{p,q}'(0)),
$$
where $\zeta_{p,q}(s)$ is the spectral zeta function of the Laplacian $\square_{p,q}$ acting on $(p,q)$-forms on $V$ in the sense of orbifolds.
As before, the analytic torsion of the trivial line bundle on $V$ is denoted by $\tau(V,\gamma)$.

\begin{lemma}
\label{lemma:BCOV:surface}
If $\dim V = 2$, then the following equality holds:
$$
T_{\rm BCOV}(V,\gamma) = \tau(V,\gamma)^{-2}.
$$
\end{lemma}

\begin{pf}
Since $\square_{p,q}$ and $\square_{2-q,2-p}$ are isospectral via the Hodge $*$-operator, we have
$\zeta_{p,q}(s) = \zeta_{2-q,2-p}(s)$.
Since $\square_{p,q}$ and $\square_{q,p}$ are isospectral via the complex conjugation, we have $\zeta_{p,q}(s)=\zeta_{q,p}(s)$.
Using these relations, we have
\begin{equation}
\label{eqn:BCOV:surface}
-\log T_{\rm BCOV}(V,\gamma) 
=
4\zeta'_{0,0}(0) - 4\zeta'_{0,1}(0) + \zeta'_{1,1}(0).
\end{equation}
Since $\zeta_{0,0}(s) - \zeta_{0,1}(s) + \zeta_{0,2}(s) =0$ and $\zeta_{1,0}(s) - \zeta_{1,1}(s) + \zeta_{1,2}(s) =0$,
we have $4\zeta'_{0,0}(0) - 4\zeta'_{0,1}(0) = -4\zeta'_{0,2}(0)$ and
$\zeta'_{1,1}(0) = \zeta'_{1,0}(0) + \zeta'_{1,2}(0) = \zeta'_{1,0}(0) + \zeta'_{0,1}(0) = 2\zeta'_{0,1}(0)$.
Substituting these into \eqref{eqn:BCOV:surface}, we get the result.
\end{pf}

Now we have the following:

\begin{theorem}
\label{thm:BCOV:logEnriques}
Let $Y$ be a good log-Enriques surface with $k$ singular points.
Let $\gamma$ be a K\"ahler form on $Y$ in the sense of orbifolds and let $\varXi\in H^{0}(Y,K_{Y}^{\otimes2})\setminus\{0\}$
be a nowhere vanishing bicanonical form on $Y$. Then 
$$
\begin{aligned}
\tau_{\rm BCOV}(Y)
&:=
T_{\rm BCOV}(Y,\gamma){\rm Vol}(Y,\gamma)^{-2}\,\|\varXi\|_{L^{1}(Y)}^{\frac{4+k}{4}}
\left\{
\prod_{{\frak p}\in{\rm Sing}(Y)}\left(\frac{\gamma^{2}/2!}{|\varXi|}\right)({\frak p})
\right\}^{-\frac{5}{16}}
\\
&\quad\times
\exp
\left[
-\frac{1}{12}\int_{Y}
\log\left(\frac{|\varXi|}{\gamma^{2}/2!}\right)\,c_{2}(Y,\gamma)
\right]
\end{aligned}
$$
is independent of the choices of $\gamma$ and $\varXi$. In fact,
$$
\tau_{\rm BCOV}(Y) = \tau_{k}(Y)^{-2} = C_{k}^{-2} \| \Phi_{V}( \overline{\varpi}(Y)) \|^{\frac{1}{2}}, 
$$
where $C_{k}$ is the same constant as in Theorem~\ref{thm:torsion:Bor:prod}.
\end{theorem}

\begin{pf}
Since $\tau_{\rm BCOV}(Y)=\tau_{k}(Y)^{-2}$ by Theorem~\ref{thm:torsion:logEnriques} and Lemma~\ref{lemma:BCOV:surface}, 
we get the first claim. The second claim follows from Theorem~\ref{thm:torsion:Bor:prod}.
\end{pf}

We call $\tau_{\rm BCOV}(Y)$ the {\em BCOV invariant of $Y$}. When $\gamma$ is Ricci-flat and $|\varXi| = \gamma^{2}/2!$, 
we have the following simple expression:
\begin{equation}
\label{eqn:BCOV:log:Enr:Ricci:flat}
\tau_{\rm BCOV}(Y) = T_{\rm BCOV}(Y,\gamma){\rm Vol}(Y,\gamma)^{\frac{k-4}{4}}.
\end{equation}
As in the case of Enriques surfaces, the BCOV invariant of good log-Enriques surfaces is expressed by the Peterssion norm
of a Borcherds product. In particular, the BCOV invariant of log-Enriques surfaces is {\em not} a birational invariant,
for the birational equivalence classes of log-Enriques surfaces consist of a single class.

\begin{problem}
For a good log-Enriques surface $Y$, there exists a log-Enriques surface $Y'$ with a unique singular point admitting a birational morphism 
$Y\to Y'$ (cf. \cite{Zhang91}). In general, the singularity of $Y'$ is worse than those of $Y$. Can one construct a holomorphic torsion 
invariant of $Y'$ using some ALE instanton instead of the Eguchi-Hanson instanton? If this is the case, compare the holomorphic torsion 
invariants between $Y$ and $Y'$. 
\end{problem}

\begin{problem}
Let $Y$ be a good log-Enriques surface. Let $p\colon\widetilde{Y}\to Y$ be a resolution such that $p^{-1}({\rm Sing}\,Y)$
is a disjoint union of smooth $(-4)$-curves. Compare the BCOV invariant of $Y$ and that of the pair $(\widetilde{Y}, p^{-1}({\rm Sing}\,Y))$ 
defined by Zhang \cite{YZhang19}.
\end{problem}

\begin{problem}
Can one construct a holomorphic torsion invariant of log-Enriques surfaces with index $\geq3$ and prove its automorphy?
\end{problem}

\subsection
{The BCOV invariant of certain Borcea-Voisin type orbifolds}
\par
Let $Y$ be a good log-Enriques surface with $k$ singular points and let $X$ be the canonical double covering of $Y$.
Then $X$ is a nodal $K3$ surface with $k$ nodes endowed with an anti-symplectic involution $\iota$ with fixed point set
${\rm Sing}\,X =\{ p_{1}, \ldots, p_{k} \}$. Let $T$ be an elliptic curve. We define 
$$
V=V_{(X,\iota,T)} := (X\times T)/(\iota\times(-1)_{T}).
$$
Then $V$ is a Calabi-Yau orbifold of dimension $3$. Let $\widetilde{V}$ be the Borcra-Voisin orbifold
$$
\widetilde{V} = \widetilde{V}_{(\widetilde{X},\theta,T)} := (\widetilde{X}\times T)/(\theta\times(-1)_{T}),
$$
where $\pi\colon \widetilde{X} \to X$ is the minimal resolution of $X$ and $\theta$ is the involution on $\widetilde{X}$ induced from $\iota$. 
As before, we set $E_{i}:=\pi^{-1}(p_{i}) \cong {\mathbf P}^{1}$.
The birational morphism from $\widetilde{V}$ to $V$ induced by $\pi$ is denoted again by $\pi$.
Then $\pi \colon \widetilde{V} \to V$ is a partial resolution such that the $k$ cyclic quotient singularities of type 
$(\frac{1}{4},\frac{1}{4},\frac{1}{2})$ of $V$ are replaced by the milder cyclic quotient singularities of type $(\frac{1}{2},\frac{1}{2},0)$.
As an application of some results in Section~8, we compare the BCOV invariants between $\widetilde{V}$ and $V$. 
\par
Let $\gamma_{X}$ (resp. $\gamma_{\widetilde{X}}$) be a Ricci-flat K\"ahler from on $X$ (resp. $\widetilde{X}$)
and let $\gamma_{T}$ be the flat K\"ahler form with ${\rm Vol}(V,\gamma_{T})=1$. 
Let $\pi_{1} \colon V \to Y=X/\iota$ and $\pi_{2} \colon V \to T/(-1)_{T}$ be the projections. Similarly, let
$\widetilde{\pi}_{1} \colon \widetilde{V} \to \widetilde{X}/\theta$ and $\widetilde{\pi}_{2} \colon \widetilde{V} \to T/(-1)_{T}$ 
be the projections. We define a Ricci-flat K\"ahler form $\gamma$ (resp. $\widetilde{\gamma}$) on $V$ (resp. $\widetilde{V}$) by
$$
\gamma := \pi_{1}^{*}\gamma_{X} + \pi_{2}^{*}\gamma_{T},
\qquad
\widetilde{\gamma} := \widetilde{\pi}_{1}^{*}\gamma_{\widetilde{X}} + \widetilde{\pi}_{2}^{*}\gamma_{T}.
$$
Since ${\rm Sing}(X\times T) = (\{ p_{1} \} \times T) \amalg\cdots\amalg (\{ p_{k} \} \times T)$, we have
$$
\begin{aligned}
{\rm Sing}\,V 
&= 
(\{ p_{1} \} \times T/(-1)_{T}) \amalg\cdots\amalg (\{ p_{k} \} \times T/(-1)_{T}) \amalg (X^{\iota} \times T[2])
\\
&=
(\{ p_{1} \} \times T/(-1)_{T}) \amalg\cdots\amalg (\{ p_{k} \} \times T/(-1)_{T}) \amalg ({\rm Sing}\,X \times T[2]),
\end{aligned}
$$
where $T[2]$ denotes the points of order $2$ of $T$. Similarly,
$$
{\rm Sing}\,\widetilde{V} = \widetilde{X}^{\theta} \times T[2] = (E_{1}\times T[2]) \amalg\cdots\amalg (E_{k}\times T[2]).
$$
Hence the $1$-dimensional strata of ${\rm Sing}\,V$ (resp. ${\rm Sing}\,\widetilde{V}$) consist of
$k$-copies of the quotient $T/(-1)_{T}$ (resp. $4$-copies of $E_{1},\ldots,E_{k}$), 
which are endowed with the flat orbifold K\"ahler form $\gamma_{T}$ (resp. K\"ahler form $\gamma_{\widetilde{X}}|_{E_{i}}$
induced from $\gamma_{\widetilde{X}}$).
\par
Recall from \cite[(6.12)]{Yoshikawa17} that the orbifold BCOV invariant of $V$ is defined by
\begin{equation}
\label{eqn:def:orb:BCOV:V}
\begin{aligned}
\tau_{\rm BCOV}^{\rm orb}(V) 
&=
T_{\rm BCOV}(V,\gamma){\rm Vol}(V,\gamma)^{-3+\frac{\chi^{\rm orb}(V)}{12}}{\rm Vol}_{L^{2}}(H^{2}(V,{\mathbf Z}),\gamma)^{-1}
\\
&\quad\times
\prod_{i=1}^{k} \tau( \{ p_{i} \}\times (T/(-1)_{T}), \gamma_{T})^{-1}{\rm Vol}(T/(-1)_{T}, \gamma_{T})^{-1}
\\
&=
T_{\rm BCOV}(V,\gamma){\rm Vol}(V,\gamma)^{-3+\frac{\chi^{\rm orb}(V)}{12}}{\rm Vol}_{L^{2}}(H^{2}(V,{\mathbf Z}),\gamma)^{-1}
2^{k}\tau( T, \gamma_{T})^{-\frac{k}{2}},
\end{aligned}
\end{equation}
where we used the facts $\tau(T/(-1)_{T},\gamma_{T}) = \tau(T,\gamma_{T})^{1/2}$ and ${\rm Vol}(T/(-1)_{T},\gamma_{T})=1/2$
to get the second equality and ${\rm Vol}_{L^{2}}(H^{2}(V,{\mathbf Z}),\gamma)$ is the covolume of the lattice 
$H^{2}(V,{\mathbf Z})_{\rm fr}:=H^{2}(V,{\mathbf Z})/{\rm Torsion}$ with respect to the $L^{2}$ metric induced by $\gamma$. 
(In what follows, for a finitely generated ${\mathbf Z}$-module $M$, we set $M_{\rm fr} := M/{\rm Tors}(M)$.)
For the definition of $\chi^{\rm orb}(V)$, see \cite[(6.2)]{Yoshikawa17}. By \cite[Prop.\,6.2]{Yoshikawa17}, $\chi^{\rm orb}(V)$ coincides
with the Euler characteristic of a crepant resolution of $V$. Similarly, we have
$$
\begin{aligned}
\tau_{\rm BCOV}^{\rm orb}(\widetilde{V})
&=
T_{\rm BCOV}(\widetilde{V},\widetilde{\gamma})
{\rm Vol}(\widetilde{V},\widetilde{\gamma})^{-3+\frac{\chi^{\rm orb}(\widetilde{V})}{12}}
{\rm Vol}_{L^{2}}(H^{2}(\widetilde{V},{\mathbf Z}),\widetilde{\gamma})^{-1}
\\
&\quad\times
\{\prod_{i=1}^{k}\tau(E_{i},\gamma_{\widetilde{X}}|_{E_{i}}){\rm Vol}(E_{i}, \gamma_{\widetilde{X}}|_{E_{i}})\}^{-4}.
\end{aligned}
$$
\par
Let $q \colon X \times T \to V$ and $\widetilde{q} \colon \widetilde{X} \times T \to \widetilde{V}$ be the quotient maps.
Let $H^{2}(X\times T, {\mathbf Z})^{+}$ (resp. $H^{2}(\widetilde{X}\times T, {\mathbf Z})^{+}$) be the invariant subspace with respect to
the $\iota\times(-1)_{T}$ (resp. $\theta\times(-1)_{T}$)-action on $X \times T$ (resp. $\widetilde{X}\times T$). 
We define $H^{2}(X, {\mathbf Z})^{+}$ and $H^{2}(\widetilde{X}, {\mathbf Z})^{+}$ in the same way.
Let $r:={\rm rk}_{\mathbf Z} H^{2}(X, {\mathbf Z})^{+}$ and $\widetilde{r} := {\rm rk}_{\mathbf Z} H^{2}(\widetilde{X}, {\mathbf Z})^{+}$.
Then $\widetilde{r}=r+k=10+k$.
The maps of cohomologies 
$$
q^{*} \colon H^{2}(V,{\mathbf Z})_{\rm fr} \to H^{2}(X\times T, {\mathbf Z})^{+}_{\rm fr} 
= H^{2}(X,{\mathbf Z})^{+}_{\rm fr} \oplus H^{2}(T, {\mathbf Z}),
$$ 
$$
\widetilde{q}^{*} \colon H^{2}(\widetilde{V},{\mathbf Z})_{\rm fr} \to H^{2}(\widetilde{X}\times T, {\mathbf Z})^{+} 
= H^{2}(\widetilde{X},{\mathbf Z})^{+} \oplus H^{2}(T, {\mathbf Z}),
$$ 
have finite cokernel. Let ${\rm disc}(H^{2}(X,{\mathbf Z})^{+}_{\rm fr})$ be the discriminant of the lattice $H^{2}(X,{\mathbf Z})^{+}_{\rm fr}$ 
with respect to the intersection pairing $\langle\cdot,\cdot\rangle$ on $H^{2}(X,{\mathbf Z})_{\rm fr}\subset H^{2}(X,{\mathbf Q})$. 
Namely, if $\{ {\mathbf e}_{1},\ldots,{\mathbf e}_{r} \}$ is a basis of $H^{2}(X,{\mathbf Z})_{\rm fr}$, then
${\rm disc}(H^{2}(X,{\mathbf Z})^{+}_{\rm fr}) := \det( \langle{\mathbf e}_{i}, {\mathbf e}_{j} \rangle )$.
Obviously, $|{\rm Coker}\, q^{*}|$, $|{\rm Coker}\, \widetilde{q}^{*}|$, ${\rm disc}(H^{2}(X,{\mathbf Z})^{+}_{\rm fr})$, 
${\rm disc}(H^{2}(\widetilde{X},{\mathbf Z})^{+})$ depend only on $k$. Recall that the constant $C(k)$ was defined in
\eqref{eqn:constant}, which is the $k$-th power of the product of the normalized analytic torsion of the Eguchi-Hanson instanton 
and the analytic torsion of ${\mathbf P}^{1}$ endowed with the Fubini-Study metric, up to a universal constant.

\begin{theorem}
\label{thm:comparison:orb:BCOV}
The following equality holds:
$$
\frac{\tau_{\rm BCOV}^{\rm orb}(V)}{\tau_{\rm BCOV}^{\rm orb}(\widetilde{V})}
=
2^{-k-4}C(k)^{8}\left(\frac{|{\rm Coker}\, q^{*}|}{|{\rm Coker}\, \widetilde{q}^{*}|}\right)^{-2}
\left( \frac{|{\rm disc}(H^{2}(X,{\mathbf Z})^{+}_{\rm fr})|}{|{\rm disc}(H^{2}(\widetilde{X},{\mathbf Z})^{+})|} \right)^{-1}.
$$
\end{theorem}

\begin{pf}
We express $T_{\rm BCOV}(V,\gamma)$ in terms of $\tau_{{\mathbf Z}_{2}}(X,\gamma_{X})(\iota)$ and $\tau(T,\gamma_{T})$.
As is easily verified, Lemmas~8.3-8.7 of \cite{Yoshikawa17} hold true for $V$ without any change. 
Since $h^{1,1}(X)=20-k$, the coefficient $21$ of $\zeta^{T,+}(s)$ in \cite[Lemma~8.8]{Yoshikawa17} should be replaced 
by $21-k$. Hence, for $V$, the equality corresponding to \cite[Eq.\,(8.28)]{Yoshikawa17} becomes 
$$
\sum_{p,q}(-1)^{p+q}pq\,\zeta_{p,q}(s) = (24-k)\zeta^{T,+}(s) + 8\{ \zeta^{X,+}(s) - \zeta^{X,-}(s)\}.
$$
As a result, we get the following equality as in the first equality of \cite[p.\,358]{Yoshikawa17}
\begin{equation}
\label{eqn:BCOV:torsion:V}
T_{\rm BCOV}(V,\gamma) = \tau_{{\mathbf Z}_{2}}(X,\gamma_{X})(\iota)^{-4}\tau(T,\gamma_{T})^{-(12-\frac{k}{2})}.
\end{equation}
By \cite[l.2-3]{Yoshikawa17}, we have
\begin{equation}
\label{eqn:BCOV:torsion:tilde:V}
T_{\rm BCOV}(\widetilde{V},\widetilde{\gamma}) 
= 
\tau_{{\mathbf Z}_{2}}(\widetilde{X},\gamma_{\widetilde{X}})(\theta)^{-4}\tau(T,\gamma_{T})^{-12}.
\end{equation}
\par
Since $\widetilde{X}^{\theta}$ consists of $k$ copies of mutually disjoint ${\mathbf P}^{1}$, we get
$\chi^{\rm orb}(V) = \chi^{\rm orb}( \widetilde{V} ) 
= \frac{1}{2}\chi(\widetilde{X}\times T) + \frac{3}{2}\chi( \widetilde{X}^{\theta} \times T[2])=12k$
by \cite[Prop.\,6.1 and (6.3)]{Yoshikawa17}. Hence
\begin{equation}
\label{eqn:volume:V}
{\rm Vol}(V,\gamma)^{-3+\frac{\chi^{\rm orb}(V)}{12}} = {\rm Vol}(V,\gamma)^{-3+k} = 2^{3-k}\,{\rm Vol}(X,\gamma_{X})^{-3+k},
\end{equation}
where we used the fact ${\rm Vol}(T,\gamma_{T})=1$ and ${\rm Vol}(V,\gamma)={\rm Vol}(X,\gamma_{X}){\rm Vol}(T,\gamma_{T})/2$.
Similarly, 
\begin{equation}
\label{eqn:volume:tilde:V}
{\rm Vol}(\widetilde{V},\widetilde{\gamma})^{-3+\frac{\chi^{\rm orb}(\widetilde{V})}{12}} 
= 
2^{3-k}\,{\rm Vol}(\widetilde{X},\gamma_{\widetilde{X}})^{-3+k}.
\end{equation}
\par
Let $\{{\mathbf f}_{1},\ldots,{\mathbf f}_{r+1}\}$ be a basis of $H^{2}(V,{\mathbf Z})_{\rm fr}$. By definition, we have
$$
{\rm Vol}_{L^{2}}(H^{2}(V,{\mathbf Z}),\gamma)=|\det( \langle{\mathbf f}_{i},{\mathbf f}_{j}\rangle_{L^{2}} )|,
$$
where $\langle\cdot,\cdot\rangle_{L^{2}}$ denotes the $L^{2}$ inner product on $H^{2}(V,{\mathbf R})$ induced from $\gamma$.
Since ${\rm Vol}_{L^{2}}(H^{2}(T,{\mathbf Z}),\gamma_{T})=1$, the same calculations as in \cite[Lemma~13.4]{FLY08} yield that
\begin{equation}
\label{eqn:covolume:V}
\begin{aligned}
{\rm Vol}_{L^{2}}(H^{2}(V,{\mathbf Z}),\gamma) 
&= 
|{\rm Coker}\,q^{*}|^{2}\,{\rm Vol}_{L^{2}}(H^{2}(X,{\mathbf Z})^{+}_{\rm fr}\oplus H^{2}(T,{\mathbf Z}),\gamma_{X}\oplus\gamma_{T})
\\
&=
|{\rm Coker}\,q^{*}|^{2}\,{\rm Vol}_{L^{2}}(H^{2}(X,{\mathbf Z})^{+}_{\rm fr},\gamma_{X})\,{\rm Vol}(X,\gamma_{X})/2
\\
&=
2^{-(r+1)}|{\rm Coker}\,q^{*}|^{2}\,|{\rm disc}(H^{2}(X,{\mathbf Z})^{+}_{\rm fr})|\,{\rm Vol}(X,\gamma_{X}).
\end{aligned}
\end{equation}
Similarly, we have
\begin{equation}
\label{eqn:covolume:tilde:V}
{\rm Vol}_{L^{2}}(H^{2}(\widetilde{V},{\mathbf Z}),\widetilde{\gamma}) 
=
2^{-(\widetilde{r}+1)}|{\rm Coker}\,\widetilde{q}^{*}|^{2}\,|{\rm disc}(H^{2}(\widetilde{X},{\mathbf Z})^{+})|\,
{\rm Vol}(\widetilde{X},\gamma_{\widetilde{X}}).
\end{equation}
\par
Substituting \eqref{eqn:BCOV:torsion:V}, \eqref{eqn:volume:V}, \eqref{eqn:covolume:V} into \eqref{eqn:def:orb:BCOV:V} and
using \eqref{eqn:4:vanishing:torsion}, we get
\begin{equation}
\label{eqn:orb:BCOV:V}
\begin{aligned}
\tau_{\rm BCOV}^{\rm orb}(V) 
&=
2^{r+4}|{\rm Coker}\,q^{*}|^{-2}\,|{\rm disc}(H^{2}(X,{\mathbf Z})^{+}_{\rm fr})|^{-1}
\\
&\quad\times
\tau_{{\mathbf Z}_{2}}(X,\gamma_{X})(\iota)^{-4}{\rm Vol}(X,\gamma_{X})^{-4+k}\tau(T,\gamma_{T})^{-12}
\\
&=
2^{r+4}|{\rm Coker}\,q^{*}|^{-2}\,|{\rm disc}(H^{2}(X,{\mathbf Z})^{+}_{\rm fr})|^{-1}
\\
&\quad\times
\tau(X,\gamma_{X})^{-4}\tau_{{\mathbf Z}_{2}}(X,\gamma_{X})(\iota)^{-4}{\rm Vol}(X,\gamma_{X})^{-4+k}\tau(T,\gamma_{T})^{-12}
\\
&=
2^{r-k}|{\rm Coker}\,q^{*}|^{-2}\,|{\rm disc}(H^{2}(X,{\mathbf Z})^{+}_{\rm fr})|^{-1}
\tau(Y,\gamma_{Y})^{-8}{\rm Vol}(Y,\gamma_{Y})^{-4+k}\tau(T,\gamma_{T})^{-12}
\\
&=
2^{r-k}|{\rm Coker}\,q^{*}|^{-2}\,|{\rm disc}(H^{2}(X,{\mathbf Z})^{+}_{\rm fr})|^{-1}C(k)^{8}
\tau_{M}(\widetilde{X},\theta)^{-4}\tau(T,\gamma_{T})^{-12},
\end{aligned}
\end{equation}
where we used the equality $\tau(Y,\gamma_{Y})^{2}=\tau(X,\gamma_{X})\tau_{{\mathbf Z}_{2}}(X,\gamma_{X})(\iota)$
to get the third equality and 
Theorem~\ref{thm:comparison:torsion:logEnriques} to get the last equality.
Similarly, substituting \eqref{eqn:BCOV:torsion:tilde:V}, \eqref{eqn:volume:tilde:V}, \eqref{eqn:covolume:tilde:V} 
into \cite[(6.12)]{Yoshikawa17}, we get
\begin{equation}
\label{eqn:orb:BCOV:tilde:V}
\begin{aligned}
\tau_{\rm BCOV}^{\rm orb}(\widetilde{V}) 
&=
\tau_{{\mathbf Z}_{2}}(\widetilde{X},\gamma_{\widetilde{X}})(\theta)^{-4}\tau(T,\gamma_{T})^{-12}
\{\prod_{i=1}^{k}\tau(E_{i},\gamma_{\widetilde{X}}|_{E_{i}}){\rm Vol}(E_{i}, \gamma_{\widetilde{X}}|_{E_{i}})\}^{-4}
\\
&\quad\times
2^{3-k}\,{\rm Vol}(\widetilde{X},\gamma_{\widetilde{X}})^{-3+k}
2^{\widetilde{r}+1}|{\rm Coker}\,\widetilde{q}^{*}|^{-2}\,|{\rm disc}(H^{2}(\widetilde{X},{\mathbf Z})^{+})|^{-1}
{\rm Vol}(\widetilde{X},\gamma_{\widetilde{X}})^{-1}
\\
&=
2^{\widetilde{r}+4-k}|{\rm Coker}\,\widetilde{q}^{*}|^{-2}\,|{\rm disc}(H^{2}(\widetilde{X},{\mathbf Z})^{+})|^{-1}\tau(T,\gamma_{T})^{-12}
\\
&\quad\times
\tau_{{\mathbf Z}_{2}}(\widetilde{X},\gamma_{\widetilde{X}})(\theta)^{-4}{\rm Vol}(\widetilde{X},\gamma_{\widetilde{X}})^{-4+k}
\{\prod_{i=1}^{k}\tau(E_{i},\gamma_{\widetilde{X}}|_{E_{i}}){\rm Vol}(E_{i}, \gamma_{\widetilde{X}}|_{E_{i}})\}^{-4}
\\
&=
2^{r+4}|{\rm Coker}\,\widetilde{q}^{*}|^{-2}\,|{\rm disc}(H^{2}(\widetilde{X},{\mathbf Z})^{+})|^{-1}
\tau_{M}(\widetilde{X},\theta)^{-4}\tau(T,\gamma_{T})^{-12}.
\end{aligned}
\end{equation}
Comparing \eqref{eqn:orb:BCOV:V} and \eqref{eqn:orb:BCOV:tilde:V}, we get the result.
\end{pf}

We define the BCOV invariant of elliptic curve $T$ as
$$
\tau_{\rm BCOV}(T)
:=
{\rm Vol}(T,\omega)^{-1}\,\tau_{\rm BCOV}(T,\omega)\,
\exp\left[-\frac{1}{12}\int_{T}
\log\left(
\frac{i\,\xi\wedge\overline{\xi}}{\omega}
\right)\,c_{1}(T,\omega)\right],
$$
where $\omega$ is an arbitrary K\"ahler from on $T$. By \cite[Th.\,8.1]{Yoshikawa17}, $\tau_{\rm BCOV}(T)$ is independent of
the choice of $\omega$ and is expressed by the Petersson norm of the Dedekind $\eta$-function. By definition, we have
$\tau_{\rm BCOV}(T) = \tau(T,\gamma_{T})^{-1}$.
By \eqref{eqn:orb:BCOV:V}, we have the following factorization of the orbifold BCOV invariant of $V$.

\begin{corollary}
The following equality of BCOV invariants holds:
$$
\tau_{\rm BCOV}^{\rm orb}(V) 
= 
2^{r-k} |{\rm Coker}\,q^{*}|^{-2} |{\rm disc}( H^{2}(X,{\mathbf Z})^{+}_{\rm fr})|^{-1}
\tau_{\rm BCOV}(Y)^{4}\tau_{\rm BCOV}(T)^{12}.
$$
\end{corollary}

\begin{pf}
The result follows from \eqref{eqn:BCOV:log:Enr:Ricci:flat} and the third equality of \eqref{eqn:orb:BCOV:V}.
\end{pf}

\begin{remark}
In \cite[p.357 l.7]{Yoshikawa17}, it seems that the equality $H^{2}(X,{\mathbf Z}) = H^{2}(S\times T,{\mathbf Z})^{+}$ does not hold
in general. As the difference of these two quantities, $|{\rm Coker}\,\widetilde{q}^{*}|$ should appear in the formula for
$\tau_{\rm BCOV}^{\rm orb}(\widetilde{V})$ as in \eqref{eqn:orb:BCOV:tilde:V}.
\end{remark}

\begin{remark}
In this subsection, for the sake of simplicity of notation, we adopt the definitions ${\rm Vol}(V,\gamma)=\int_{V}\gamma^{3}/3!$ and 
$\langle \alpha, \beta \rangle_{L^{2}} = \int_{V} ({\mathcal H}\alpha) \wedge \overline{*}({\mathcal H}\beta)$ etc., where
${\mathcal H}(\cdot)$ denotes the harmonic projection.
If we follow the tradition in Arakelov geometry, it is more natural to define the $L^{2}$-inner product by
${\rm Vol}(V,\gamma)=(2\pi)^{-\dim V}\int_{V}\gamma^{3}/3!$ and 
$\langle \alpha, \beta \rangle_{L^{2}} = (2\pi)^{-\dim V}\int_{V} ({\mathcal H}\alpha) \wedge \overline{*}({\mathcal H}\beta)$ etc.
Notice that in \cite{Yoshikawa17}, this latter definition is adopted.
\end{remark}

\begin{problem}
Is the orbifold BCOV invariant \cite{Yoshikawa17} a birational invariant of Calabi-Yau orbifolds?
(To our knowledge, it is still open that the BCOV invariant of KLT Calabi-Yau varieties \cite{FuZhang20}
coincides with the orbifold BCOV invariant \cite{Yoshikawa17}.)
If the answer is affirmative, then it follows from Theorem~\ref{thm:comparison:orb:BCOV} that
the normalized analytic torsion of the Eguchi-Hanson instanton will essentially be given by 
the analytic torsion of ${\mathbf P}^{1}$ with respect to the Fubini-Study metric.
Once a comparison formula for the BCOV invariants for birational Calabi-Yau orbifolds is obtained, then one will get
a formula for the normalized analytic torsion of the Eguchi-Hanson instanton through Theorem~\ref{thm:comparison:orb:BCOV}.
\end{problem}


\end{document}